\documentclass{amsart}
\usepackage[utf8]{inputenc}
\usepackage[margin=1in]{geometry}

\title[{SL2}-like Properties and Generalizations of Markov Numbers]{
\texorpdfstring{$\SL_2$}{SL2}-like Properties of Matrices Over Noncommutative Rings and Generalizations of Markov Numbers}

\author{Zachary Greenberg}
\address{Max Planck Institute for Mathematics in the Sciences\\
Inselstr. 22\\04103 Leipzig, 
Germany }
\email{zachary.greenberg@mis.mpg.de}

\author{Dani Kaufman}
\address{University of Copenhagen \\
         Department of Mathematical Sciences \\
             2100 Copenhagen \o, Denmark }
\email{dk@math.ku.dk}

\author{Anna Wienhard}
\address{Max Planck Institute for Mathematics in the Sciences\\
Inselstr. 22\\04103 Leipzig, 
Germany }
\email{wienhard@mis.mpg.de}

\usepackage{amsmath}
\usepackage{amsthm}
\usepackage{amssymb}
\usepackage{mathtools}
\usepackage{shuffle}
\usepackage{tikz}
\usepackage{tikz-cd}
\usepackage{enumerate}
\usetikzlibrary {arrows.meta}
\usetikzlibrary{decorations.pathreplacing,decorations.markings}
\usepackage{xifthen}
\usepackage{subcaption}
\usepackage{hyperref}
\hypersetup{
	colorlinks=true,
	linkcolor=blue,
	hypertexnames=false,
	filecolor=black,      
	urlcolor=blue,
	citecolor=blue
}
\usepackage[nameinlink]{cleveref}

\newcommand{\st}{\mid}
\newcommand{\keyword}[1]{\textbf{\emph{#1}}}
\newtheorem{theorem}{Theorem}[section]
\newtheorem{lemma}[theorem]{Lemma} % Same counter as theorems

\newtheorem{proposition}[theorem]{Proposition}
\newtheorem{corollary}[theorem]{Corollary}
\newtheorem{conjecture}[theorem]{Conjecture}
\theoremstyle{definition}
\newtheorem{remark}[theorem]{Remark}
\newtheorem{definition}[theorem]{Definition}
\newtheorem{example}[theorem]{Example}

\tikzstyle{base}=[circle,draw,fill=black,inner sep=0pt, minimum width = 4pt]
\tikzstyle{frozen}=[rectangle,draw,fill=blue,inner sep=0pt, minimum width = 4pt,minimum height=4pt]
\tikzstyle{mutableBig}=[circle,draw=black,fill=black,inner sep=0pt,minimum size=10]
\tikzstyle{frozenBig}=[circle,draw=blue,fill=blue,inner sep=0pt,minimum size=10]
\tikzstyle{mutable}=[circle,draw=black,fill=black,inner sep=0pt,minimum size=6]
\tikzstyle{affine}=[circle,draw,fill=red,inner sep=0pt, minimum width = 4pt]
\tikzstyle{affine2}=[circle, draw, fill=red,inner sep=0pt, minimum width=8pt]
\tikzstyle{affine3}=[circle, draw, fill=red,inner sep=0pt, minimum width=12pt]
\tikzstyle{affine4}=[circle, draw, fill=red,inner sep=0pt, minimum width=16pt]
\tikzstyle{invis}=[circle,inner sep=0pt, minimum width=4pt]
\tikzstyle{fat2}=[circle,draw,fill=black,inner sep=0pt, minimum width=8pt]
\tikzstyle{fat3}=[circle,draw,fill=black,inner sep=0pt, minimum width=12pt]
\tikzstyle{fat4}=[circle,draw,fill=black,inner sep=0pt, minimum width=16pt]

\tikzstyle{vertex}=[circle,draw=black,fill=black,inner sep=0pt,minimum size=0]
\tikzstyle{line}=[thick]

\definecolor{vc0}{RGB}{0,0,0}
\definecolor{vf}{RGB}{209,242,242}
\definecolor{transpose}{named}{purple}
\definecolor{normal}{named}{black}

\newcommand{\Z}{\mathbb{Z}}
\newcommand{\Q}{\mathbb{Q}}
\newcommand{\R}{\mathbb{R}}
\newcommand{\K}{k}
\newcommand{\C}{\mathbb{C}}
\newcommand{\Proj}{\mathbb{P}}

\DeclareMathOperator{\trace}{tr}
\newcommand{\traceTilde}{\trace_\sigma}

\newcommand{\sigmaT}[1]{\sigma(#1)}
\newcommand{\Adj}{\mathrm{Adj}}
\newcommand{\tensor}{\otimes}
\renewcommand{\vec}[1]{\mathbf{#1}}

\newcommand{\id}{\mathrm{Id}}
\newcommand{\PSL}{\mathrm{PSL}}
\newcommand{\GL}{\mathrm{GL}}
\newcommand{\SL}{\mathrm{SL}}
\newcommand{\PSP}{\mathrm{PSp}}
\newcommand{\SP}{\mathrm{Sp}}
\newcommand{\SPA}{\SP_2(\mathcal{A},\sigma)}
\newcommand{\PSPA}{\PSP_2(\mathcal{A},\sigma)}

\newcommand{\matSP}[1]{\textsc{\MakeLowercase{#1}}}
\newcommand{\matSPgreek}[1]{#1}
\newcommand{\alg}[1]{\mathcal{#1}}
\newcommand{\A}{\alg{A}}
\newcommand{\B}{\alg{B}}

\newcommand{\CoordRing}{\mathcal{O}}
\newcommand{\CoordSPA}{\CoordRing(\mathrm{SP}^\sigma(2))}

\DeclareMathOperator{\Mat}{Mat}
\DeclareMathOperator{\End}{End}

\newcommand{\Diag}[2]{\mathrm{Diag}^#1(#2)}

\newcommand{\boundary}{\partial}

\newcommand{\twistedCommutatorLeft}[2]{\prescript{\sigma}{}{[#1,#2]}}
\newcommand{\twistedCommutatorRight}[2]{[#1,#2]^\sigma}

\newcommand{\halfInvariant}{\Phi}
\renewcommand{\hat}{\widehat}

\thanks{This project was partially funded by the Deutsche Forschungsgemeinschaft (DFG, German Research Foundation) Project-ID 281071066 – TRR 191.
D.K. was supported by the Danish National Research Foundation (CPH-GEOTOP-DNRF151) and the Alexander von Humboldt Foundation. 
Z.G. and A.W. were supported by the European Research Council under ERC-Advanced Grant 101018839. A.W. thanks the Hector Fellow Academy for support. \\
\newline
2020 {\em Mathematics Classification.} Primary 15B33, Secondary 16W10, 13F60, 20G42.
}

\begin{document}

\begin{abstract}
    We study $2\times 2$ matrices over noncommutative rings with anti-involution, with a special focus on the symplectic group $\SPA$. We define traces and determinants of such matrices and use them to prove a Cayley Hamilton identity and trace relations which generalize well known relations for elements of $\SL_2(R)$ over a commutative ring. We compare the structure of elements of $\SPA$ with Manin matrices over general noncommutative rings; this naturally leads to a quantization $\SPA_q$. In contrast to the usual definition of the quantum group as a deformation of the ring of matrix functions on $\SL_2(R)$, this quantization produces a group of matrices over a new noncommutative ring with involution. We finish the comparison by constructing a generalization of a Hopf algebra structure on the noncommutative ring of matrix functions of our quantum group. Finally, we use the noncommutative surface-type cluster algebras of Berenstein and Retakh to give a geometric interpretation of our Hopf algebra structure and to produce noncommutative generalizations of Markov numbers over many rings with involution including the complex numbers, dual numbers, matrix rings, and group rings.
\end{abstract}

\maketitle
\setcounter{tocdepth}{1}

\tableofcontents

\section{Introduction}
The ring of $2\times 2$ matrices over a field $\K$, $\Mat_2(\K)$ is one of the most ubiquitous groups in mathematics. The two functions, trace ($\trace$) and determinant ($\det$), defined from this ring to $\K$ are the most fundamental tools for analyzing matrices. This is showcased by the classical Cayley Hamilton theorem, for any $M\in \Mat_2(\K)$ 
$$ M^2 - \trace(M)M + \det(M)\id = 0.$$ 
%Importantly, both of these functions are conjugation invariant.

%The ring of $2\times 2$ real matrices, $\Mat_2(\R)$ is one of the most ubiquitous groups in mathematics. The two functions, trace ($\trace$) and determinant ($\det$), defined from this ring to $\R$ are the most fundamental tools for analyzing matrices. This is showcased by the classical Cayley Hamilton theorem, for any $X\in \Mat_2(\R)$ 
%$$ X^2 - \trace(X)X + \det(X)\id = 0.$$ 
%Importantly, both of these functions are conjugation invariant.

Restricting to $\SL_2(\K)$, one derives the first example of a ``trace identity'' from the Cayley Hamilton theorem by taking the trace of the whole equation:
$$ \trace(M^2) = \trace(M)^2 -2  $$ for all $M \in \SL_2$. This can be generalized and leads to the well known trace identity for matrices in $\SL_2(\K)$, going back to Fricke and Klein's work \cite{Fricke-Klein-traces}: For any $M, N \in \SL_2(\K)$ we have 
$$ \trace(MN)+ \trace(MN^{-1}) = \trace(M)\trace(N).$$ 
and for the commutator:  
\begin{equation}\label{eqn:TraceCommutator}
   \trace(MNM^{-1}N^{-1}) = \trace(M)^2+\trace(N)^2+\trace(MN)^2 - \trace(M)\trace(N)\trace(MN) - 2 
\end{equation}

These trace identities have played a major role in studying character varieties, i.e.\ homomorphisms $\pi_1(S) \to \SL_2(\K)$ up to conjugation. Since the trace is a conjugation invariant function, when $S$ is a cylinder, $\pi_1(S) = \mathbb{Z}$, any representation $\pi_1(S) \to \SL(2,\K)$ is essentially determined by the trace of a generator. For general surfaces, the traces of generators and some products of generators of the fundamental group give explicit coordinates for the character variety \cite{Procesi-traces}.

In the case of a punctured torus, $\pi_1(S) = \mathbb{Z}*\mathbb{Z}$, these representations are determined by the traces of the images of the two generators and their product.  The trace of their commutator is then determined by the trace identity for the commutator.

When considering representations $\pi_1(S) \to \SL(2,\R)$ that arise from hyperbolic structure on $S$, it is natural to consider the case when the 
trace of the commutator is $-2$ the identity. This corresponds to complete hyperbolic structures of finite volume. In that case \Cref{eqn:TraceCommutator} becomes 
\begin{equation}\label{eqn:TraceMarkov}
    x^2 + y^2+z^2 - xyz = 0. 
\end{equation}
This relates the study of hyperbolic structures on a punctured torus to the classical theory of Markov triples \cite{Markoff-Markoff}, which are of interest in number theory. Markov triples are integer solutions to 
\begin{equation}
    x^2+y^2+z^2- 3xyz =0
\end{equation}
A number appearing in a Markov triple is a Markov number. Markov numbers are well studied in number theory \cite{Aigner-Markov100Years,GS-IntegralPointsMarkov}. They are closely related to Lagrange numbers, which relate to bounds on approximating irrational numbers by rationals \cite{cusick-LagrangeMarkov,BGS-MarkovTriplesAndStrongApproximation}. 

Many attempts have been made to extend this classical story to $2\times2$ matrices over a noncommutative ring. Motivated by the theory of quantum groups, Manin studied matrices whose entries generally do not commute, but with certain commutativity conditions. His study naturally leads to a study of the quantum group via the language of Hopf algebras \cite{Manin-QuantumGroupsNoncomGeometry}. The naïve attempt to define the quantum group as an actual group fails; instead, one proceeds by quantizing the ring of functions $\mathcal{O}(\SL_2)$ by deforming its Hopf algebra structure according to Manin's commutativity relations.

In this paper we consider a different noncommutative generalization. Specifically, we consider matrices over a noncommutative ring $\A$ equipped with an anti-involution $\sigma \colon \A \rightarrow \A, \sigma(AB) = \sigma(B)\sigma(A)$. Our motivating example of such a ring is the ring $\Mat_n(\K)$ of $n\times n$ matrices over $\K$ with $\sigma$ being the transpose, but there are many other natural examples, including group rings. We think of such a ring as being `$\sigma$' commutative; When $\sigma$ is the identity, $\A$ is forced to be commutative. However, even in the commutative case one can take $\sigma$ to be a non-trivial involution, for example the complex conjugation on $\C$.  

% We will see that by replacing commutativity with the appropriate notion of $\sigma$ commutativity allows us to extend many results of linear algebra to matrices over involutive rings.

In \cite{ABRRW-SymplecticGroups}
the authors introduced the group $\SPA$ as a noncommutative generalization of $\SL(2,\R)$. 
This group is a subgroup of the ring of $2\times2$ matrices over $\A$ consisting of elements that preserve a specific symplectic form.  When $\sigma$ is the identity map we find that $\SPA \simeq \SL_2(\A)$, and in the case that $\A = \Mat_n(\K)$ we find that $\SPA \simeq \SP_{2n}(\K)$. Many of the properties of $\SL_2(\R)$ generalize to $\SPA$. In \cite{ABRRW-SymplecticGroups} the focus has been on exploring the $\SPA$-homogeneous space, in particular generalized flag varieties as well as symmetric spaces. 

Our focus begins by exploring  generalizations of the above mentioned trace identities. 
We introduce appropriate generalizations of the determinant and trace functions, taking values in $\A$, and find an analogue of the Cayley Hamilton theorem. From this we deduce generalization of the trace relations that follow from it. 
Note that the trace and determinant function we introduce are no longer conjugation invariant. However the trace function generalizes the stronger commutativity property of the classical trace function. 
For all $M,N \in \SL_2(\K)$ the trace satisfies $\trace(MN) = \trace(NM)$, this implies the conjugacy invariance. 
Here, extending the anti-involution $\sigma$ to a map on $2\times 2$-matrices over $\A$, the $\A$-valued trace function  we introduce satisfies a $\sigma$-commutative version: 
For all $\matSP{M},\matSP{N} \in \SPA$ the $\A$-valued trace function satisfies 
$\trace(\matSP{M}\matSP{N}) = \sigma(\trace(\sigma(\matSP{N}) \sigma(\matSP{M})))$.

Here the extension of $\sigma$ to matrices is by applying $\sigma$ pointwise to the entries of $\matSP{M}$. The terms $\sigma(\matSP{M})$ for $\matSP{M}\in \SPA$ are generally not elements of the symplectic group, but rather elements of $\SP_2(\A^{op},\sigma)$ when $\sigma$ is non-trivial. 

The entries of matrices in $\SPA$ satisfy analogous relations to Manin's commutation relations.  We show that there are characterizations of the group $\SPA$ which mirror the characterizations of Manin matrices. Surprisingly the set of matricies satisfying the analogous ``quantum symplectic relations'' is an honest group, contrasting with usual dilemma that the quantum group is not actually a group. This new group is essentially the same as $\SPA$, just over a different ring with a slightly modified symplectic form.  

The elements of the group $\SL_2(\K)$ are the $\K$ points of an algebraic group, $\SL(2)$. Its ring of functions, $\CoordRing(\SL(2))$, has a further structure of a Hopf algebra. The coalgebra structure and antipode of come from the group structure. We generalize this by  constructing an algebra whose $\A$ points give the group $\SPA$. This algebra is itself a noncommutative algebra with an involution. To construct and study it, we briefly introduce a basic theory of involutive algebras. We define a generalization of the quantized Hopf algebra of matrix functions, $\CoordSPA_q$. The usual quantized bialgebra of matrix functions $\CoordRing(\SL(2))_q$ is recovered as the $\sigma$ fixed points of $\CoordSPA_q$.  The only piece of the Hopf algebra structure which is not recovered directly is the antipode.

The introduction of $\SPA$ in \cite{ABRRW-SymplecticGroups} was motivated by the construction of noncommutative coordinates on symplectic character varieties in \cite{AGRW-SymplecticClusters}. 
These coordinates are generalizations of Fock-Goncharov's cluster coordinates on higher Teichmüller spaces, and provide a noncommutative generalization of Thurston's and Penner coordinates on classical (decorated) Teichm\"uller space. 
The important underlying structure is a cluster atlas on the classical moduli space of decorated representations, which is given by the usual cluster algebra associated to a surface. The positive points of this cluster of this cluster structure are then precisely the points in Penner's decorated Teichmüller space. 
%Another viewpoint on $\SL_2(\R)$ comes from the work of  \cite{AGRW-SymplecticClusters} on  theta positive representation spaces. These representation spaces are generalizations of Fock and Goncharov higher Teichmüller spaces. Fock and Goncharov associate to an oriented surface $S$ with punctures and a split real lie group $G$ a moduli space $\A_{S,G}$ whose positive points parameterize discrete and faithful representations $\pi_1(S) \to G$. When $G = \SL_2(\R)$ the positive points of this moduli space parameterize Penner's decorated Teichmüller space which in turn parameterizes hyperbolic structures on $S$ along with a choice of horocycles around each punctures. Importantly there is a cluster atlas on this moduli space which is given by the usual cluster algebra associated to a surface. 
 
By replacing $\SL_2(\R)$ with $\SPA$ and appropriately generalizing the notion of positivity, it is shown in  \cite{AGRW-SymplecticClusters} and later via spectral networks in \cite{KR-FramedLocalSystems}  that the positive points now parameterize maximal representations, and that the noncommutative coordinates give rise to a cluster atlas which gives a geometric realization of the noncommutative the noncommutative surfaces of \cite{BR-noncommutativeClusters}.  

Here we show that the noncommutative cluster structure in a surface $S$, encode some of the structures we introduce. We focus on two examples, the case when $S$  is a square and $S$ is the once punctured torus. 

When $S$ is a square, we find that the noncommutative cluster algebra is closely related to our Hopf algebra of matrix functions and gives many of its operations a geometric interpretation.  For example, the coproduct is given by gluing two squares together and performing 3 mutations.

When $S$ is a punctured torus, we recover a noncommutative generalization of the Markov numbers. For general $(\A,\sigma)$ we describe examples of integer triples $(A,B,C)$ in $\A$ that satisfy a generalized Markov equation for some constant $F_0 \in \A$:
\begin{align*}
    F_0 =& \hphantom{+}\sigma(A^{-1})B\sigma(C^{-1}) + \sigma(B^{-1})C\sigma(A^{-1})+\sigma(C^{-1})A\sigma(B^{-1})\\
    &+\sigma(C^{-1})B\sigma(A^{-1})+\sigma(A^{-1})C\sigma(B^{-1})+\sigma(B^{-1})A\sigma(C^{-1}) 
\end{align*}
When $\sigma =\id$ and $F_0 = 6$ we recover the usual Markov equation.

The cluster algebra structures gives rise to an explicit noncommutative mutation rule, which generates a new Markov triple from a given one. 

In some examples the noncommutative cluster mutation rule splits into a standard Markov mutation and a new ``shadow'' part. For example for  $\A = \C$ with $\sigma(a+bi) = a-bi$ we see:
\[ (Ae^{2\pi i a}, Be^{2\pi i b}, Ce^{2\pi i c }) \mapsto \left(\frac{B^2+C^2}{A^2}e^{2\pi i (b-c)}, Ce^{2\pi i c}, Be^{2\pi i b} \right)\]
Here the modulus of the complex number changes like a real Markov number while the argument follows a new rule. 
The set of Markov numbers has been intensely studied \cite{B-MarkovTriples,FG-DualTeichmuller,rabideau-continuedFractionsMarkov}. It would be very interesting to investigate the new shadows of complex Markov numbers. Experimentally we find when the argument is rational it achieves all possible values throughout the Markov tree. We obtain similar results for Markov numbers over the dual numbers $\{a + b \epsilon | \epsilon^2=0~ a,b\in \R\} $ relating to the work of \cite{Bonin-DualMarkovNumbers}.

Specializing the noncommutative generalization of the Markov equation to algebras like the algebra of matrices $\A = \Mat_n(\R)$  we find even more new phenomena. Instead of shadows we produce families of matrix Markov triples ``deforming'' the classic triples. The traces and determinants of these matrices are integer polynomials $T(t),D(t)$ that specialize at $t=0$ to simple functions of real Markov numbers. Over group rings we find examples of Markov triples where the starting triple consists of three nontrivially related elements. This allows to break the symmetry of the Markov tree fully. 

\subsection{Outline of the Paper}
After recalling the relevant definitions of $\SPA$ in \Cref{sec:SP2Asigma}, we define appropriate generalizations of the determinant and trace functions in \Cref{sec:invariants}. In \Cref{sec:CayleyHamilton} we establish several identities for the trace functions following from a generalization of the Cayley Hamilton theorem. %Then in \Cref{sec:quotient} we describe the twisted Abelianization of an involutive ring, a natural image of the trace function.
In \Cref{sec:Quantization} we discuss a quantization of $\SPA$ analogous to the quantization of $\SL_2(\R)$ and in \Cref{sec:NoncomMatrixFunctions} establish an analogous Hopf algebra structure. In \Cref{sec:Cluster Structure} we recall the noncommutative cluster structure on a surface. We relate our Hopf algebra of matrix functions to the cluster structure on a square. Finally in \Cref{sec:NoncomMarkov} we discuss nontrivial generalizations of Markov numbers using the noncommutative cluster structure on a punctured torus. Using this we define Markov like sequences in  matrix rings, complex numbers, dual numbers, and group rings of finite cyclic groups.

\section{Symplectic Groups over Noncommutative Rings}\label{sec:SP2Asigma}
In \cite{ABRRW-SymplecticGroups} the authors introduced symplectic groups over noncommutative rings with anti-involution. In this section we recall the definitions, prove some basic properties and introduce some new invariants. One key difference with our treatment is that we do not assume that our rings are finite dimensional algebras over a field; this leads to the notion of the left and right symplectic monoids whose intersection is the symplectic group. 

\subsection{Involutive Rings}\label{sec:rings}
An involutive ring $(\A, \sigma)$ is a pair consisting of a unital ring $\A$ with an anti-involution $\sigma \colon \A \rightarrow \A$, i.e. $ \sigma(AB) = \sigma(B) \sigma(A)$, and $\sigma(\sigma(A)) = A$ for all $A,B \in \A$. We denote by $\A^\sigma$ the set of fixed points of $\sigma$. 
\begin{definition}\label{def:SigmaCommute}
    We say two elements of $A,B\in\A$ \keyword{left $\sigma-$commute} if $\sigmaT{A}B = \sigmaT{B}A$ and  \keyword{right $\sigma-$commute} if $A~\sigmaT{B} = B~\sigmaT{A}$. To measure $\sigma-$commutativity we define the corresponding left and right $\sigma-$twisted commutators:
    \[ \twistedCommutatorLeft{A}{B} = \sigmaT{A}B-\sigmaT{B}A \hspace{3pc} \twistedCommutatorRight{A}{B} = A\sigmaT{B}-B\sigmaT{A}.\]
\end{definition}

\begin{example}
    Every commutative ring $R$ is involutive with the identity map chosen as $\sigma$. We can also take more interesting involutions on commutative rings. For example on the complex numbers we can take $\sigma$ to be complex conjugation. In this case the fixed points set $\C^\sigma$ is $\R$. \\
    A similar example arises for the dual numbers $\A=\{a + b \epsilon\st a,b\in \R,~ \epsilon^2=0\}$. 
    We can take $\sigma = \id$ and have $\A^\sigma = \A$. 
    But if we consider the conjugation involution 
    $\sigma(a + b\epsilon) = a - b\epsilon$ we again get $\A^\sigma \simeq \R$.
\end{example}
\begin{example}
    Our key example of an involutive ring is $\Mat_n(k)$ the ring of $n\times n$ matrices over a field $k$ with transpose as the anti-involution. Here the set of fixed points is the set of symmetric matrices.
\end{example}
\begin{example}
A further interesting family  of examples arises when $R$ is a commutative ring,  $G$ a group, and $\A = R[G]$ the associated group ring. 
There is a natural anti-involution $\sigma\colon  R[G]\rightarrow R[G]$ which is induced by $g \mapsto g^{-1}$. 
%    \[\sigma \left(\sum_g a_g g \cdot \sum_h b_h h\right) = \sum_{g,h} a_gb_h \sigma(gh) =\sum_{g,h}a_gb_h h^{-1}g^{-1} =  \sum_{h}b_h h^{-1}\cdot \sum_g a_g g^{-1} = \sigma\left(\sum_h b_h h\right)\sigma\left(\sum_g a_g g\right) \]\textit{}
This allows us to realize a rich class of examples, for example the Laurent polynomials over a commutative ring can be realized as group ring $R[\Z]$, and the  polynomials modulo $x^n -1$ as group ring $R[\Z/n\Z]$. 
\end{example}

We denote by $\Mat_2(\A,\sigma)$ the ring of $2\times 2$ matrices over $\A$. We extend $\sigma$ to a map  
$\sigma \colon \Mat_2(\A,\sigma) \rightarrow \Mat_2(\A,\sigma)$, by setting 
$$\sigma\left(\begin{bmatrix}A & B \\ C & D\end{bmatrix}\right) = \begin{bmatrix}\sigma(A) & \sigma(B) \\ \sigma(C) & \sigma(D)\end{bmatrix}.$$ 
Note that the map $\sigma$ is of course an involution, but does not respect the ring structure on $\Mat_2(\A,\sigma)$ in any way. 
In particular, for $\matSP{M},\matSP{N} \in \Mat_2(\A,\sigma) $, $\sigma(\matSP{M}\matSP{N})$ is not equal to $\sigma(\matSP{N})\sigma(\matSP{M})$.
We further denote by $\matSP{M}^T$ the matrix transpose, i.e.\ \[\begin{bmatrix}A & B \\ C & D\end{bmatrix}^T = \begin{bmatrix}A & C \\ B & D \end{bmatrix}.\]
Note that when $\A$ is noncommutative, matrix transpose also doesn't respect the multiplicative structure on $\Mat_2(\A,\sigma)$. Only the composite transform $\matSP{M} \mapsto \sigmaT{\matSP{M}}^T$ is an anti-homomorphism.\\

We consider $\A^2 = \A\times \A$ as an $\A$ module. Given $$\matSPgreek{\Omega}_1 = \begin{bmatrix}0 & 1 \\ -1 & 0\end{bmatrix} \in \Mat_2(\A,\sigma),$$
we can define an $\A$ valued skew symmetric right-bilinear form $\omega\colon \A^2\times \A^2 \rightarrow \A$ by 
\[\omega(\vec{v},\vec{w}) = \sigma(\vec{v})^T\matSPgreek{\Omega}_1 \vec{w}.\] This form is right bilinear in the sense that 
\[\omega(\vec{v}B,\vec{w}A) = \sigma(B)\omega(\vec{v},\vec{w})A\] for $A,B \in \A$.

The symplectic group is given by the elements of $\Mat_2(\A)$ which preserve this form. 

\begin{definition}\label{def:symplectic}
The left-$(\A,\sigma)$ symplectic monoid $\SP_2(\A,\sigma)^l$ is the set of right $\A-$linear endomorphisms of $\A^2$ preserving $\omega$. 
\[\SP_2(\A,\sigma)^l = \{f\in \End(A^2) \st \forall \vec{v},\vec{w} \in \A^2. \omega(f(\vec{v}),f(\vec{w})) = \omega(\vec{v},\vec{w})\}\]
In matrix form this condition is equivalent to:
\[\SPA^l = \{\matSP{M} \in \Mat_2(\A) \st \sigmaT{\matSP{M}}^T\matSPgreek{\Omega}_1 \matSP{M} = \matSPgreek{\Omega}_1\}. \] 
The right-$(\A,\sigma)$ symplectic monoid $\SP_2(\A,\sigma)^r$ is defined to be those elements $\matSP{M}\in \Mat_2(\A)$ such that $\sigmaT{\matSP{M}}^T \in \SP_2(\A,\sigma)^l$. In matrix form, $\SP_2(\A,\sigma)^r = \{\matSP{M}\in \Mat_2(\A) \mid \matSP{M}\matSPgreek{\Omega}_1\sigmaT{\matSP{M}}^T = \matSPgreek{\Omega}_1\}$.
The symplectic group $\SPA$ is given by the intersection of these two monoids. 
\end{definition}

It is clear that these two sets are monoids under composition, but it is necessary to check that their intersection is a group, which is done in the next section.

First we show an example illustrating the distinction between the two symplectic monoids. 

\begin{example}
    Let $(\A,\sigma)$ be a ring where there exists an element $A\in \A$ that has a right inverse, $B$ but not a left inverse, so that $AB =1$ but $BA \neq 1$. Then the element of $\Mat_2(\A) $ given by $\begin{bmatrix}
        \sigma(A) & 0 \\
        0& B 
    \end{bmatrix}$ is in $\SPA^l$ but not $\SPA^r$.
    There are many examples of such elements in general rings. For example take $\A$ to be the ring of endomorphisms of $\prod_{n=0}^\infty \R$. This ring has anti-involution $\sigma$ given by ``infinite matrix transpose.'' Now let $A$ be the right shift operator sending $(a_0,a_1,\dots)\to (0,a_0,a_1,\dots)$. Then $\sigma(A)$ is the left shift operator sending $(a_0,a_1,\dots)\to (a_1,a_2,\dots)$. It is a simple computation to check that $\sigma(A)A = \id$ but $A\sigma(A) \neq \id$. 
\end{example}

% \begin{theorem}
%     If $\A$ satisfies that every invertable element has both a left inverse and a right inverse then $\SPA^r = \SPA^l$. 
% \end{theorem}
% \begin{proof}
%     Let $M \in \SPA^l = \begin{bmatrix}
%         A & B \\
%         C & D 
%     \end{bmatrix}$ First assume that $A$ is invertable
%     Then we can write 
%     \[ \begin{bmatrix}
%         A^{-1} & 0 \\
%         0 & \sigma(A)
%     \end{bmatrix} \begin{bmatrix}
%         1 & 0 \\
%         -CA^{-1} & 1
%     \end{bmatrix} M
%     =
%     \begin{bmatrix}
%         1 & A^{-1}B \\
%         0 & 1
%     \end{bmatrix}\]
% Thus since this is a product of elements in $\SPA^l$ we have that $\begin{bmatrix}
%         1 & A^{-1}B \\
%         0 & 1
%     \end{bmatrix} \in \SPA^l$ Therefore we have that $A^{-1}B$ is symmetric, which implies the rest of the left symplectic relations. 

%     If only one of $B,C,D$ is invertable we can use this trick, first by replacing $M$ with $M\Omega, \Omega M $ or $\Omega M \Omega$ as nessesary. 

% \end{proof}

% We can define a slight twist of this definition simply by changing $\Omega$ to a different symplectic form (i.e. satisfying $\sigmaT{\Omega}^T = -\Omega$.

We realize familiar groups for simple choices of $\A$.
\begin{example}\label{ex:MatrixRing}
Let $\A = \Mat_n(k)$ be the ring of $n\times n$ matrices over a field $k$, and $\sigma$ be the matrix transpose. In this case $\SPA \simeq \SP(2n,k)$. 
See \cite{ABRRW-SymplecticGroups} for more examples involving rings of matrices. 
\end{example}

The choice of anti-involution is crucial and dramatically affects the structure of $\SPA$ even if the underlying algebra is commutative. 
\begin{example}\label{ex:involution}
 Consider $\A = \C$, the complex numbers.
 In the case with trivial involution the group $\SPA$ is simply $\SL_2(\mathbb{C})$ but in the case that the $\sigma$ is complex conjugation the group $\SPA$ can be computed to be $\SL_2(\R)\times \text{U}(1)$
\end{example}

\subsection{Basic Properties of the Symplectic Group}

Let us look at some basic properties of the symplectic group $\SP_2(\A,\sigma) \subset \Mat_2(\A,\sigma)$. Our descriptions of the group come from analysing the inverse of an element. This is also where we find the subtleties surrounding the distinction between the two monoids and the symplectic group.

\begin{definition}
    The \keyword{adjoint matrix} of  $\matSP{M} \in \Mat_2(\A,\sigma)$ is given by \[ \Adj(\matSP{M}) =\begin{bmatrix}
    \sigma(D) & -\sigma(B) \\
    -\sigma(C) & \sigma(A)
\end{bmatrix}.\]
\end{definition}

\begin{lemma}\label{thm:SP2AInverse}
Let $\matSP{M} \in \Mat_2(\A,\sigma)$. We have the following
\begin{itemize}
    \item $\matSP{M} \in \SPA^l$ if and only if $\matSP{M}$ has a left inverse given by the adjoint matrix.
    \item $\matSP{M} \in \SPA^r$ if and only if $\matSP{M}$ has a right inverse given by the adjoint matrix.
    \item $\matSP{M} \in \SPA$ if and only if $\matSP{M}$ has a two-sided inverse given by the adjoint matrix.
\end{itemize}
Moreover, the adjoint of an element in $\SP_2(\A,\sigma)^l$ is an element of $\SP_2(\A,\sigma)^r$. 
\end{lemma}
\begin{proof}
As $\matSPgreek{\Omega}_1^{-1} = -\matSPgreek{\Omega}_1$ it is a simple computation to see 
\[\Adj(\matSP{M}) = \matSPgreek{\Omega}_1^{-1}\sigmaT{\matSP{M}}^T\matSPgreek{\Omega}_1 = \matSPgreek{\Omega}_1\sigmaT{\matSP{M}}^T\matSPgreek{\Omega}_1^{-1}.\]
Then $\Adj(\matSP{M})$ being a left inverse for $\matSP{M}$ is equivalent to $1 = \matSPgreek{\Omega}_1^{-1}\sigmaT{\matSP{M}}^T\matSPgreek{\Omega}_1 \cdot \matSP{M}$. Left multiplying by $\matSPgreek{\Omega}_1$ yields $\matSPgreek{\Omega}_1 = \sigmaT{\matSP{M}}^T\matSPgreek{\Omega}_1 \matSP{M}$ which is the definition of $\matSP{M} \in \SPA^l$. A similar calculation shows that $\matSP{M}\Adj(\matSP{M}) = 1$ is equivalent to the defining equation of $\SPA^r$.\\
Furthermore we have 
\[ \Adj(\matSP{M})\matSPgreek{\Omega}_1\sigmaT{\Adj(\matSP{M})}^{T} = (-\matSPgreek{\Omega}_1\sigma(\matSP{M})^{T}\matSPgreek{\Omega}_1)\matSPgreek{\Omega}_1 (-\matSPgreek{\Omega}_1 \matSP{M}\matSPgreek{\Omega}_1) = \matSPgreek{\Omega}_1\] 
showing that the adjoint of $\matSP{M} \in \SPA^l$ is an element of $\SPA^r$. 

% The reverse implication comes from computing $\Adj(f)f$. For this product to be the identity matrix, $f$ must  satisfy the conditions of \Cref{thm:SymplecticEquations}.
% \[
% \begin{bmatrix}
%     \sigmaT{D} & -\sigmaT{B}\\
%     -\sigmaT{C} & \sigmaT{A}
% \end{bmatrix} \begin{bmatrix}
%     A & B\\ C & D
% \end{bmatrix}
% = \begin{bmatrix}
%     \sigmaT{D}A-\sigmaT{B}C & \sigmaT{D}B-\sigmaT{B}D\\
%     \sigmaT{A}C-\sigmaT{C}A & \sigmaT{A}D-\sigmaT{C}B
% \end{bmatrix} = \begin{bmatrix}
%     1 & 0\\0 & 1
% \end{bmatrix}
% \]
\end{proof}

\begin{corollary}\label{thm:SymplecticEquations}
An element $\matSP{M} = \begin{bmatrix} A& B\\ C &D\end{bmatrix} \in \Mat_2(\A,\sigma)$ is in $\SP_2(\A,\sigma)^l$ if and only if $\matSP{M}$ satisfies the left symplectic equations:
\begin{equation}\label{eq:sp2a_quations_l}
    \begin{aligned}
        \sigma(A)C-\sigma(C)A =& 0\\
        \sigma(B)D-\sigma(D)B =& 0\\
        \sigma(A)D-\sigma(C)B =& 1
    \end{aligned}
\end{equation}

Similarly $\matSP{M}$ is in $\SP_2(\A,\sigma)^r$ if and only if $\matSP{M}$ satisfies the right symplectic equations:
    \begin{equation}\label{eq:sp2a_quations_r}
    \begin{aligned}
        A\sigma(B)-B\sigma(A) =& 0\\
        C\sigma(D)-D\sigma(C) =& 0\\
        A\sigma(D)-B\sigma(C) =& 1
    \end{aligned}
    \end{equation}
    Thus $\matSP{M}$ is in $\SP_2(\A,\sigma)$ if and only if $\matSP{M}$ satisfies both the left and right symplectic equations.
\end{corollary}
\begin{proof}
    The left symplectic equations come from computing $\Adj(\matSP{M})\matSP{M}$ and setting the result equal to the identity. To obtain the right symplectic equations we use the fact that $\sigmaT{\matSP{M}}^T \in \SP_2(\A,\sigma)^l$.
\end{proof}

With very minor assumptions on $\matSP{M}$ we can conclude that $\matSP{M} \in \SP_2(\A,\sigma)$ by checking a subset of the equations.
\begin{lemma}\label{thm:SymplecticEquationsAinv}
    Let $\matSP{M} = \begin{bmatrix} A& B\\ C &D\end{bmatrix} \in \Mat_2(\A,\sigma)$ and assume $A$ is invertible. Then $\matSP{M} \in \SPA$ if $\matSP{M}$ satisfies the following three equations:
    \begin{equation}\label{eq:sp2a_quations_A}
        \begin{aligned}
            \sigma(A)C-\sigma(C)A =& 0\\
            A\sigma(B)-B\sigma(A) =& 0\\
            \sigma(A)D-\sigma(C)B =& 1
        \end{aligned}
    \end{equation}

\end{lemma}
\begin{proof}
    We obtain the missing equations by conjugating our assumed equations by $A$. First conjugate $\sigmaT{A}D-\sigmaT{C}B = 1$ by $\sigmaT{A}$. This simplifies using the other two equations:
    \begin{align*}
         \sigma(A^{-1})\sigma(A)D\sigma(A)-\sigma(A^{-1})\sigma(C)B\sigma(A) =& \sigmaT{A^{-1}}\sigmaT{A}\\
         D\sigma(A)-C\sigma(B) =& 1
%         A\sigmaT{D}-B\sigmaT{C} =& 1
    \end{align*}
    Similarly we multiply the first equation on the right by $A^{-1}$ and on the left by $\sigmaT{D}$.
    \begin{align*}
        \sigma(A^{-1})\sigma(A)C\sigma(D) - \sigma(A^{-1})\sigma(C)A\sigma(D) = 0 \\
    C\sigma(D) - \sigma(A^{-1})\sigma(C)(1+B\sigma(C)) = 0 \\
    C\sigma(D) - \sigma(A^{-1})(1+\sigma(C)B)\sigma(C) = 0\\
    C\sigma(D) - D\sigma(C) = 0  
    \end{align*} 
    The final equation is obtained from the second equation by right multiplication by $A^{-1}$ and left multiplication by $D$.
\end{proof}
\begin{remark}
    A similar lemma holds if any of one of $B, C, D$ is invertible. 
\end{remark}

\begin{theorem}\label{thm:left=right}
    If $\A$ satisfies that every endomorphism of $\A^2$ is injective if and only if it is also surjective, then $\SPA^r=\SPA^l$
\end{theorem}
\begin{proof}
    Let $\matSP{M}\in \SPA^l$. Then $\matSP{M}$ has a left inverse, call it $\Mat{N}= \Adj(\matSP{M})$ and $\matSP{N}\matSP{M} = \id$. Since $\matSP{M}$ has a left inverse it is injective. By assumption this implies $\matSP{M}$ is surjective and thus has a right inverse $\matSP{N}'$. A standard fact is that if $\matSP{M}$ has a left and a right inverse they must agree, so $\matSP{N} = \matSP{N}'$. Since $\matSP{N} = \Adj(\matSP{M})$ is also a right inverse to $\matSP{M}$ we have that $\matSP{M} \in \SPA^r$. A similar argument proves the reverse inclusion.  
\end{proof}

\begin{remark}
    The hypothesis of Theorem~\ref{thm:left=right} is satisfied in many relevant cases. For example, if $\A$ is a finite dimensional algebra over a field. 
\end{remark}
\begin{corollary}
    When the hypothesis of \Cref{thm:left=right} is satisfied a matrix $\matSP{M}$ is in $\SPA$ if and only if $\matSP{M}$ satisfies the left symplectic equations if and only if $\matSP{M}$ satisfies the right symplectic equations.
\end{corollary} 

The element $\sigma(M)$ for $\matSP{M} \in \SPA$ is not generally back in $\SPA$. Nevertheless this element does appear in many of our identities proved in later sections. 
\begin{proposition}\label{prop:sigma(f)}
    The element $\sigma(\matSP{M})$ is an element of $\SP_2(\A^{op},\sigma)$ if and only if $\matSP{M} \in \SPA$. Moreover, the map $\sigma: \SPA \to \SP_2(\A^{op},\sigma)$ is an isomorphism of groups. 
\end{proposition}
\begin{proof}
    The first statement is clear from our characterization of elements of $\SPA$. The second proposition follows from checking that $\sigma(\matSP{M}\matSP{N}) = \sigma(\matSP{M})*\sigma(\matSP{N})$ remembering that the second multiplication happens in $\A^{op}$.
\end{proof}

% \begin{remark}
%     One might hope to prove that the monoids $\SPA^r$ and $\SPA^l$ actually coincide. This is equivalent to proving that the left symplectic equations imply the right ones, but this seems impossible. However, for many choices of $\A$ this is actually true. For example the cases considered in \cite{ABRRW-SymplecticGroups} satisfy this property, since $\A$ is assumed to be a finite dimensional algebra over a field. 
% \end{remark}

% \begin{proof}
% It suffices to see when $\sigma(f)^{T}\Omega_1 f = \Omega_1$:
% \[
% \begin{bmatrix}
%     \sigma(A) & \sigma(C)\\
%     \sigma(B) & \sigma(D)
% \end{bmatrix}
% \begin{bmatrix}0 & 1 \\ -1 & 0\end{bmatrix}
% \begin{bmatrix}
%     A & B\\
%     C & D
% \end{bmatrix}
% = \begin{bmatrix}
%     \sigma(A)C-\sigma(C)A & \sigma(A)D-\sigma(C)B\\
%     -\sigma(D)A+\sigma(B)C & \sigma(B)D-\sigma(D)B
% \end{bmatrix}
% \]
% There are only three independent equations as the top right and bottom left entries are equivalent after applying $\sigma$.
% \end{proof}

\subsection{Comparison with Manin Matrices}
Here we remark on the analogy to Manin matrices over a noncommutative ring $R$. First, we recall the definitions from \cite{Manin-QuantumGroupsRemarks}.
\begin{definition}
     A matrix $\matSP{M} = \begin{bmatrix}
        A &B\\C &D
    \end{bmatrix} \in \Mat_2(R)$ is \keyword{Manin} if it satisfies the Manin commutativity equations:
   \begin{equation}\label{eqn:ManinCommutivityEquations}
       AC-CA=0 \hspace{2pc} BD-DB=0 \hspace{2pc} AD-CB = DA-BC
   \end{equation}
\end{definition}
Note that the first two left symplectic equations correspond to the first two Manin commutativity equations by replacing products $xy$ with $\sigmaT{x}y$. If we add the additional constraint to a Manin matrix that $AD-CB = 1$ then the third equations correspond. Since $\sigmaT{1} = 1$,  applying $\sigma$ to the third symplectic equation yields the corresponding equation $\sigmaT{D}A-\sigmaT{B}C = 1 = \sigmaT{A}D-\sigmaT{B}C$.

In \cite{CFR-ManinMatrices} they show the analogous characterization of Manin matrices by their inverse.
\begin{proposition}
    A matrix $\begin{bmatrix}
        A & B\\ C & D
    \end{bmatrix}$ is Manin if and only if its left inverse is given by $\frac{1}{AD-CB}\begin{bmatrix}D & -B\\ -C & A\end{bmatrix}$
\end{proposition}

Here we find an interesting departure between our group $\SPA$ and the theory of Manin matrices. If a matrix and its transpose are both Manin, then the entries are all forced to pairwise commute. In our situation this condition would simply mean that $\matSP{M}\in \SPA$, which in many cases is implied already by $\matSP{M}\in \SPA^l$ due to \Cref{thm:left=right}.

\subsection{The \texorpdfstring{$(A,\sigma)$}{a sigma} Plane}\label{sec:AsigmaPlane}
Another characterization of matrices in $\SPA$ is via the matrix action on the $\sigma-$affine plane. This space and action is defined in \cite{ABRRW-SymplecticGroups} and we recall the definition here.
%We already saw that matrices in $\SPA$ satisfy properties very similar to Manin matrices. In particular the ``$\sigma$-commutivity'' condition of \Cref{thm:SymplecticEquations} correspond to the commutativity conditions of \Cref{eqn:ManinCommutivityEquations}. Similarly the formulas for inversion and Cayley Hamilton theorem correspond. Now we recall the action of $\SPA$ on the Affine plane to complete the correspondence.\\

\begin{definition}
    Two vectors $\vec{v}$ and $\vec{w}$ in $\A^2$ are \keyword{transverse} if $\{\vec{v},\vec{w}\}$ is a basis of $\A^2$ as an $\A$ module. We then say a  vector $\vec{v}\in \A^2$ is \keyword{regular} if there is another vector $\vec{w}\in \A^2$ transverse to it.
\end{definition}

\begin{definition}
    A vector $\vec{v} \in \A^2$ is \keyword{isotropic} if $\omega(\vec{v},\vec{v}) = 0$. The \keyword{left $\sigma-$affine plane} $\A^2_\omega$ is the set of regular isotropic vectors. 
    \[\A^2_\omega = \{\vec{v} \in \A^2 \mid \text{$v$ regular and } \omega(\vec{v},\vec{v}) = 0 \}\] 
\end{definition}

\begin{lemma}
    $\vec{v} = (V_1,V_2) \in \A^2$ is isotropic if and only if $\sigma(V_1)V_2 = \sigma(V_2)V_1$
\end{lemma}
\begin{proof}
This is a simple computation given by expanding $\omega(\vec{v},\vec{v})$.
    \begin{align*}
        0= \omega(\vec{v},\vec{v}) =& \begin{bmatrix} \sigma(V_1) & \sigma(V_2)\end{bmatrix} \begin{bmatrix} 0 & 1\\ -1 & 0 \end{bmatrix} \begin{bmatrix}  V_1\\V_2 \end{bmatrix} = -\sigma(V_2)V_1+\sigma(V_1)V_2
    \end{align*}
\end{proof}

In this way the left $\sigma-$affine plane can be thought of as the space of vectors in $\A^2$ with left $\sigma-$commuting coordinates.

\begin{remark}
We stress that the positions of $\sigma$ are important when $\sigma$ is nontrivial. The right $\sigma-$affine plane
    \[\A^2_{\omega_{\sigma}} = \{\vec{v} \in \A^2 \mid \omega(\sigma (\vec{v}),\sigma(\vec{v})) =0\} = \{ \vec{v}\in \A^2 \mid V_1\sigma(V_2) = V_2\sigma(V_1) \} \]
    is not the same set as $\A^2_\omega$ if $\sigma$ is nontrivial. 
\end{remark}

There is a natural left action of $\Mat_2(\A)$ on $\A^2$ given by 
\[\begin{bmatrix} A & B\\C &D\end{bmatrix} \begin{bmatrix}V_1 \\ V_2\end{bmatrix} = \begin{bmatrix}
    AV_1 + BV_2\\CV_1+DV_2
\end{bmatrix}\] and a right action by 
\[\begin{bmatrix}V_1 & V_2\end{bmatrix} \begin{bmatrix} A & B\\C &D\end{bmatrix} = \begin{bmatrix}
    V_1A + V_2C &  V_1B+V_2D
\end{bmatrix}\]

%This action also induces an action on $\Proj(\A^2) = \{v\A \mid v\in \A^2 \text{ regular} \}$ the space of $\A$ lines in $\A^2$. In \cite{ABRRW-SymplecticGroups} they prove the following theorem:
%\begin{theorem}
%    The symplectic group $\SP_2(\A,\sigma)$ acts transitively on the space of isotropic lines $\Proj_\omega(\A^2) = \{v\A \in \Proj(\A^2) \mid \omega(v,v) = 0\}$.
%\end{theorem}
%\begin{proof}
%    This is Theorem 1.1 of \cite{ABRRW-SymplecticGroups}.
%\end{proof}

\begin{theorem}\label{thm:SP2ACharacterization_ActionOnSigmaPlane}
    A matrix $\matSP{M}$ is in $\SPA^l$ if and only if $\sigmaT{A}D-\sigmaT{C}B = 1$ and $\matSP{M}$ preserves the left $\sigma$-affine plane under the left matrix action. Similarly $\matSP{M} \in \SPA^r$ if and only if $A\sigma(D)- B\sigma(C)=1$ and $\matSP{M}$ preserves the right $\sigma-$affine plane via the right action.
\end{theorem}
\begin{proof}
    The forward implication follows from Theorem 1.1 of \cite{ABRRW-SymplecticGroups}. They show that $\SPA$ preserves the space of regular isotropic lines $\Proj_\omega(\A^2) = \{\vec{v}\A \mid \vec{v} \in \A^2 \text{ regular and } \omega(\vec{v},\vec{v}) = 0\}$. \\
    We now show the reverse implication. Assume $\vec{v} = (V_1,V_2) \in \A^2_\omega$. We compute the action of $\matSP{M}$ on $\vec{v}$ and assume the result is again in $\A^2_\omega$. This implies:
    \begin{align*}
        0 =& \sigmaT{AV_1+BV_2}(CV_1+DV_2)-\sigmaT{CV_1+DV_2}(AV_1+BV_2)\\
        =& \sigmaT{V_1}(\sigmaT{A}C-\sigma(C)A)V_1 + \sigma(V_2)(\sigma(B)D-\sigma(D)B)V_2 \\
        & + \sigma(V_1)(\sigma(A)D-\sigma(C)B)V_2 + \sigma(V_2)(\sigma(B)C-\sigma(D)A)V_1\\
        =& \sigmaT{V_1}(\sigmaT{A}C-\sigmaT{C}A)V_1 + \sigmaT{V_2} (\sigmaT{B}D-\sigmaT{D}B)V_2 + \sigmaT{V_1}V_2-\sigmaT{V_2}V_1\\
        =&  \sigmaT{V_1}(\sigmaT{A}C-\sigmaT{C}A)V_1 + \sigmaT{V_2} (\sigmaT{B}D-\sigmaT{D}B)V_2
    \end{align*}
    The last two simplifications are due to the assumption that $\sigmaT{A}D-\sigmaT{C}B= 1$ and that $(V_1,V_2)$ left $\sigma-$commute.\\
    To finish the proof we evaluate the previous equality at the points $(1,0)$ and $(0,1)$ each of which belong to $\A_\omega^2$. This yields:
    \[0 = \sigmaT{A}C-\sigmaT{C}A \hspace{1pc}\text{and}\hspace{1pc} 0=\sigmaT{B}D-\sigmaT{D}B \]
    Thus $\matSP{M}$ satisfies the left symplectic equations and so by \Cref{thm:SymplecticEquations} $\matSP{M}\in \SPA^l$.
    The same proof for the right action follows by applying $\sigmaT{}^T$ everywhere.
\end{proof}

We now compare this result to the case of Manin matrices. Manin matrices don't preserve the plane of commuting elements $\{\vec{v}= (v_1,v_2) \in R^2 \mid v_1v_2 = v_2v_1 \}$. However they can be characterized by their action on the plane of points whose entries not only commute with each other, but also commute with the entries of $\matSP{M}$. 
\begin{proposition}[\cite{CFR-ManinMatrices,Manin-QuantumGroupsNoncomGeometry}]
    The matrix $\begin{bmatrix} A & B\\C & D\end{bmatrix}$ is Manin if and only if for all $v_1,v_2\in R$ such that $v_1$ and $v_2$ commute with each other and $A,B,C,D$ the pair $\widetilde{v}_1 = Av_1+Bv_2$ and $\widetilde{v}_2 = Cv_1+Dv_2$ commute. 
\end{proposition}

%%Not a group:
% For this reason we define $\omega_\sigma \colon \A^2\times \A^2 \rightarrow \A$ by 
% \[\omega_\sigma(v,w) = \omega(\sigma(v),\sigma(w)) = v^T\Omega_I\sigma(w).\]
% We write $\SPsigma(\A,\sigma)$ for the group preserving $\omega_\sigma$.\\
% \begin{lemma}
%     An element $f = \begin{bmatrix}
%         A' & B'\\
%         C' & D'
%     \end{bmatrix} \in \SPsigma(\A,\sigma)$ if and only if
%     \begin{align*}
%         A'\sigma(C') - C'\sigma(A') =& 0\\
%         B'\sigma(D') - D'\sigma(B') =& 0\\
%         A'\sigma(D') - C'\sigma(B') =& 1
%     \end{align*}
% \end{lemma}
% \begin{proof}
%     This is an identical computation to the analogous statement for the standard Symplectic group (\Cref{thm:SymplecticEquations})
% \end{proof}
% The action of the involution $\sigma$ on $\Mat_2(\A,\sigma)$ exchanges the role of $\omega$ and $\omega_\sigma$. Thus we have the following lemma:
% \begin{lemma}
%     An element $f$ is in $\SP_2(\A,\sigma)$ if and only if  $\sigma(f)$ is in $\SPsigma_2(\A,\sigma)$. Moreover $\SPsigma_2(\A,\sigma)$ acts (nontransitively) on $\A^2_{\omega_\sigma}$.
% \end{lemma}

%\begin{definition}
%    For $f \in \Mat_2(\A,\sigma)$ the transpose $f^{T}$ is given  by
%    \[\begin{bmatrix}
%        \sigmaT{A} & \sigmaT{C}\\
%        \sigmaT{B} & \sigmaT{D}
%    \end{bmatrix} \]
%\end{definition}
%\begin{lemma}
%    If $f \in \SPA$ then so is $f^T$.
%\end{lemma}
%\begin{proof}
%    This is a simple calculation from conditions in %\Cref{thm:SymplecticEquations}.
%\end{proof}

\section{Invariants of matrices}\label{sec:invariants}
In this section we introduce general functions on $\Mat_2(\A,\sigma)$, which we call traces and determinants. When $\A$ is commutative and $\sigma =\id$ these are precisely the trace and determinants of $2 \times 2$ matrices. However, when $\sigma$ is not the identity and $\A$ is potentially noncommutative they have quite different properties.

\begin{definition}
    The \keyword{left determinant} $\det_l \colon \Mat_2(\A,\sigma) \rightarrow \A$ is given by
    \[\begin{bmatrix}A& B\\ C & D\end{bmatrix} \mapsto \sigma(A)D - \sigma(C)B\]
    The corresponding \keyword{right determinant} $\det_r \colon \Mat_2(\A,\sigma) \rightarrow \A$ is given by 
    \[\begin{bmatrix}A& B\\ C & D\end{bmatrix} \mapsto A \sigmaT{D} - B\sigmaT{C}\]
\end{definition}
Note that $\det_r(\matSP{M}) = \det_l(\sigmaT{\matSP{M}}^T)$. However $\det_r(\sigma(\matSP{M}))$ has no relation to $\det_r(\matSP{M})$ or $\det_l(\matSP{M})$ and in fact $\det_r(\matSP{M})$ and $\det_l(\matSP{M})$ can both be invertible while $\det_r(\sigma(\matSP{M}))$ is not.
\begin{definition}
    The \keyword{trace} $\trace\colon \Mat_2(\A,\sigma) \rightarrow\A$ is given by 
    \[ \begin{bmatrix}A& B\\ C & D\end{bmatrix} \mapsto A+D\]
    We further define the \keyword{twisted trace} $\traceTilde := \sigma(A) + D$
\end{definition}
% \begin{remark}
%     Note that $\trace_\sigma = \trace \circ \sigma$. However for the determinant the function  determinant $ \det \circ \sigma $ is a distinct  function. When $x = \begin{bmatrix}
%         A & B \\
%         C & D
%     \end{bmatrix}$, then $\det \circ \sigma(x) = A\sigma(D) - C\sigma(B)$.  In fact in general $\det(x)$ being invertible does not imply that $\det(\sigma(x))$ is invertible.  
%     However we have $\det_\sigma(x) = \det(\sigma(\Adj(x))^T).$
%    \end{remark}
Unlike the determinant, both the trace and the twisted trace commute with $\sigma$, $\trace(\sigmaT{\matSP{M}}) = \sigma(\trace(\matSP{M}))$.\\
We now establish the first property of the trace:

\begin{proposition}\label{prop:sigma_commute_trace}
For all $\matSP{M},\Mat{N} \in \Mat_2(\A,\sigma)$ we have 
\[\trace(\sigma(\matSP{M}\matSP{N})) = \trace(\sigma(\matSP{N})\cdot \sigma(\matSP{M}))\]
\end{proposition}
\begin{proof}
    This follows by direct matrix computation.
\end{proof}
\begin{remark}
    Note that in general $\sigma(\matSP{M}\matSP{N})$ is not equal to $\sigma(\matSP{N})\sigma(\matSP{M})$. If $\sigma$ is trivial and thus $\A$ is commutative, the above identity reduces to the classical commutativity property of the trace $\trace(\matSP{M}\matSP{N}) = \trace(\matSP{N}\matSP{M})$, which further implies the conjugation invariance of the trace function  $\trace(\matSP{N}\matSP{M}\matSP{N}^{-1}) = \trace(\matSP{M})$.
\end{remark}

Even though the trace is not conjugation invariant when $\sigma$ is non-trivial, we obtain generalizations of well known trace identities for $\SL_2(\R)$. The first of which is the following.  
\begin{proposition}
    If $\matSP{M} \in \SP_2(\A,\sigma)$ then $\trace(\matSP{M}^{-1}) = \sigma(\trace(\matSP{M}))$
\end{proposition}
\begin{proof}
    This is clear from the explicit formula for $\matSP{M}^{-1}$.
\end{proof}
\begin{remark}
    When $\sigma$ is trivial this reduces to the well known identity $\trace(\matSP{M}) = \trace(\matSP{M}^{-1})$ for $\matSP{M}\in \SL_2(\R)$. 
\end{remark}

We can compute how the trace  is affected by conjugation.

Let $\matSP{M} = \begin{bmatrix}
    A & B\\C & D
\end{bmatrix}$ and $\matSP{N} = \begin{bmatrix}
    E & F\\G & H
\end{bmatrix}$. Then $\matSP{N}^{-1} = \begin{bmatrix}
    \sigma(H) & -\sigma(F)\\-\sigma(G) & \sigma(E)
\end{bmatrix}$. So 
\begin{align*}
    \trace(\matSP{N}^{-1}\matSP{M}\matSP{N}) =& +\sigma(H)AE+\sigma(H)BG-\sigma(F)CE-\sigma(F)DG\\
    &-\sigma(G)AF-\sigma(G)BH+\sigma(E)CF+\sigma(E)DH
\end{align*}

% In order to get an identity we have to place some conditions on the entries of  $x$.

% \begin{lemma}
%     If $B \ in \A^\sigma$ and $A + \sigma(D) $ is in the center of $\A$ then $A +\sigma(D) \in \A^\sigma $.
% \end{lemma}
% \begin{proof}
%     $$(A+\sigma(D))B=(A\sigma(B)+\sigma(D)B) =B\sigma(A) =\sigma(B)D= B(\sigma(A)+D)$$
% \end{proof}

In particular from this we obtain the following proposition. 

\begin{proposition}\label{thm:traceConjugation}
    If $A+\sigmaT{D}$ is in the center of $\A$ and $A+\sigmaT{D},B,C \in \A^\sigma$ and $\matSP{N} \in \SP_2(\A,\sigma)$ then $\trace(\matSP{M}) + \sigma(\trace(\matSP{M})) = \trace(\matSP{N}^{-1}\matSP{M}\matSP{N}) + \sigma(\trace(\matSP{N}^{-1}\matSP{M}\matSP{N}))$
\end{proposition}
\begin{proof}
    We compute and gather terms by first and last factor
    \begin{align*}
        \trace(\matSP{N}^{-1}&\matSP{M}\matSP{N})+\sigma(\trace(\matSP{N}^{-1}\matSP{M}\matSP{N})) \\
        =&+\sigmaT{H}[A+\sigmaT{D}]E+\sigmaT{H}[B-\sigma(B)]G-\sigmaT{F}[C-\sigmaT{C}]E-\sigmaT{F}[\sigmaT{A}+D]G\\
        &-\sigmaT{G}[A+\sigmaT{D}]F-\sigmaT{G}[B-\sigmaT{B}]H+\sigmaT{E}[C-\sigmaT{C}]F+\sigmaT{E}[\sigmaT{A}+D]H\\       
 %       =& +\sigma(H)[A+\sigma(D)]E + \sigma(E)[\sigma(A)+D]H+\sigma(H)[B-\sigma(B)]G+\sigma(E)[C-\sigma(C)]F\\
%    &-\sigma(F)[\sigma(A)+D]G-\sigma(G)[A+\sigma(D)]F-\sigma(G)[B-\sigma(B)]H-\sigma(F)[C-\sigma(C)]E\\
    =&(\sigma(A)+D)(\sigma(E)H-\sigma(F)G) + (A+\sigma(D))(\sigma(H)E-\sigma(G)F)\\
    =& (A+\sigmaT{D})(\sigmaT{E}H-\sigmaT{F}G + \sigmaT{H}E-\sigmaT{G}F)\\
    =&(A+\sigmaT{D})(I + I)\\
    =& A+D + \sigma(A) + \sigma(D)
    \end{align*}
\end{proof}

We note that $\traceTilde$ also has a complicated action under conjugation.
\begin{align*}
    \traceTilde(\matSP{N}^{-1}\matSP{M}\matSP{N}) =& +\sigmaT{H}AE+\sigmaT{H}BG-\sigmaT{F}CE-\sigmaT{F}DG\\
    & -\sigmaT{F}\sigmaT{A}G -\sigmaT{H}\sigmaT{B}G +\sigmaT{F}\sigmaT{C}E+\sigmaT{H}\sigmaT{D}E\\
    =& \sigmaT{H}(A+\sigmaT{D})E+\sigmaT{H}(B-\sigmaT{B})G -\sigmaT{F}(C-\sigmaT{C})E -\sigmaT{F}(\sigmaT{A}+D)G
\end{align*}
We then get a similar, albeit a bit simpler, proposition to \Cref{thm:traceConjugation} for the twisted trace. 
\begin{proposition}
    If $A+\sigmaT{D}$ is in the center of $\A$ and $B,C,A+\sigmaT{D}\in \A^\sigma$ then for any $\matSP{N}\in \SPA$, $\traceTilde(\matSP{N}^{-1}\matSP{M}\matSP{N}) = \traceTilde(\matSP{M})$.
\end{proposition}
\begin{proof}
    With the assumptions on the entries of $\matSP{M}$, $\traceTilde(\matSP{N}^{-1}\matSP{M}\matSP{N})$ simplifies to $(A+\sigmaT{D})(\sigmaT{H}E-\sigmaT{F}G)$. Thus when $\matSP{N} \in \SPA$ the proposition holds.
\end{proof}

The determinant is also ill behaved under multiplication. The following proposition follows by direct calculation.

\begin{proposition}
    Let $\matSP{M} = \begin{bmatrix}
    A & B\\C & D
\end{bmatrix}$ and $\matSP{N} = \begin{bmatrix}
    E & F\\G & H
\end{bmatrix}$ $\in \Mat_2(A)$ Then 
$$\det\nolimits_l(\matSP{M}\matSP{N}) = \sigma(E)\det\nolimits_l(\matSP{M})H - \sigma(G)\sigma(\det\nolimits_l(\matSP{M}))F + \sigma(E)\twistedCommutatorLeft{A}{C}F + \sigma(F)\twistedCommutatorLeft{B}{D}H $$
\end{proposition}

\begin{corollary}
    If $\matSP{M} \in \Mat_2(\A,\sigma)$ with left $\sigma$-commuting columns and $\det_l(\matSP{M}) \in \A^\sigma$ and in the center of $\A$ then $\det_l(\matSP{M}\matSP{N}) = \det_l(\matSP{M})\det_l(\matSP{N})$.
\end{corollary}

\begin{remark}
    We can take the conditions of the previous Corollary to define a group $\GL_2(\A,\sigma)$ as an analogue to the general linear group. However the assumptions that the determinant is central and fixed by $\sigma$ are very strong. Since $\GL_2(\A,\sigma) \simeq \SPA \times Z(\A^\sigma)^\times$, we focused on $\SPA$ to obtain cleaner formulas.
\end{remark}

In the commutative situation a well known relation between the determinant and the trace is that the trace is the first order approximation of the determinant. The twisted trace $\traceTilde$ is the first order approximation of $\det_l$. Explicitly \[\det\nolimits_l(\matSP{M}+\epsilon \matSP{Y}) = \det\nolimits_l(\matSP{M}) + \traceTilde(\Adj(\matSP{M})\matSP{N})\epsilon + O(\epsilon^2).\]

\section{Trace Identities}\label{sec:CayleyHamilton}
In this section we prove several identities for the trace function. We begin by proving a Cayley-Hamilton theorem for $2\times 2$ matrices over involutive rings. From this we establish identities for the commutators and powers of matrices.\\
We begin with a few simple lemmas which make the subsequent calculations easier.

\begin{lemma}
    Let $\matSP{M}\in \Mat_2(\A,\sigma)$ Then we have 
    $$ \sigma(\matSP{M}) + \Adj(\matSP{M}) = \sigma(\trace(\matSP{M}))\id$$
\end{lemma}

\begin{corollary}[Classical Trace Identity]\label{thm:TraceIdentityCayleyHamilton}
    For $\matSP{M},\matSP{N}\in \Mat_2(\A)$ 
    \[\trace(\Adj(\matSP{M})\matSP{N}) + \trace(\sigma(\matSP{M})\matSP{N})-\sigma(\trace(\matSP{M}))\trace(\matSP{N}) =0\]
\end{corollary}

\begin{lemma}
     $$\Adj(\matSP{M})\matSP{M} = \begin{bmatrix} \sigma(\det_l(\matSP{M})) & \twistedCommutatorLeft{D}{B}\\ \twistedCommutatorLeft{A}{C} & \det_l(\matSP{M})\end{bmatrix} $$
\end{lemma}

\begin{theorem}[(left) Cayley-Hamilton]\label{thm:CayleyHamilton}
Let $\matSP{M} \in \Mat_2(\A,\sigma)$ then 
\[\sigma(\matSP{M})\matSP{M}-\sigma(\trace(\matSP{M}))\matSP{M} = -\begin{bmatrix} \sigma(\det_l(\matSP{M})) & \twistedCommutatorLeft{D}{B}\\ \twistedCommutatorLeft{A}{C} & \det_l(\matSP{M})\end{bmatrix}
\]

In particular an element $\matSP{M} \in \Mat_2(\A,\sigma)$ lies in $\SP_2(\A,\sigma)^l$ if and only if 
\begin{equation}\label{eqn:SymplecticCayleyHamilton}
    \sigma(\matSP{M})\matSP{M} - \sigma(\trace(\matSP{M}))\matSP{M} + \id = 0
\end{equation}
\end{theorem}
\begin{proof}
This follows from direct computation using the previous two lemmas.
%     \begin{align*}
%         \begin{bmatrix} \sigma(A) & \sigma(B)\\ \sigma(C) & \sigma(D)\end{bmatrix} & \begin{bmatrix} A & B\\ C & D \end{bmatrix} - (\sigma(A) + \sigma(D))\begin{bmatrix}A & B\\ C & D\end{bmatrix}\\
%         =& \begin{bmatrix}
%             \sigma(A)A+\sigma(B)C & \sigma(A)B+\sigma(B)D\\
%             \sigma(C)A+\sigma(D)C & \sigma(C)B+\sigma(D)D
%         \end{bmatrix} - \begin{bmatrix}
%             \sigma(A)A + \sigma(D)A & \sigma(A)B+\sigma(D)B\\
%             \sigma(A)C+\sigma(D)C & \sigma(A)D + \sigma(D)D
%         \end{bmatrix}\\
%         =& \begin{bmatrix}
%             \sigma(B)C-\sigma(D)A & \sigma(B)D-\sigma(D)B\\
%             \sigma(C)A-\sigma(A)C & \sigma(C)B-\sigma(A)D
%         \end{bmatrix}
%     \end{align*}
\end{proof}

\begin{remark}
    Note that although $\sigmaT{\matSP{M}}$ appears in the Cayley-Hamilton formula, $\sigmaT{\matSP{M}}$ is not an element of $\SPA$.
\end{remark}
When $\sigma$ is trivial \Cref{eqn:SymplecticCayleyHamilton} is exactly the classic Cayley-Hamilton theorem for $\SL_2(\R)$. Furthermore the characterization of $\SPA^l$ via this equation mirrors the characterization of Manin matrices as matrices $m$ which satisfy the classic Cayley Hamilton theorem:
\[m^2 - (A+D)m + (AD-CB)\id = 0.\]

\begin{remark}
    We also have a (right) version of Cayley-Hamilton 
    \begin{equation}\label{eqn:CayleyHamiltonRow}
        \matSP{M}\sigmaT{\matSP{M}}-\matSP{M}\sigma(\trace(\matSP{M})) = -\begin{bmatrix} \det(\sigma(\matSP{M})^T) & \twistedCommutatorRight{A}{B}\\ \twistedCommutatorRight{D}{C} & \det_\sigma(\sigma(\matSP{M})^T)\end{bmatrix}
    \end{equation}
    This can be seen by replacing $\matSP{M} $ with $\sigma(\matSP{M})^T$ and applying $\sigma()^T$ to the entire equation. 
\end{remark}

% We can now deduce several trace identities from the Cayley-Hamilton Theorem. Each of these equations have a corresponding right version
% The first is the classical $\SL_2$-trace identity.  
%     For $\SL_2(\R)$ this identity is key for writing the trace of any word in $x,y$ as a polynomial in $\trace(x),\trace(y),\trace(xy)$.

We can similarly obtain a trace relation for appropriate twisted commutators.
\begin{corollary}[Trace Identity for commutators]
Let $\matSP{X},\matSP{Y} \in \SPA^r$
\begin{align*}
    \trace(\matSP{X}\sigmaT{\sigmaT{\matSP{Y}}\sigmaT{\matSP{X}^{-1}}}\matSP{Y}^{-1}) = &\trace(\matSP{X})\sigmaT{\trace(\matSP{X})}+\trace(\matSP{Y})\sigmaT{\trace(\matSP{Y})}+\trace(\matSP{X}\sigmaT{\sigmaT{\matSP{Y}}\matSP{X}}\sigmaT{\matSP{Y}}) \\ & -\trace(\matSP{X})\trace(\sigmaT{\matSP{X}}\matSP{Y})\sigmaT{\trace(\matSP{Y})}
\end{align*}
\end{corollary}
\begin{proof}
Let $X=\trace(\matSP{X})$, $Y=\trace(\matSP{Y})$ and $Z=\trace(\sigma(\matSP{X})Y)$. We are then trying to show that 
\[\trace(\matSP{X}\sigmaT{\sigmaT{\matSP{Y}}\sigmaT{\matSP{X}^{-1}}}\matSP{Y}^{-1}) = X\sigmaT{X}+Y\sigmaT{Y}+\trace(\matSP{X}\sigmaT{\sigmaT{\matSP{Y}}\matSP{X}}\sigmaT{\matSP{Y}})-XZ\sigmaT{Y}\]
Now using our previous lemmas
\begin{align*}
    \matSP{X}&\sigmaT{\sigmaT{\matSP{Y}}\sigmaT{\matSP{X}^{-1}}}\matSP{Y}^{-1} = \matSP{X}\sigmaT{\sigmaT{\matSP{Y}}(X\id-\matSP{X})}\matSP{Y}^{-1} \\
    &= \matSP{X}\sigmaT{X} - \matSP{X}\sigmaT{\sigmaT{\matSP{Y}}\matSP{X}}\matSP{Y}^{-1} \\
    &= \matSP{X}\sigmaT{X} - \matSP{X}\sigmaT{\sigmaT{\matSP{Y}}\matSP{X}}(\sigmaT{Y}\id-\sigmaT{\matSP{Y}}) \\
    &= \matSP{X}\sigmaT{X} - \matSP{X}\sigmaT{\sigmaT{\matSP{Y}}\matSP{X}}\sigmaT{Y} + \matSP{X}\sigmaT{\sigmaT{\matSP{Y}}\matSP{X}}\sigmaT{\matSP{Y}}) \\
    &= \matSP{X}\sigmaT{X} - (X\id-\sigmaT{\matSP{X}^{-1}})\sigmaT{\sigmaT{\matSP{Y}}\matSP{X}}\sigmaT{Y} + \matSP{X}\sigmaT{\sigmaT{\matSP{Y}}\matSP{X}}\sigmaT{\matSP{Y}}) \\
    &= \matSP{X}\sigmaT{X} - X\sigmaT{\sigmaT{\matSP{Y}}\matSP{X}}\sigmaT{Y}+ \sigmaT{\matSP{X}^{-1}}\sigmaT{\sigmaT{\matSP{Y}}\matSP{X}}\sigmaT{Y} +\matSP{X}\sigmaT{\sigmaT{\matSP{Y}}\matSP{X}}\sigmaT{\matSP{Y}})
\end{align*}
To complete the proof we take trace and focus on the term  \[\trace(\sigmaT{\matSP{X}^{-1}}\sigmaT{\sigmaT{\matSP{Y}}\matSP{X}}\sigmaT{Y}) = \trace(\sigmaT{\matSP{X}^{-1}}\sigmaT{\sigmaT{\matSP{Y}}\matSP{X}})~\sigmaT{Y}.\]
Now we apply \Cref{prop:sigma_commute_trace} to obtain
\[\trace(\sigmaT{\matSP{X}^{-1}}\sigmaT{\sigmaT{\matSP{Y}}\matSP{X}})~\sigmaT{Y} = \sigma(\trace(\sigmaT{\matSP{Y}}\matSP{X} \matSP{X}^{-1}))~\sigmaT{Y} = Y\sigmaT{Y}\]
We then have that 
\begin{align*}
    \trace(\matSP{X}\sigmaT{\sigmaT{\matSP{Y}}\sigmaT{\matSP{X}^{-1}}}\matSP{Y}^{-1}) =&X\sigmaT{X}+Y\sigmaT{Y}+\trace(\matSP{X}\sigmaT{\sigmaT{\matSP{Y}}\matSP{X}}\sigmaT{\matSP{Y}}) -XZ\sigmaT{Y}
\end{align*}
Finally note that we only used $\matSP{X}(\matSP{X}^{-1}) = \id$ and $\matSP{Y}(\matSP{Y}^{-1}) = \id$ in our calculations justifying that we only need $\matSP{X},\matSP{Y}$ to be in the right symplectic group. 
%% Proof using X as matrices and x as trace
%Let $x=\trace(X)$ and $y=\trace(Y)$. Using equation 7 we see that $\sigmaT{X}=\sigmaT{x}Id-X^{-1}$. Now using this fact we write 
%\begin{align*}
%    X&\sigmaT{\sigmaT{Y}\sigmaT{X^{-1}}}Y^{-1} = X\sigmaT{\sigmaT{Y}(xId-X)}Y^{-1} \\
%   &= X\sigmaT{x} - X\sigmaT{\sigmaT{Y}X}Y^{-1} \\
%    &= X\sigmaT{x} - X\sigmaT{\sigmaT{Y}X}(\sigmaT{y}\id-\sigmaT{Y}) \\
%    &= X\sigmaT{x} - X\sigmaT{\sigmaT{Y}X}\sigmaT{y} + X\sigmaT{\sigmaT{Y}X}\sigmaT{Y}) \\
%    &= X\sigmaT{x} - (xId-\sigmaT{X^{-1}})\sigmaT{\sigmaT{Y}X}\sigmaT{y} + X\sigmaT{\sigmaT{Y}X}\sigmaT{Y}) \\
%    &= X\sigmaT{x} - x\sigmaT{\sigmaT{Y}X}\sigmaT{y}+ \sigmaT{X^{-1}}\sigmaT{\sigmaT{Y}X}\sigmaT{y} \\ & \quad+X\sigmaT{\sigmaT{Y}X}\sigmaT{Y})
%\end{align*}
%Now since $\trace(\sigmaT{X^{-1}}\sigmaT{\sigmaT{Y}X}\sigmaT{y}) = \sigmaT{\trace(\sigmaT{Y}XX^{-1}})\sigmaT{y}=y\sigmaT{y}$
%We have that 
%\begin{align*}
%    \trace(X\sigmaT{\sigmaT{Y}\sigmaT{X^{-1}}}Y^{-1}) =&x\sigmaT{x}+y\sigmaT{y}+\trace(X\sigmaT{\sigmaT{Y}X}\sigmaT{Y}) \\ 
%    &-x~\trace(\sigmaT{X}Y)\sigmaT{y}
%\end{align*}
\end{proof}

We also obtain a relation for the trace of powers of $\matSP{M}$. For convenience we write $\matSP{M}^{\sigma k}$ to mean the alternating product of $\sigmaT{\matSP{M}}$ and $\matSP{M}$ ending with $\matSP{M}$. 
\[\matSP{M}^{\sigma 1} = \matSP{M} \hspace{1pc} \matSP{M}^{\sigma 2} = \sigmaT{\matSP{M}}\matSP{M} \hspace{1pc} \matSP{M}^{\sigma 3} = \matSP{M}\sigmaT{\matSP{M}}\matSP{M} \text{ etc.}\]
\begin{corollary}[Trace Identity for Powers]
    We have $\trace(\matSP{M}^{\sigma k}) = U_k(\trace(\matSP{M}))$ where $U_k(x)$ are the Chebyshev polynomials of the second kind with powers interpreted as alternating $x$, $\sigmaT{x}$ as above.
\end{corollary}
\begin{proof}
    It suffices to show that $\trace(\matSP{M}^{\sigma k})$ satisfies the recurrence for the Chebyshev polynomials of the second kind, $U_{k+1}(x) = U_1(x)U_k(x)-U_{k-1}(x)$. We have two cases depending on whether $k$ is even or odd. If $k$ is odd then $\matSP{M}^{\sigma(k+1)} = \sigmaT{\matSP{M}}\matSP{M}^{\sigma k} = \sigmaT{\matSP{M}}\matSP{M} \matSP{M}^{\sigma(k-1)}$. Then using \Cref{thm:TraceIdentityCayleyHamilton} we have
    \[\trace(\sigmaT{\matSP{M}}\matSP{M}^{\sigma k}) = \trace(\sigmaT{\matSP{M}})\trace(\matSP{M}^{\sigma k}) -\trace(\matSP{M}^{-1}\matSP{M} \matSP{M}^{\sigma(k-1)}) = \sigmaT{\trace(\matSP{M})}\trace(\matSP{M}^{\sigma k}) - \trace(\matSP{M}^{\sigma(k-1)}) = \sigmaT{U_1}U_k - U_{k-1}\]
    Similarly if $k$ is even $\matSP{M}^{\sigma(k+1)} = \matSP{M} \matSP{M}^{\sigma k} = \matSP{M}\sigmaT{\matSP{M}}\matSP{M}^{\sigma(k-1)}$ and 
    \[\trace(\sigmaT{\matSP{M}}\matSP{M}^{\sigma k}) = \trace(\matSP{M})\trace(\matSP{M}^{\sigma k}) - \trace(\matSP{M}^{\sigma(k-1)}) = U_1 U_k - U_{k-1}\]
\end{proof}

Note that we really need $\matSP{M}\in \SPA$ and not just in $\SPA^l$ since we use both left and right trace relations.

% \begin{definition}
%     We say $f\in\SP_2(\A,\sigma)$ is \keyword{$\sigma$-stable} if $\sigma(f)$ is also in $\SP_2(A,\sigma)$
% \end{definition}
% \begin{remark}
%     A matrix $\begin{bmatrix}
%         A & B\\C & D
%     \end{bmatrix}$ in $\SP_2(\A,\sigma)$ is $\sigma$-stable if
%     \begin{align*}
%         A\sigma(C) - C\sigma(A) =& 0\\
%         B\sigma(D) - D\sigma(B) =& 0\\
%         A\sigma(D)-C\sigma(B) =& 1
%     \end{align*}
% \end{remark}
% Note that $\sigma$-stablility is not preserved under product. However $\sigma$-stable matrices satisfy even more $\SL_2(\R)$ style trace identities.

\subsection{Twisted  Abelianization}\label{sec:quotient}

It might be a bit unsatisfactory that we have the different trace and determinant functions which are all different. 
Some of this can be remedied if we pass to an appropriate quotient object. For this we recall the left/right $\sigma$-twisted commutators (\Cref{def:SigmaCommute}). 

\begin{definition}
We denote subset of $\A$ that is generated as an abelian group by right $\sigma$-twisted commutators by $\twistedCommutatorRight{\A}{\A}$. Similarly we write $\twistedCommutatorLeft{\A}{\A}$ for the subset of $\A$ generated as an abelian group by left $\sigma-$twisted commutators. 
\end{definition}

We can then consider the quotient  
$$\overline{A} = \A/\twistedCommutatorRight{\A}{\A}.$$ Note that this is not a quotient as a ring, but just as abelian groups. 

Since $\A$ is unital, in this quotient any element $A \in \A $ is identified with $\sigma(A)$. 
\begin{remark}
    If we take the quotient using left $\sigma-$twisted commutators the quotient abelian group is the same.
\end{remark}
When $\A$ is an algebra over a field of characteristic not equal to $2$, we can in fact identify $\A/\twistedCommutatorRight{\A}{\A}$ with $A^\sigma$. 

Under this identification, the natural quotient map is given by $\A \rightarrow\A/\twistedCommutatorRight{\A}{\A} \simeq \A^\sigma$ is given by $A \mapsto \frac{1}{2} (A + \sigma(A))$.

\begin{remark}
When $\sigma$ is trivial (and hence $\A$ is commutative), the twisted commutator is just the commutator, and $\A/[\A, \A]= \A$. 
However when $\sigma$ is non-trivial the twisted commutator is quite different, and $\A/[\A, \A]^\sigma = \A^\sigma$ is strictly smaller than $\A$. 

On the other hand when $\A$ is noncommutative, $\A/[\A, \A]^\sigma  = A^\sigma$ can be much bigger than $\A/[\A, \A]$, take for example take the case $\A = \Mat_n(k)$.  
\end{remark}

We introduce notation for the composition of the determinant or the trace function with the quotient map $\A \rightarrow\A/\twistedCommutatorRight{\A}{\A}$. If $F:\Mat_2(\A,\sigma) \rightarrow A$ is a function, we denote the corresponding map from $\Mat_2(\A,\sigma) \rightarrow \A/\twistedCommutatorRight{\A}{\A}$ by $\overline{F}$.  

\begin{lemma}\label{lem:quotientequal}
We have 
$\overline{\det_l} = \overline{\sigma \circ \det_l}$, $\overline{\det_r} = \overline{\sigma \circ \det_r}$ and 
$\overline{\trace} = \overline{\sigma\circ \trace}$. 
If $\A$ is an algebra over a field not of characteristic $2$, we further have 
$\overline{\trace}  = \overline{\traceTilde} = \overline{\sigma \circ \trace}$. In particular, $\overline{\trace}$ is the first order approximation of $\overline{\det_l}$. 
\end{lemma}

Now consider the space $\Mat_2(\overline{\A}, \sigma)$. There is a natural quotient map 
$\Mat_2(\A, \sigma) \rightarrow \Mat_2(\overline{\A}, \sigma)$, by passing to the quotient componentwise. 

We then have an obvious corollary to the Cayley Hamilton theorem established in $\Mat_2(\A,\sigma)$.
\begin{corollary}[Cayley-Hamilton in quotient]\label{cor:CayleyHamilton}
Let $\matSP{M} \in \Mat_2(\A, \sigma)$ then 
$$\overline{[\sigma(\matSP{M})\matSP{M}-\sigma(\trace(\matSP{M}))\matSP{M} ]}= -\overline{\det\nolimits_l}(\matSP{M})\id.$$
\end{corollary}

\section{Quantization}\label{sec:Quantization}

In this section we will define the quantization of $\SPA$ which parallels the quantization of $\SL_2$  \cite{Manin-QuantumGroupsRemarks,Manin-QuantumGroupsNoncomGeometry,Kassel-QuantumGroups}. 
\subsection{Quantum \texorpdfstring{$\SL_2$}{SL2}}\label{sec:QuantumSL2}
First we review the quantization of $\SL_2$. Informally this is a deformation in a parameter $q$ of the ring of algebra of matrix functions on $\SL_2$. \\
Formally we fix a field $k$. Let $k_q = k(q)$ be the extension of $k$ by an element $q$. Often $q$ is taken to be an element of $k$, although this introduces additional subtlety if $q$ is a root of unity.
\begin{definition}
    The \keyword{(coordinate ring of the) quantum plane} $\A_q$ is the algebra $k_q[x_1,x_2]$ such that $x_1$ and $x_2$ q-commute, i.e. $x_1x_2 = q^{-1}x_2x_1$.
\end{definition}
Intuitively, the quantization of $\SL_2(k)$ should act like the endomorphisms of the quantum plane. 
\begin{definition}
    The \keyword{algebra of quantum matrix functions} $\CoordRing(\SL_2)_q$ is the ring $k_q\{a,b,c,d\}/I$ where $I$ is generated by the quantum $\SL_2$ relations:
    \begin{align}\label{eqn:QuantumSl2}
\begin{aligned}
     ac-q^{-1}ca =& 0  \hspace{3pc}& ab-q^{-1}ba =& 0\\
     bd-q^{-1}db =& 0 & cd-q^{-1}dc =& 0\\
     ad-q^{-1}cb =& da-qbc  &  bc=cb \\
     ad-q^{-1}cb =& 1
\end{aligned}
\end{align}
\end{definition}
These equations can be obtained by adding the appropriate power of $q$ to the equations for a matrix and its transpose to be Manin (\ref{eqn:ManinCommutivityEquations}).  To better compare with $\SPA$ we can rewrite the quantum $\SL_2$ relations as follows
\begin{align}\label{eqn:QuantumSL2rephrased}
    \begin{aligned}
        ac-q^{-1}ca =& 0 \hspace{3pc} & ab-q^{-1}ba =& 0\\
        bd-q^{-1}db =& 0 & cd-q^{-1}dc =& 0\\
        ad-q^{-1}cb =& 1  & ad-q^{-1}bc =&1 \\
        da-qbc =&1 & da-qcb =&1
    \end{aligned}
\end{align}

This algebra is actually a Hopf algebra with a coproduct $\Delta : \SL_2(R)_q \rightarrow \SL_2(R)_q \tensor \SL_2(R)_q$ and antipode $S: \SL_2(R)_q \to \SL_2(R)_q $ inspired by matrix multiplication and inverse. This means the coproduct on the generators is defined by reading the following componentwise:
\[\Delta\left( \begin{bmatrix}
    a & b \\ c & d
\end{bmatrix} \right)= \begin{bmatrix}
    a & b \\ c & d
\end{bmatrix} \tensor \begin{bmatrix}
    a & b \\c & d
\end{bmatrix} = \begin{bmatrix}
    a \tensor a + b\tensor c & a\tensor b + b \tensor d\\
    c\tensor a + d\tensor c & c \tensor b + d \tensor d
\end{bmatrix}
\] and the antipode defined by
\[S\left(\begin{bmatrix}
    a & b \\ c & d
\end{bmatrix}\right)=  \begin{bmatrix}
    d & -qb \\ -q^{-1}c & a
\end{bmatrix}.
\] 
This coproduct allows Manin to classify the algebra of quantum matrix functions by the coaction on the quantum plane in an analogous manner to the classical case.

\subsection{Quantum Symplectic Group}

Usually when we define a quantum deformation of an algebra, we  are constructing a noncommutative algebra out of a commutative one. Here we already start with a noncommutative algebra, so its quantization is a much less drastic change. In fact the corresponding quantization can be seen simply as a change of symplectic form.\\
First we add a quantum parameter $q$ to our algebra. Formally let $\A_q = \A[q^{\pm}]$ with the additional assumption that $q$ commutes with all of $\A$. We then extend $\sigma$ to $\A_q$ by setting $\sigma(q) = q^{-1}$.

\begin{definition}
    Two elements $X,Y \in \A_q$ \keyword{left $q-\sigma-$commute} if $\sigmaT{X}Y = q^{-1}\sigmaT{Y}X  $.
\end{definition}
 Following the quantization of $\SL_2(\R)$ we define a system of equations on the matrix entries of a $2\times 2$ matrix specifying the quantum commutativity and quantum determinant of the elements.
 \begin{definition}
    The \keyword{quantum determinant} of $M \in \Mat_2(\A_q)$ is $\sigmaT{A}D-q^{-1}\sigmaT{C}D$.
\end{definition}
 \begin{definition}
     The \keyword{quantum symplectic equations} for a matrix $\begin{bmatrix}
         A & B\\ C &D
     \end{bmatrix} \in \Mat_2(\A_q)$ are
       \begin{align}\label{eqn:QuantumSymplectic}
    \begin{aligned}
         \sigma(A)C-q^{-1}\sigma(C)A =& 0  \hspace{3pc}& A\sigma(B)-q^{-1}B\sigma(A) =& 0\\
         \sigma(B)D-q^{-1}\sigma(D)B =& 0 &  C\sigma(D)-q^{-1}D\sigma(C) =& 0\\
         \sigma(A)D-q^{-1}\sigma(C)B =& 1 &  A\sigma(D)-q^{-1}B\sigma(C) =& 1
    \end{aligned}
    \end{align}
 \end{definition}
 \begin{remark}
     These equations correspond both to the standard symplectic equations (\ref{eq:sp2a_quations_l},\ref{eq:sp2a_quations_r}) and to the rephrased quantum $\SL_2(\R)$ equations (\ref{eqn:QuantumSL2rephrased}). Note that these equations don't imply any commutativity between $B$ and $C$. 
 \end{remark}
As with $\SPA$, these equations are equivalent to $\matSP{M}$ having a two sided inverse and preserving a particular skew symmetric bilinear form. Given 
\[\matSPgreek{\Omega}_q = \begin{bmatrix} 0& 1\\ -q^{-1} & 0\end{bmatrix} \in \Mat_2(\A_q)\]
we define an $\A_q$ valued skew symmetric right-bilinear form $\omega_q : \A_q^2 \times \A_q^2 \rightarrow \A_q$ by
\[\omega_q(\vec{v},\vec{w}) = \sigmaT{\vec{v}}^{T}\Omega_q \vec{w}\]
\begin{remark}
    If we take $q =1$ we recover $\Omega_1$ and the same form we defined in \Cref{sec:SP2Asigma}.
\end{remark}

\begin{definition}
    The \keyword{left quantum symplectic monoid}, $\SPA_q^l$ is the set of endomorphisms $f$ of $\A_q^2$ preserving $\omega_q$. Equivalently 
    \[\SPA_q^l = \{f\in\End(\A_q^2) \mid \forall \vec{v},\vec{w}\in\A_q^2. \omega_q(f(\vec{v}),f(\vec{w})) = \omega(\vec{v},\vec{w})\} = \{\matSP{M} \in\Mat_2(\A_q) \mid \sigmaT{\matSP{M}}^T\matSPgreek{\Omega}_q \matSP{M} = \matSPgreek{\Omega}_q\}\]
    As in the classical case, the \keyword{right quantum symplectic monoid} $\SPA_q^r$ are the elements for which $\sigmaT{\matSP{M}}^T \in \SPA_q^l$. 
    Then the \keyword{quantum symplectic group} $\SPA_q$ is again the intersection of the left and right quantum symplectic monoids.  
\end{definition}

With this definition it is clear that $\SPA_q$ is an honest group. It can be characterized in all the same ways that $\SPA$ could. As such we define the quantum plane as follows.

\begin{definition}
    The \keyword{left quantum plane}, $(\A^2_\omega)_q$, is the set of generic left $q$-commuting elements of $\A_q^2$.
    \[(\A^2_\omega)_q = \{\vec{v}=(V_1,V_2) \in \A_q^2 \mid \text{$\vec{v}$ generic and } q^{-1}\sigmaT{V_2}V_1 = \sigmaT{V_1}V_2\}\]
\end{definition}

\begin{theorem}
    The following are equivalent for a matrix $\matSP{M} = \begin{bmatrix}
        A & B\\ C & D
    \end{bmatrix} \in \Mat_2(\A_q)$:
    \begin{enumerate}
        \item\label{thm:SPAqClassification_Defn} $\matSP{M} \in \SPA_q$.
%        \item\label{thm:SPAqClassification_Tranpose} $\sigma(M)^T \in \SPA_q$.
        \item\label{thm:SPAqClassification_Equations} The entries of $\matSP{M}$ satisfy the quantum symplectic equations (\ref{eqn:QuantumSymplectic}).
%        \item\label{thm:SPAqClassification_HalfEquations} The entries of $M$ satisfy the left half of the quantum symplectic equations.
        \item\label{thm:SPAqClassification_InverseFormula} $\matSP{M}$ has a two sided inverse: $\matSP{M}^{-1} = \begin{bmatrix}
            \sigmaT{D} & -q \sigmaT{B}\\
            -q^{-1}\sigmaT{C} & \sigmaT{A}
        \end{bmatrix}$.
        \item\label{thm:SPAqClassification_ActionOnPlane} $M$ has quantum determinant $1$ and $\matSP{M}$ preserves the quantum plane.
    \end{enumerate}
\end{theorem}
\begin{proof}
    These proofs are identical to proofs for $\SPA$ already presented. See \Cref{thm:SymplecticEquations} for (\ref{thm:SPAqClassification_Defn}) equivalent to (\ref{thm:SPAqClassification_Equations}). \Cref{thm:SP2AInverse} shows the equivalence of (\ref{thm:SPAqClassification_InverseFormula}). %and \Cref{thm:SP2A_Transpose} shows the equivalence of (\ref{thm:SPAqClassification_Tranpose}). %The proof that (\ref{thm:SPAqClassification_HalfEquations}) implies (\ref{thm:SPAqClassification_Equations}) follows from the fact that applying the left half the quantum symplectic equations to $\sigma(M)^T$ yields the right half. 
    Finally the equivalence of (\ref{thm:SPAqClassification_ActionOnPlane}) corresponds to \Cref{thm:SP2ACharacterization_ActionOnSigmaPlane} noting that $(1,0)$ and $(0,1)$ are both in the quantum plane.
\end{proof}

 \begin{remark}
     In \cite{Manin-QuantumGroupsRemarks} to obtain an action on the affine plane, the additional assumption that $V_1$ and $V_2$ commuted with $(A,B,C,D)$ was needed. However with the $\sigma$ commuting plane we don't need this assumption and $\SP_2(\A,\sigma)_q$ acts on the entire quantum plane. 
 \end{remark}

   The following theorem shows that $\SPA_q$ has many elements if we assume that $q$ has a square root $q^{1/2}$. This is very different from the usual quantization of $\SL_2(\R)$ where the corresponding group is just $\left\{\begin{bmatrix}
       a & 0\\0 & a^{-1}
   \end{bmatrix} \mid a \in \R\right\}.$
\begin{theorem}\label{thm:SP2A_elements}
    If $\matSP{M} = \begin{bmatrix}
        A & B\\ C &D 
    \end{bmatrix} \in \SPA$ then $\begin{bmatrix}
        A & q^{1/2}B\\q^{-1/2}C & D
    \end{bmatrix} \in \SPA_q$
\end{theorem}
\begin{proof}
    This can be verified by checking the quantum symplectic equations (\ref{eqn:QuantumSymplectic}). We check the first equation here and the other computations are similar.
    \begin{align*}
        \sigmaT{A}q^{-1/2}C  - q^{-1}\sigmaT{q^{-1/2}C}A =& q^{-1/2}\sigmaT{A}C - q^{-1}q^{1/2}\sigmaT{C}A = q^{-1/2}(\sigmaT{A}C-\sigmaT{C}A) = 0
    \end{align*}
\end{proof}

To highlight the difference consider the example $(\A,\sigma) = (\R,\id)$. Let $q = e^{ih} \neq 1$. Then $\A_q \simeq \C$ the extension of  $\sigma$ to $\A_q$ is complex conjugation. We check that $\sigma(q) = \overline{q} = e^{-ih} = q^{-1}$ as needed. Following the recipe of \Cref{thm:SP2A_elements} we see that $\begin{bmatrix}
    a & be^{ih/2}  \\ ce^{-ih/2} & d\\
\end{bmatrix}$ is an element of $\SP(\R,\id)_q$. However this matrix does not satisfy the standard quantum equations. In particular 
\[ a(ce^{-ih/2}) - e^{-ih}(ce^{-ih/2})a = ace^{-ih/2}(1-e^{-ih}) \neq 0\]
For this reason quantum $\SL_2$ is defined by quantizing the algebra of matrix functions instead of the matrices themselves. To define the analogous objects for the quantum symplectic group we need to extend the definition of an involutive ring to an involutive algebra.

\section{The Ring of Noncommutative Matrix Functions}\label{sec:NoncomMatrixFunctions}

In \Cref{sec:QuantumSL2} we saw that the quantization of $\SL_2$ is defined as a quotient of a noncommutative algebra. We think of this algebra as the algebra of matrix functions on $(\SL_2)_q$. In order to define the analogous object for $\SPA_q$ we must first be precise about the extra structure the algebra of quantum symplectic matrix functions should have.

\subsection{Involutive  Algebras}\label{sec:KsigmaAlgebra}

Let $\K$ be a commutative ring and $\sigma \colon \K \rightarrow \K$ an involution.

\begin{definition}
    A \keyword{$(k,\sigma)$-algebra} is a (noncommutative) $k$ algebra $\B$ along with a ring homomorphism $\sigma_B: \B \to \B^{op}$ restricting to $\sigma$ on $k$.\\
    For two $(k,\sigma)$-algebras $(B,\sigma_B),(C,\sigma_C)$ a \keyword{homomorphism of $(k,\sigma)$-algebras} is a map $\phi \colon B \rightarrow C$ such that the following diagram commutes:
    \begin{center}
    \begin{tikzcd}
        B\arrow[r,"\phi"]\arrow[d,"\sigma_B"] & C \arrow[d,"\sigma_C"]\\
        B^{op}\arrow[r,"\phi"] & C^{op}
    \end{tikzcd}
    \end{center}
    I.e. $\forall b \in B$ $\phi(\sigma_B(b)) = \sigma_C(\phi(b))$.
\end{definition}
\begin{definition}
   A \keyword{module over a $(k,\sigma)$-algebra}, $\B$, is a $\B$ bimodule $M$ along with a linear map $\sigma_M:M \to M$ such that for $b_1,b_2\in \B,m \in M$ $$\sigma_M(b_1mb_2)= \sigma_B(b_2)\sigma_M(m)\sigma_B(b_1)$$ 
\end{definition}
We will suppress the subscripts on $\sigma$ when it is understood which algebra/module they are acting on. 
\begin{definition}
    The \keyword{free $(k,\sigma)$-algebra}, $k\{a_1,\dots,a_n\}$, is the free algebra $k\{a_1,\dots,a_n,\hat{a}_1,\dots,\hat{a}_n\}$ with the anti-involution $\sigma$ induced by $\sigmaT{a_i} = \hat{a}_i$ and $\sigmaT{\hat{a}_i} = a_i$. 
\end{definition}
\begin{lemma}
    A $(k,\sigma)$-homomorphism $f$ from a free algebra $k\{a_1,\dots,a_n\}$ to $(\B,\sigma)$ is uniquely defined by specifying the image of $a_1,\dots,a_n$.
\end{lemma}
\begin{proof}
    For a free algebra without involution the image of $\hat{a}_1,\dots,\hat{a}_n$ would be freely chosen. However to be a $(k,\sigma)$-homomorphism we must have $f (\hat{a}_i) = f(\sigmaT{a_i}) = \sigmaT{f(a_i)}$.
\end{proof}

\begin{definition}
    A \keyword{$(k,\sigma)$-ideal} of a $(k,\sigma)$-algebra is a usual two-sided ring ideal $I$ with the additional property that $\sigmaT{I} \subseteq I$.
\end{definition}
\begin{lemma}
    The kernel of a $(k,\sigma)$-homomorphism is a $(k,\sigma)$-ideal.
\end{lemma}
\begin{proof}
    It suffices to check that if $f(a) = 0$ then so does $f(\sigmaT{a})$. Since $f$ is a $(k,\sigma)$-morphism $f(\sigmaT{a}) = \sigmaT{f(a)} = \sigmaT{0} = 0$.
\end{proof}
We can now define $(k,\sigma)$-algebras by taking quotients of free $(k,\sigma)$-algebras by a $(k,\sigma)$-ideal as usual. Such quotient modules satisfy the usual universal property:
\begin{lemma}\label{thm:KsigmaQuotientHomimorphisms}
    Let $I$ be a $(k,\sigma)$-ideal of the free $(k,\sigma)$-algebra $k\{a_1,\dots,a_n\}$ generated by relations $R_i$. Then a $(k,\sigma)$-homomorphism $f$ from $k\{a_1,\dots,a_n\} \rightarrow \B$ is uniquely specified choosing the images $b_1=f(a_1)$,\dots,$b_n=f(a_n)$ satisfying the relations $R_i$ in $\B$.
\end{lemma}
In view of this lemma we make the following definition:
\begin{definition}
    For a $(k,\sigma)$-algebra $(\B,\sigma)$, a \keyword{$(\B,\sigma)$-point of $k\{a_1,\dots,a_n\}/I$} is a $(k,\sigma)$-homomorphism from $k\{a_1,\dots,a_n\}/I \rightarrow \B$.
\end{definition}
This construction lets us define the algebra of matrix functions corresponding to $\SPA$ and $\SPA_q$.
\begin{definition}
    The \keyword{algebra of symplectic matrix functions} $\CoordSPA$ is the $(k,\sigma)$-algebra $k\{A,B,C,D\}/I$ where $I$ is the ideal generated by the defining equations of $\SPA$ (\Cref{eq:sp2a_quations_l,eq:sp2a_quations_r}):
    \begin{align*}
        \begin{aligned}
        \hat{A}C - \hat{C}A =&0 \hspace{3pc} & A\hat{B}-B\hat{A} =&0\\
        \hat{B}D - \hat{D}B =&0 & C\hat{D}-D\hat{C} =&0\\
        \hat{A}D-\hat{C}B =& 1 & A\hat{D}-B\hat{C} =&1\\
        \end{aligned}
    \end{align*}
\end{definition}
Note that as $1$ is fixed by $\sigma$, we also have $\hat{D}A-\hat{B}C = 1$ and so the ideal is closed under $\sigma$.\\
\begin{remark}\label{rem:points_of_coordinate_ring}
    
We are justified in the notation $\CoordSPA$ as we recover $\SPA$ as the $(\A,\sigma)$-points of the algebra. By \Cref{thm:KsigmaQuotientHomimorphisms} we know that a map from $\CoordSPA$ to $\A$ is given by finding four elements $(A,B,C,D)$ in $\A$ satisfying the symplectic equations, i.e. $\begin{bmatrix}
    A &B \\C &D
\end{bmatrix} \in \SPA$. We will see in the next section that by endowing $\CoordSPA$ with the structure of a generalization of a Hopf algebra, we recover the group structure on $\SPA$ as well. 
\end{remark}

We are now going to quantize this algebra by defining the analogous $(k_q,\sigma_q)$ algebra. Recall $k_q = k(q)$ is the extension of $k$ by a formal variable $q$. When $k$ has an anti-involution $\sigma$, there are two choices of extension to $k_q$ given by setting $\sigma_q(q) = q^{\pm 1}$. To define the algebra of quantum symplectic matrix functions it will be necessary to chose $\sigma_q(q) = q^{-1}$.
\begin{definition}
    The \keyword{algebra of quantum symplectic matrix functions}, $\CoordSPA_q$ is the $(k_q,\sigma_q)$ algebra $k_q\{A,B,C,D\}/I_q$ where $I_q$ is the ideal generated by the quantum symplectic equations (\ref{eqn:QuantumSymplectic}):
    \begin{align*}
        \begin{aligned}
         \hat{A}C-q^{-1}\hat{C}A =& 0  \hspace{3pc}& A\hat{B}-q^{-1}B\hat{A} =& 0\\
         \hat{B}D-q^{-1}\hat{D}B =& 0 &  C\hat{D}-q^{-1}D\hat{C} =& 0\\
         \hat{A}D-q^{-1}\hat{C}B =& 1 &  A\hat{D}-q^{-1}B\hat{C} =& 1
        \end{aligned}
    \end{align*}
    The \keyword{left quantum algebra}, $\CoordRing(\A^2_\omega)_q$ is the subalgebra of $\CoordSPA_q$ generated by $A,C$. 
\end{definition}
Note that for $I_q$ to be preserved by $\sigma_q$ we must have chosen $\sigma_q(q) = q^{-1}$.\\
As for $\SPA$ we recover $\SPA_q$ by considering the $(\A_q,\sigma_q)$ points.\\

\subsection{Hopf Algebras}
As $\SL_2(\R)_q$ is a Hopf algebra, to compare with $\CoordSPA_q$ we need to understand it's Hopf algebra structure. However since $(\A,\sigma)$ is itself a noncommutative algebra, the algebra of functions from $\SPA_q$ into $\A$  does not have a Hopf algebra structure in the usual sense. To help understand our new construction we first recall the classic Hopf algebra structure.

 \begin{definition}
     A \keyword{bialgebra} $\alg{H}$ is both a $k$ algebra and a $k$ coalgebra with compatible structures. Explicitly $\alg{H}$ is equipped with the following maps:
     \begin{itemize}
         \item (Multiplication) An associative map $\nabla \colon \alg{H}\tensor \alg{H} \rightarrow \alg{H}$
         \item (Unit) A map  $\eta \colon k \rightarrow \alg{H}$
         \item (Comultiplication) A coassociative map $\Delta \colon \alg{H} \rightarrow \alg{H}\tensor \alg{H}$
         \item (Counit) A map  $\epsilon \colon \alg{H} \rightarrow k$
     \end{itemize}
     The required compatibility conditions are that $\Delta \colon (\alg{H},\nabla,\eta) \rightarrow (\alg{H}\tensor \alg{H}, \nabla_2,\eta_2)$ and $\epsilon\colon (\alg{H},\nabla,\eta) \rightarrow (K,\nabla_0,\eta_0)$ are algebra homomorphisms. The algebra structure on $k$ is clear, $\nabla_0$ is the multiplication on $k$ and $\eta_0$ is the identity. On  $\alg{H}\tensor \alg{H}$ we define $\eta_2 = \eta \tensor \eta$ and $\nabla_2((x_1\tensor y_1) \tensor (x_2 \tensor y_2)) = (x_1x_2)\tensor(y_1y_2)$. 
 \end{definition}
 \begin{definition}
     A \keyword{Hopf algebra} $\alg{H}$ is bialgebra over $k$ with the  additional data of an \keyword{antipode} $S \colon \alg{H} \rightarrow \alg{H}$. The antipode is an anti-homomorphism that makes the following diagram commute:
\begin{center}
    \begin{tikzcd}
                                                                    & \alg{H}\tensor \alg{H} \arrow[rr, "{(S,1)}"] &                      & \alg{H}\tensor \alg{H} \arrow[rd, "\nabla"]  &   \\
\alg{H} \arrow[rr, "\epsilon"] \arrow[rd, "\Delta"'] \arrow[ru, "\Delta"] &                                & k \arrow[rr, "\eta"] &                                & \alg{H} \\
                                                                    & \alg{H}\tensor \alg{H} \arrow[rr, "{(1,S)}"] &                      & \alg{H}\tensor \alg{H} \arrow[ru, "\nabla"'] &  
\end{tikzcd}
\end{center}
 \end{definition}
\begin{remark}
    The antipode if it exists is uniquely determined by the bialgebra structure. 
\end{remark}

\begin{example}\label{ex:HopfAlgebraOfGroup}
    The motivating example of a Hopf algebra is the ring of finitely supported functions on a group, $\CoordRing(G) = \{f\colon G \rightarrow \R\mid f(g) \neq 0 \text{ for finitely many $G$}\}$. This has a natural multiplication inherited from $\R$ and a natural comultiplication induced by the multiplication on $G$ as follows:
    \[\Delta f (g_1, g_2) = f(g_1g_2).\]
    Similarly the antipode $S$ is induced by inversion in $G$, $(S f)(g) = f(g^{-1})$.
\end{example}

\begin{remark}
    In the previous example we implicitly use that fact that $\CoordRing(G\times G) \simeq \CoordRing(G)\tensor \CoordRing(G)$ to upgrade the dual group multiplication map from $\CoordRing(G) \rightarrow \CoordRing(G \times G)$ to $\CoordRing(G) \rightarrow \CoordRing(G) \tensor \CoordRing(G)$. We need $\R$ to be commutative for the map $\Phi \colon (\CoordRing(G)\tensor \CoordRing(G),\nabla_2) \rightarrow (\CoordRing(G\times G),\nabla_\R)$ given by $\Phi(f_1\tensor f_2)(g_1,g_2) = f_1(g_1)f_2(g_2) $ to preserve the multiplicative structure. On the one hand
    \[\Phi((\eta_1 \tensor \eta_2)(f_1\tensor f_2))(g_1,g_2) = \Phi(\eta_1f_1\tensor \eta_2f_2) = (\eta_1f_1)(g_1) \cdot (\eta_2f_2)(g_2) = \eta_1(g_1)f_1(g_1)\eta_2(g_2)f_2(g_2).\]
    On the other hand:
    \[\Phi(\eta_1\tensor \eta_1)(g_1,g_2)\cdot \Phi(f_1\tensor f_2)(g_1,g_2) = \eta_1(g_1)\eta_2(g_2) \cdot f_1(g_1) f_2(g_2).\]
\end{remark}

The heart of the problem is that the multiplication $\nabla_2$ on $\CoordRing(G)\tensor \CoordRing(G)$ doesn't preserve the order of the tensor products. Thus when we generalize to coordinate rings landing in an involutive ring $(\A,\sigma)$ we must replace $\alg{H}\tensor \alg{H}$ with a new object with the following properties:
\begin{itemize}
    \item Multiplication in the new object respects the ordering of tensor terms.
    \item When $\alg{H}$ is the ring of functions $G \rightarrow \A$ with $\A$  noncommutative, the new object should be isomorphic to the ring of functions $G\times G \rightarrow \A$.
\end{itemize}
We call this object the interleaved tensor algebra $T_2^*(\alg{H})$ and define it in \Cref{sec:KsigmaHopfAlgebra}.

\begin{remark}
    In the case of the usual quantum group, as Manin describes, the algebra $\alg{H}$ itself is non commutative, but as verified in property 5 of the quantum group in \cite{Manin-QuantumGroupsRemarks} it is necessary to assume some extra commutativity properties in order to verify a form of multiplicative structure on the quantum group. This commutativity is implied by the compatibility of the coproduct and product on $\alg{H}$.
\end{remark}

 \subsection{\texorpdfstring{($k,\sigma$)}{k sigma} Hopf Algebras}\label{sec:KsigmaHopfAlgebra}

Before we can define a $(k,\sigma)$-Hopf algebra we must define how the $\sigma$ structure on a $(k,\sigma)$-algebra extends to the tensor product.
Let $\alg{H}$ be a $(k,\sigma)$-algebra. Consider the $n-$fold tensor product over $k$, $T^{(n)}\alg{H} = \alg{H} \tensor \cdots \tensor \alg{H}$. 
\begin{lemma}
    $T^{(n)} \alg{H}$ is an $\alg{H}$ bimodule with the following actions and anti-involution $\sigma$.
    \begin{align*}
        a * (c_1\otimes \dots \otimes c_n) * b &= ac_1\otimes \dots \otimes c_nb \\
        \sigma((c_1\otimes \dots \otimes c_n)) &= \sigma(c_n)\otimes \dots \otimes \sigma(c_1) 
    \end{align*}
\end{lemma}
\begin{proof}
    It is clear that $\sigma$ is $k$ linear. So it suffices to check how $\sigma$ respects the left and right actions. We compute
    \begin{align*}
        \sigma(a * (c_1\otimes c_2 \otimes \dots \otimes c_{n-1} \otimes c_n) *b) =& \sigma(ac_1\otimes c_2 \otimes\dots \otimes c_{n-1}\otimes c_nb)\\
        =& \sigma(c_n b) \otimes \sigma(c_{n-1}) \otimes \dots \sigma(c_2) \otimes \sigma(a c_1)\\
        =& \sigma(b)\sigma(c_n) \otimes \sigma(c_{n-1}) \otimes \dots \otimes \sigma(c_2) \otimes \sigma(c_1)\sigma(a)\\
        =& \sigma(b) * \sigma(c_n)\otimes \sigma(c_{n-1})\otimes \dots \otimes \sigma(c_1) \otimes \sigma(c_1) * \sigma(a)\\
        =& \sigma(b) * \sigma(c_1\otimes \dots \otimes c_n) * \sigma(a)
    \end{align*}
\end{proof}

We now consider the full tensor algebra $T^*(\alg{H})$ and as an immediate corollary we get:
\begin{corollary}
    The tensor algebra $T^*(\alg{H})$ is a $(k,\sigma)$-algebra. 
\end{corollary}

\begin{remark}
    For any $n$ we have the induced multiplication map $T^{(n)}\alg{H} \rightarrow \alg{H}$ given by applying $\nabla$ $n-1$ times. This map provides an identification of $T^{(n)}\alg{H}$ with $\alg{H}$ as algebras.
\end{remark}

For our notion of a $(k,\sigma)$-Hopf algebra we will need to change the image of the coproduct from $\alg{H} \tensor \alg{H}$ to a slight modification of the tensor algebra that preserves the ordering of terms in the product while keeping the two ``sides'' of the coproduct separate.\\

Consider two $(k,\sigma)$-algebras $\alg{M}$ and $\alg{N}$. Then tensor algebra $T^*(\alg{M}) \tensor T^*(\alg{N})$ contains all possible interleaved tensor products of elements of  $\alg{M}$ and $\alg{N}$.  The multiplication map in $\alg{M}$ (respectively $\alg{N}$) that sends $\alg{M} \tensor \alg{M} \rightarrow \alg{M}$ via $m_1\tensor m_2 \rightarrow m_1m_2$ induces a map
\begin{align*}
    T^*(\alg{M})\tensor T^*(\alg{N}) \rightarrow& T^*(\alg{M}) \tensor T^*(\alg{N})\\
    s\tensor m_1\tensor m_2 \tensor t \mapsto & s\tensor m_1m_2 \tensor t
\end{align*}
\begin{definition}
    The \keyword{Interleaved Tensor algebra} $T_2^*(\alg{M},\alg{N})$ of is the quotient of $T^*(\alg{M})\tensor T^*(\alg{N})$ by the multiplication maps applied to repeated tensor factors in $\alg{M}$ and $\alg{N}$. For a single $(k,\sigma)$-algebra $\alg{H}$ we write $T_2^*(\alg{H})$ to mean $T_2^*(\alg{H}',\alg{H}'')$ where $\alg{H}'$ and $\alg{H}''$ are two formal copies of $\alg{H}$. 
\end{definition}

In other words every tensor in $T_2^*(\alg{H})$ is represented by $h_1\tensor h_2\tensor \dots \tensor h_n$ with $h_i \in \alg{H}'$ for $i$ odd and $h_i \in \alg{H}''$ for $i$ even (or vice versa).

\begin{example}
    In $T_2^*(\alg{H})$ $x'\tensor y'' \tensor (y'')^{-1} \tensor w' = x'w'$
\end{example}
The interleaved tensor algebra is functorial. A map $f\colon \alg{H}_1\rightarrow \alg{H}_2$ induces a map $T_2^*(f) \colon T_2^*(\alg{H}_1) \rightarrow T_2^*(\alg{H}_2)$ by applying $f$ to each factor. Similarly given two maps $f': \alg{H}\rightarrow \alg{H} $ and $f'' \colon \alg{H}\rightarrow \alg{H}$ there is a natural map $f' \otimes f'' \colon T_2^*(\alg{H}) \rightarrow T_2^*(\alg{H})$ given by applying $f'$ to $\alg{H}'$ factors and $f''$ to $\alg{H}''$ factors.
\begin{remark}
    Our construction of this tensor algebra is essentially realizing the ring of functions $\CoordRing(\alg{H}_1,\alg{H}_2)$
with values in $\A$ as $T_2^*(\alg{H})$. Since our functions are into a non commutative ring, we no longer have the classical fact that $\CoordRing(\alg{H}_1 \times \alg{H}_2) \simeq \CoordRing(\alg{H}_1) \otimes \CoordRing(\alg{H}_2)$
    
\end{remark}
Note that we also identify elements $A' \otimes c = cA'$ for $c\in k$.\\

We denote by $\tau$ the map $T_2^*(\alg{H}) \rightarrow T_2^*(\alg{H})$ the map swapping primes and double primes. This is the analogue of the map reversing the tensor factors in $\alg{H} \tensor \alg{H}$. However the order in $T_2^*(\alg{H})$ is preserved. For example, $\tau(a_1' \tensor a_2'' \tensor a_3') = a_1'' \tensor a_2' \tensor a_3''$.

Now we are ready to present our definition of a $(k,\sigma)$-Hopf algebra. The important take away from this definition is the interplay with the usual maps of a Hopf algebra and the map $\sigma$. We add the adjective $\sigma-$commutative to this definition to emphasise that this is a generalization of a commutative Hopf algebra. 

\begin{definition}
    A $\sigma-$commutative $(k,\sigma)$-Hopf algebra is a $(k,\sigma)$-algebra $\alg{H}$ along with 
    \begin{enumerate}
        \item (coproduct) a $(k,\sigma)$-algebra map $\Delta: \alg{H} \to T_2^*(\alg{H})$
        \item (counit) a $(k,\sigma)$-algebra map $\epsilon:\alg{H} \to k$
        \item (antipode) a $(k,\sigma)$-algebra map $S: \alg{H} \to \alg{H}$
    \end{enumerate}
    such that the following diagrams commute
\begin{center}
    \begin{tikzcd}
         & T_2^*(\alg{H}) \arrow[rr, "{(S,1)}"] &                      & T_2^*(\alg{H}) \arrow[rd, "\nabla"]  &   \\
\alg{H} \arrow[rr, "\epsilon"] \arrow[rd, "\Delta"'] \arrow[ru, "\Delta"] &                                & k \arrow[rr, "\eta"] &                                & \alg{H} \\
        & T_2^*(\alg{H}) \arrow[rr, "{(1,S)}"] &                      & T_2^*(\alg{H}) \arrow[ru, "\nabla"'] &  
\end{tikzcd}
\hspace{2pc}
\begin{tikzcd}
\alg{H} \arrow[rr, "S"] \arrow[d, "\Delta"']  &  & \alg{H} \arrow[d, "\Delta"] \\
T_2^*(\alg{H}) \arrow[rr, "\tau \circ S\otimes S"] &  & T_2^*(\alg{H})          
\end{tikzcd}
\end{center}

\end{definition}

   The fact that $\Delta$ is a $(k,\sigma)$-algebra map is encoded in the following diagrams commuting:
     \begin{center}
             \begin{tikzcd}
T_2^*(\alg{H}) \arrow[r, "\nabla"] \arrow[d, "T_2^*(\Delta)"']  & \alg{H} \arrow[d, "\Delta"]    \\
T_2^*(T_2^*(\alg{H})) \arrow[r, "T_2^*(\nabla)"]               &  T_2^*(\alg{H})       
\end{tikzcd}\hspace{2pc}
\begin{tikzcd}
         & k \arrow[rd, "T_2^*(\eta)"] \arrow[ld, "\eta"'] &                        \\
 \alg{H}\arrow[rr, "\Delta"] &                                                 & T_2^*(\alg{H})
\end{tikzcd}\hspace{2pc}
\begin{tikzcd}
    \alg{H} \arrow[d,"\Delta"']\arrow[r,"\sigma"] & \alg{H} \arrow[d,"\Delta"]\\
    T_2^*(\alg{H})\arrow[r,"\sigma"] & T_2^*(h)
\end{tikzcd}
        \end{center}
     Similarly the fact that the counit is a $(k,\sigma)$-algebra map is encoded as:
        \begin{center}
             \begin{tikzcd}
T_2^*(\alg{H}) \arrow[d, "T_2^*(\epsilon)"'] \arrow[r, "\nabla"] &   \alg{H} \arrow[d, "\epsilon"] \\
 T_2^*(k) \arrow[r,"\nabla_0"] & k &                         
\end{tikzcd}\hspace{2pc}
           \begin{tikzcd}
           & k\arrow[dr,"\id"]\arrow[dl,"\eta"'] & \\
            \alg{H}\arrow[rr,"\epsilon"] &   & k\\
       \end{tikzcd}\hspace{2pc}
\begin{tikzcd}
    \alg{H} \arrow[d,"\epsilon"']\arrow[r,"\sigma"] & \alg{H} \arrow[d,"\epsilon"]\\
    k\arrow[r,"\sigma"] & k
\end{tikzcd}
        \end{center}
These diagrams are analogous to the diagrams for a standard Hopf algebra with $T_2^*(\alg{H})$ replacing $\alg{H}\tensor \alg{H}$.

\begin{remark}
    The antipode is defined to be an algebra homomorphism not an anti-homomorphism, as opposed to the usual Hopf algebra definition. In the case of a commutative Hopf algebra, this is not an important distinction. Our requirement that $\sigma$ commutes with $S$ means that the composition $S\circ \sigma$ can be considered to be the analogue to the antipode in the usual Hopf algebra.
\end{remark}

\begin{theorem}
    The algebra $\CoordSPA_q$ is a $(k_q,\sigma)$-Hopf algebra with maps defined on generators:  
    \begin{enumerate}
        \item $\Delta\left(\begin{bmatrix} A & B\\ C & D\end{bmatrix}\right) = \begin{bmatrix} A' & B'\\ C' & D'\end{bmatrix} \otimes \begin{bmatrix} A'' & B''\\ C'' & D''\end{bmatrix} $ 
        \item $\epsilon \left( \begin{bmatrix} A & B\\ C & D\end{bmatrix} \right) = \begin{bmatrix} 1 & 0\\ 0 & 1\end{bmatrix}$
        \item  $S\left( \begin{bmatrix} A & B\\ C & D\end{bmatrix}\right) = \begin{bmatrix} \hat{D} & -q\hat{B}\\ -q^{-1}\hat{C} & \hat{A} \end{bmatrix}$
    \end{enumerate}
\end{theorem}
\begin{proof}
   First we must check these maps are well defined as $\CoordSPA_q$ is a quotient algebra. For example we must check that $\Delta(\hat{A}C- q^{-1}\hat{C}A)=0$:
   \begin{align*}
       \Delta&(\hat{A}C- q^{-1}\hat{C}A) \\ 
       & = \hat{A}''\otimes\hat{A}'C' \otimes A'' + \hat{A}''\otimes\hat{A}'D' \otimes C'' + \hat{C}'' \otimes \hat{B}'C' \otimes A'' +  \hat{C}'' \otimes \hat{B}'D' \otimes C'' \\ -& q^{-1}  \left( \hat{A}''\otimes\hat{C}'A' \otimes A'' + \hat{A}''\otimes\hat{C}'B' \otimes C'' + \hat{C}'' \otimes \hat{D}'A' \otimes A'' +  \hat{C}'' \otimes \hat{D}'B' \otimes C'' \right) \\
       &= \hat{A}''\otimes( \hat{A}'D' - q^{-1}\hat{C}'B' ) \otimes C'' + \hat{C}'' \otimes ( \hat{B}'C' - q^{-1} \hat{D}'A') \otimes A'' \\
       &= \hat{A}''C'' - q^{-1}\hat{C}''A''\\
       &=0
   \end{align*}
   The calculations for the other relations are similar.\\
   Similarly we check the antipode respects the relations. For example:
   \begin{align*}
       0 & = S(A\hat{B}-q^{-1}B\hat{A}) \\
       & = -\hat{D}q^{-1}B +\hat{B}D =0  
   \end{align*}
   The compatibility of multiplication and comultiplication is trivial as this is used to comultiplication extend beyond the generators. The compatibility between the unit/counit and other maps is also clear.\\
   Next we must check the compatibility of the antipode on each generator. We show the computation for the top path on the generator $\hat{A}$ and omit the remaining checks for brevity.
   \begin{align*}
      \nabla(S\otimes 1 (\Delta(\hat{A}))) =& \nabla (S\otimes 1\left( \hat{A}'' \otimes \hat{A}' + \hat{C}''\otimes \hat{B}'\right)) \\
       &= D\hat{A} - qC\hat{B} \\
       &=1 = \eta(\epsilon(\hat{A}))
   \end{align*}
   Finally we check the compatibility of the antipode and coproduct 
   \begin{align*}
       \Delta(S(A)) &=  \Delta(\hat{D}) = \hat{D}''\otimes \hat{D}' + \hat{B}''\otimes \hat{C}' \\ 
       &= \tau(\hat{D}'\otimes \hat{D}'' + \hat{B}'\otimes \hat{C}'') = \tau( (S \otimes S) (A' \otimes A'' + B' \otimes C'')) = \tau(S\otimes S(\Delta(A)))
   \end{align*}
\end{proof}

There is a coaction of $\CoordSPA_q$ on the algebra of the quantum plane mimicking the group action of the quantum group on the quantum plane. The proof of the next proposition is clear.
\begin{proposition}\label{prop:coaction}
    There is a $(k_q,\sigma)$-map $\CoordRing(\A^2_\omega)_q \to T_2^*(\CoordSPA_q,\CoordRing(\A^2_\omega)_q) $ defined by 
    $$ \delta\left(\begin{bmatrix} X \\ Y \end{bmatrix}\right) = \begin{bmatrix} A & B\\ C & D\end{bmatrix} \otimes \begin{bmatrix} X\\ Y \end{bmatrix} $$
\end{proposition}

Now assume that $\A$ is any $(k,\sigma)$-algebra.
Given two $(k,\sigma)$-algebra morphisms, $f,g: \alg{H} \to \A$ we can define their convolution product $f*g$ by the following composition:
$$ f*g:\begin{tikzcd}
\alg{H} \arrow[r, "\Delta"] & T_2^*(\alg{H}) \arrow[r, "{(f\otimes g)}"] & \A
\end{tikzcd} $$
\begin{proposition}
    The group of algebra morphisms $\CoordSPA_q \to \A_q$ with convolution product is isomorphic to the group $\SPA_q$
\end{proposition}
\begin{proof}
    Recalling \Cref{rem:points_of_coordinate_ring}, such an algebra morphism is an $\A_q$ point of $\CoordSPA_q$ i.e. a collection of 4 elements of $\A_q$ satisfying the quantum symplectic relations. We see that, by design, the convolution product gives matrix multiplication. It is easily checked that the inverse is given by precomposition with the antipode and the unit of the group is the counit of the Hopf algebra. 
\end{proof}

\begin{remark}
    The $k-$algebra generated by these algebra morphisms is the group ring $k[\SPA_q]$. The map $\sigma$ on $\A_q$ does not extend to a map on $k[\SPA_q]$. In fact the composition $\sigma \circ f: k[\SPA] \to \A_q $ is no longer an algebra morphism to $\A_q$ but is rather an algebra morphism to $\A_q^{op}$. This explains part of our difficulty in having $\sigma(f)$ be an element of the symplectic group.
\end{remark}

\subsection{Comparison with \texorpdfstring{$SL_2$}{SL2}}

We can now compare $\CoordSPA_q$ with $\CoordRing(\SL_2)_q$. 
\begin{lemma}
    The algebra of quantum matrix functions, $\CoordRing(\SL_2)_q$ has a structure of a $(k_q,\sigma_q)$ algebra with anti-involution $\widetilde{\sigma}$ induced by taking $\widetilde{\sigma}$ to be the identity on the generators $a,b,c,d$.
\end{lemma}
\begin{proof}
    First note that $\widetilde{\sigma}$ extends $\sigma_q : k_q \rightarrow k_q$ and so takes $q$ to $q^{-1}$. Then it suffices to check that $\widetilde{\sigma}$ vanishes on the quantum $\SL_2$ relations. We check one case here as the others are similar.
    \[\widetilde{\sigma}(ac-q^{-1}ca) = \widetilde{\sigma}(c)\widetilde{\sigma}(a) - q \widetilde{\sigma}(a)\widetilde{\sigma}(c) = ca-qac = -q(-q^{-1}ca+ac) = 0 \]
\end{proof}
We then can define a surjective $(k_q,\sigma_q)$ map,  $\Phi \colon \CoordSPA_q \rightarrow \CoordRing(\SL_2)_q$ induced on generators by \[\Phi(A) = a, \Phi(B) = b, \Phi(C) = c, \Phi(D)=d.\] This map is obviously surjective as it hits every generator. However, to check it is well defined we must verify that the image of the quantum symplectic equations vanish. Since $\Phi$ is a $(k_q,\sigma_q)$ map we have that $\Phi(\hat{A}) = a$, $\Phi(\hat{B}) = b$, $\Phi(\hat{C}) =c$ and $\Phi(\hat{D}) = d$. Furthermore the only difference between the quantum $\SL_2$ equations and the quantum symplectic equations is the hats. Thus $\Phi$ directly maps each quantum symplectic equation onto a quantum $\SL_2$ equation and is well defined.\\
\begin{lemma}
    The kernel of $\Phi$ is the ideal $\hat{S}$ generated by $A-\hat{A},B-\hat{B},C-\hat{C},D-\hat{D}$.
\end{lemma}
\begin{proof}
    It is clear that $\hat{S} \subseteq \ker(\Phi)$. To show the other inclusion it suffices to show that for each quantum $\SL_2$ relation the preimage under $\Phi$ is contained in $S$. Each such preimage corresponds to one of the quantum symplectic relations although with potentially the wrong choice of hat or not. The choice of hats in each term can be modified using terms from $\hat{S}$ to reduce the preimage to the corresponding quantum symplectic relation completing the proof. For example consider $\hat{A}D-q^{-1}\hat{B}C -1 \in \Phi^{-1}(ad-q^{-1}bc-1)$. Then
    \begin{align*}
        \hat{A}D-q^{-1}\hat{B}C-1 =& \hat{A}D-AD+AD-A\hat{D}+A\hat{D}-q^{-1}\hat{B}C+q^{-1}BC-q^{-1}BC+q^{-1}B\hat{C}-q^{-1}B\hat{C}-1\\
        =&(\hat{A}-A)D+A(D-\hat{D})+A\hat{D}-q^{-1}(\hat{B}-B)C-q^{-1}B(C-\hat{C})-q^{-1}B\hat{C}-1\\
        =& (\hat{A}-A)D+A(D-\hat{D})-q^{-1}(\hat{B}-B)C-q^{-1}B(C-\hat{C})+A\hat{D}-q^{-1}B\hat{C}-1
    \end{align*}
    This is clearly in $\hat{S}$ up to an element of the quantum symplectic relations and so is in $\hat{S}$ in $\CoordSPA$.
\end{proof}
Therefore $\CoordRing(\SL_2)_q$ is isomorphic to $\CoordSPA_q/\hat{S}$. This leads to the following theorem:

\begin{theorem}
    The $(\A_q,\sigma_q)$ points of $\CoordRing(\SL_2)_q$ are exactly the $(\A_q,\sigma_q)$ points of $\CoordSPA$ with $\sigma_q$ fixed entries.
\end{theorem}
\begin{proof}
    The $(\A_q,\sigma_q)$ points of $\CoordSPA/\hat{S}$ correspond to choices of four values $A,B,C,D\in \A_q$ that satisfy the quantum symplectic relations and are fixed by $\sigma_q$. Under the isomorphism $\Phi$ these are exactly the $(\A_q,\sigma_q)$ points of quantum $\SL_2$.
\end{proof}
Note that even when the original $\sigma$ is trivial and $\A$ commutes the condition that the matrix entries are $\sigma_q$ fixed points is nontrivial. In fact this property isn't preserved under matrix multiplication, which explains part of the difficulty for $\SL_2(\A)_q$ to be a nontrivial group. 
\begin{example}
    The elements of $\SPA_q$ constructed in \Cref{thm:SP2A_elements} of the form $\begin{bmatrix}
        A & q^{1/2}B\\ q^{-1/2}C & D
    \end{bmatrix}$ for $\begin{bmatrix}
        A & B\\ C &D
    \end{bmatrix} \in \SPA$ are $(\A_q,\sigma_q)$ points of $\CoordSPA_q$ but not of $\CoordRing(\SL_2)_q$. This is clear as $q^{1/2}B$ and $q^{-1/2}C$ aren't fixed by $\sigma_q$ even if $\sigma(B) = B$ and $\sigma(C) = C$.
\end{example}

The Hopf algebra structure we have constructed on $\CoordRing_q$ descends to the
 usual bialgebra structure on $\CoordRing(\SL_2)$ after a few considerations. 

To obtain $\CoordRing(\SL_2)_q$ we imposed that $\sigma$ acted trivially on the generators. Now to obtain the usual bialgebra structure we will also impose that $\sigma$ acts trivially on the generators of $T_2^*(\CoordRing(\SL_2)_q)$ 

This means that in $T_2^*(\CoordRing(\SL_2)_q)$ we have that for all $x,y\in \{a,b,c,d\}$ 
$$x'\otimes y'' = \sigma(x'\otimes y'') = y''\otimes x'.$$

It follows that the image of $\Delta$ in $T_2^*(\CoordRing(\SL_2)_q)$ becomes simply $\CoordRing(\SL_2)_q \otimes \CoordRing(\SL_2)_q'$ with the product structure given simply by $$(x_1' \otimes y_1'')(x_2' \otimes y_2'') = (x_1'x_2'\otimes y_1''y_2'') $$

Now the coproduct defined on $\CoordSPA_q$ descends to the usual coproduct on $\CoordRing(\SL_2)_q$ and the following diagram commutes: 
\begin{center}
    \begin{tikzcd}
\alg{H}\otimes \alg{H} \arrow[r, "m"] \arrow[d, "\Delta \otimes \Delta" description] & \alg{H} \arrow[r, "\Delta"] & \alg{H}' \otimes \alg{H}''                                                          \\
\alg{H}'\otimes \alg{H}'' \otimes \alg{H}'\otimes \alg{H}'' \arrow[rr, "1 \otimes \sigma \otimes 1"] &                       & \alg{H}' \otimes \alg{H}' \otimes \alg{H}'' \otimes \alg{H}'' \arrow[u, "m\otimes m" description]
\end{tikzcd}
\end{center}

The antipode map $S$ defined earlier, however is not well defined on the quotient, explaining some of the difference between our Hopf algebra and the usual Hopf algebra on the quantum group.

\section{Cluster Structure}\label{sec:Cluster Structure}

In the commutative setting the work of \cite{BFZ-ClusterAlgebrasIII} give a cluster algebra structure on the large open cell in $\SL_2(\R)$. Using the work of \cite{BR-noncommutativeClusters} we show there is a corresponding noncommutative cluster structure on $\SPA$.

Explicitly we produce a noncommutative cluster structure on the ring $\B = \CoordSPA[A^{-1},B^{-1},C^{-1},D^{-1}]$. The $(\A,\sigma)$-points of $\B$ are the elements of $\SPA$ with invertible entries.\\

From this cluster structure we obtain functions $\Phi_L,\Phi_R$ which up to multiplication by $B$ and $C$ respectively are equivalent to the standard trace in $\SL_2$. We also interpret matrix multiplication, transpose, and inverse on $\SPA$ as operations on the surfaces controlling the cluster structure.\\
First we review the cluster structure on triangulations of a surface in both the commutative and noncommutative case. This general theory will be used in \Cref{sec:NoncomMarkov} to produce noncommutative generalizations of Markov numbers.

\subsection{Commutative Cluster Structure}\label{sec:CommutativeClusterStructure}
First we review the classic cluster structure associated to a triangulated surface. A cluster algebra is an algebra with a rich additional combinatorial structure. Informally, the algebra can be generated from an initial `seed' via a process called `mutation' which produces new seeds from old seeds. The family of seeds is called cluster structure of the algebra.\\
Each seed consists of a set of cluster variables and some combinatorial information used to describe how the cluster variables mutate. In full generality this extra information is encoded by a quiver or exchange matrix (see \cite{FWZ-IntroToClusters} for full details). However in this paper we are only concerned with cluster algebras associated to surfaces and thus will define a seed simply using a  triangulation of a surface.
\begin{definition}
    A \keyword{bordered surface} $(S,M)$ is a surface $S$ with boundary components or punctures along with a finite set of \keyword{marked points} $M\subset S$ which contains at least one element on each component of $\boundary(S)$ and all punctures. We require that $M$ is non-empty. 
\end{definition}

\begin{definition}
    An \keyword{arc} on a bordered surface is a path on $S$ with endpoints at marked points, considered up to homotopy fixing the marked points. Paths which are homotopic to a puncture are not considered arcs. An arc homotopic to a boundary component of $S$ is called a \keyword{boundary arc}. 
    A \keyword{triangulation}, $T$ of a bordered surface is a maximal collection of noncrossing arcs. 
\end{definition}
\begin{remark}
    If $S$ is a surface with punctures other than a closed surface with one puncture the cluster structure actually requires a ``tagged triangulation'' (\cite{FST-ClusterTriangulatedSurface}). We won't need this generality so we omit this technical definition. 
\end{remark}
\begin{definition}
    A \keyword{seed} for the cluster structure of the surface $(S,M)$ consists of a triangulation $T$ along with a choice of \keyword{cluster variable} $x_\gamma$ for each arc $\gamma$ of $T$. The variables associated to boundary arcs are called \keyword{frozen} while the variables associated to interior arcs are \keyword{unfrozen} or \keyword{mutable}.
\end{definition}
When the endpoints $\{i,j\} \subset M$ of $\gamma$ uniquely determine it as an arc we write $x_{i,j}$ for $x_{\gamma}$.
\begin{definition}
    Each seed $(T,\vec{x})$ can be mutated at an interior arc $\gamma$ to produce a new seed $\mu_\gamma(T,\vec{x}) = (T',\vec{x'})$. The triangulation $T'$ is obtained by removing $\gamma$ and replacing it with the arc $\gamma'$ corresponding to the opposite quadrilateral containing $\gamma$ (\Cref{fig:QuadrilateralFlip}). We name the edges of the bounding quadrilateral cyclically by $\alpha,\beta,\eta,\delta$. The cluster variables in $\vec{x'}$ are the same as the variables of $\vec{x}$ except $x_\gamma$ is replaced by $x_{\gamma'}$ which satisfies the \keyword{exchange relation}:
    \begin{equation}
        x_\gamma x_{\gamma'} = x_\alpha x_\eta + x_\beta x_\delta
    \end{equation}
\end{definition}
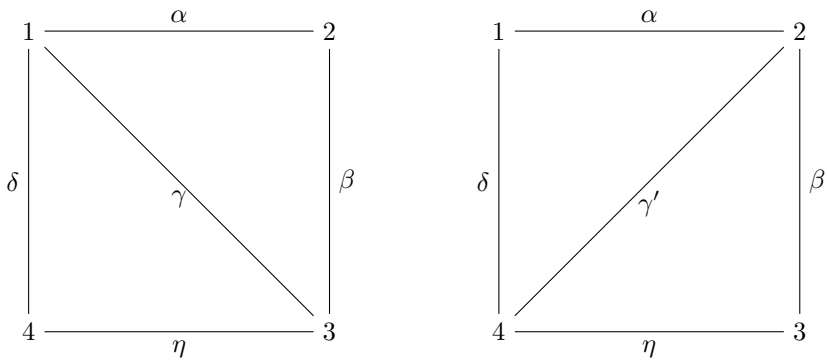
\begin{figure}[hb]
    \centering
    \begin{tikzpicture}
        \node[] at (4,4) (2) {2};
        \node[] at (0,0) (4) {4};
        \node[] at (0,4) (1) {1};
        \node[] at (4,0) (3) {3};

        %Square
        \path[] (4) edge [-] node[left] {$\delta$}  (1);
        \path[] (1) edge [-] node[below] {$\gamma$} (3);
        \path[] (1) edge [-] node[above] {$\alpha$} (2);
        \path[] (2) edge [-] node[right] {$\beta$} (3);
        \path[] (3) edge [-] node[below] {$\eta$} (4); 
    \end{tikzpicture}\hspace{3pc}
    \begin{tikzpicture}
        \node[] at (4,4) (2) {2};
        \node[] at (0,0) (4) {4};
        \node[] at (0,4) (1) {1};
        \node[] at (4,0) (3) {3};

        %Square
        \path[] (4) edge [-] node[left] {$\delta$}  (1);
        \path[] (2) edge [-] node[below] {$\gamma'$} (4);
        \path[] (1) edge [-] node[above] {$\alpha$} (2);
        \path[] (2) edge [-] node[right] {$\beta$} (3);
        \path[] (3) edge [-] node[below] {$\eta$} (4); 
    \end{tikzpicture}
    \caption{Mutation at the arc $\gamma$.}
    \label{fig:QuadrilateralFlip}
\end{figure}
\begin{definition}
    The \keyword{cluster algebra} associated to the surface $(S,M)$ with initial seed $(T,x_1, \dots x_{|T|})$ is the subalgebra of the ring of rational functions in the variables of the an initial seed, $\mathbb{Q}(x_1,\dots, x_{|T|})$, generated by all cluster variables obtained by repeated mutation from the initial seed.
\end{definition}

The cluster algebra associated to a triangulated surface can be geometrically realized as the coordinate ring of decorated Teichm\"uller space. 
\begin{definition}
    The \keyword{decorated Teichm\"uller space} associated to a bordered surface $(S,M)$ is the moduli space of hyperbolic metrics on $S$ with geodesic boundary and cusps at the marked points, along with with a decoration consisting of a choice of horocycle around each cusp. 
\end{definition}
Under this isomorphism the cluster coordinate $a_\gamma$ associated to an arc $\gamma$ is the function on the Teichmüller space computing the $\lambda-$length of the geodesic in the homotopy class of $\gamma$ for a given hyperbolic metric. (See \cite{FT-LambdaLengths} for full details). Importantly, \cite{FG-DualTeichmuller} find that the decorated Teichmüller space is parameterized by the positive $\R$-points of this cluster algebra\footnote{Since we decorate the space, we describe the type $A$ cluster structure of \cite{FG-DualTeichmuller}. }.

The key result for us is the following lemma which relates the length of the portion of a horocycle inside a triangle to the cluster coordinates of its edges.
\begin{lemma}\label{thm:MagicHorocycleLength}
    The length of the portion of the horocycle at vertex $i$ of the triangle $ijk$ is \[t_i^{jk} = \frac{x_{jk}}{x_{ij}x_{ki}}\]
\end{lemma}
\begin{proof}
    The proof is an explicit calculation in hyperbolic geometry. See Lemma 7.9 of \cite{FT-LambdaLengths}.
\end{proof}
\begin{remark}
    The exchange relation for $\lambda-$lengths is equivalent to the fact that lengths of horocycles are additive. 
    \begin{align*}
        t_i^{j\ell} =& t_i^{jk} + t_i^{k\ell}\\
        \frac{x_{j\ell}}{x_{ij}x_{i \ell }} = & \frac{x_{jk}}{x_{ij}x_{ik}} + \frac{x_{k\ell}}{x_{ik}x_{i\ell}}\\
        x_{j\ell}x_{ik} =& x_{jk}x_{i\ell} + x_{ij}x_{k\ell}
    \end{align*}
\end{remark}

The choice of a positive real number for each cluster variable is equivalent to a point in the decorated Teichmüller space which in turn is equivalent to a decorated representation $\pi_1(S) \rightarrow \PSL_2(\R)$ up to conjugation.
\begin{definition}
    A \keyword{$\PSL_2(\R)$ decorated representation} of a bordered surface $(S,M)$ is a homomorphism $\rho \colon \pi_1(S) \rightarrow \PSL_2(\R)$ along with a $\rho$ equivariant choice of pairwise transverse vectors $\vec{v}_i \in \R^2$ for each marked point $i \in M$. In particular if $\gamma$ is a loop homotopic to a puncture $i$, we require that $\rho(\gamma)\vec{v}_i = \vec{v}_i$.
\end{definition}
For our study, the important detail of the correspondence is the construction of a representation from the cluster data, a choice of real number $a_\gamma$ for each cluster variable $x_\gamma$ in a fixed triangulation. We sketch the construction here and refer to Section 4.2 of \cite{FG-DualTeichmuller} for full detail. Given a loop $\gamma\in \pi_1(S)$ we build the image by taking the product of prescribed matrices for each triangle $\gamma$ passes through. Concretely we draw a hexagon inside each triangle with ``long'' sides along the edges and ``short'' sides cutting off each angle. Orient the edges of the hexagon clockwise with respect to the orientation of the triangle. We assign a matrix to each edge of the hexagon by assigning each long edge a matrix $S_{ij}$ and each short edge a matrix $G_i^{kj}$ (See \Cref{fig:RepInTriangle}). 
\[G_{i}^{kj} = \begin{bmatrix}1 & t_i^{kj}\\ 0 & 1 \end{bmatrix} = \begin{bmatrix}1 & \frac{a_{kj}}{a_{ik} a_{ij}}\\ 0 & 1 \end{bmatrix} \hspace{3pc} S_{ij} = \begin{bmatrix}
    0 & -a_{ij}^{-1}\\ a_{ij} & 0
\end{bmatrix}\]
Then the image of $\gamma$ is given by deforming $\gamma$ onto the edges of the hexagons and taking the corresponding product. Note that when following an edge against its orientation we take the inverse matrix, $G_i^{kj} = (G_i^{jk})^{-1}$.
\begin{remark}
    We use a slightly different convention than in \cite{FG-DualTeichmuller}. They use an opposite orientation of the hexagon and use the transpose of each matrix.
\end{remark}
\begin{figure}[hb]
    \centering
    \begin{tikzpicture}
        \node[] at (4,4) (3) {3};
        \node[] at (0,0) (1) {1};
        \node[] at (0,4) (2) {2};
        \path[] (3) edge [->] node[below] {}  (1);
        \path[] (2) edge [->] node[below] {} (3);
        \path[] (1) edge [->] node[left] {} (2);

        \path[red,thick] (0,1/1.14) edge [<-] node[above] {$G_1^{32}$} (1/1.14,1/1.14);
        \path[red,thick] (1/1.14,1/1.14) edge [<-] node[left] {$S_{31}$} (4-1/1.14,4-1/1.14);
        \path[red,thick] (4-1/1.14,4-1/1.14) edge [<-] node[left] {$G_3^{21}$} (4-1/1.14,4);
        \path[red,thick] (4-1/1.14,4) edge [<-] node[below] {$S_{23}$} (1/1.14,4);
        \path[red,thick] (1/1.14, 4) edge [<-] node[below] {$G_2^{13}$} (0,4-1/1.14);
        \path[red,thick] (0, 4-1/1.14) edge [<-] node[right] {$S_{12}$} (0,1);
    \end{tikzpicture}
    \caption{Hexagon inside generic triangle.}
    \label{fig:RepInTriangle}
\end{figure}
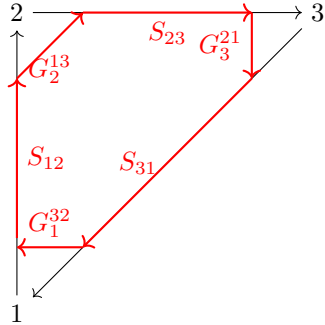
\begin{remark}
    Under this definition $S_{12}S_{21} = \begin{bmatrix}
        -1 & 0\\0 & -1
    \end{bmatrix} = G_1^{32}S_{12}G_2^{13}S_{23}G_3^{21}S_{31}$. As such the representation is well defined into $\PSL_2(\R)$, but not into $\SL_2(\R)$. To handle this sign one can define a twisted local system as in \cite{FG-DualTeichmuller}.
\end{remark}

The inverse construction is easier. Given a decorated representation $(\rho,\vec{v}_i)$ and a triangulation of the surface we assign the coordinate $a_{ij} = |\vec{v}_i \wedge \vec{v}_j |$ to the edge $i$ --- $j$.

\subsection{Noncommutative Surface Cluster Structure}\label{sec:surfaceClusterStructure}
  The noncommutative surface cluster structure from \cite{BR-noncommutativeClusters} produces a noncommutative generalization of the cluster algebra associated to a surface as defined above.

  % on the functions on the space of decorated representations $\pi_1(S) \rightarrow \PSP_2(\A,\sigma)$. We note that if $\A = \Mat(n,\R)$ and $\sigma$ is transpose, this construction is the symplectic noncommutative cluster algebras of \cite{AGRW-SymplecticClusters}. 
  
  We now recall the relevant definitions.
 \begin{definition}
     A \keyword{seed} for a noncommutative cluster structure consists of an ideal triangulation of a bordered surface. Each possible orientation of an edge $i \rightarrow j$ of the triangulation is assigned a variable $X_{ij}$ in a noncommutative algebra. 
     There is an additional condition for each triangle cyclically labeled $i,j,k$. We define a \keyword{noncommutative angle} $T_i^{jk} = X_{ji}^{-1}X_{jk}X_{ik}^{-1}$ and we require that $T_i^{jk} = T_i^{kj}$. We call this equation a \keyword{triangle condition}.

      Again, we grow the noncommutative cluster algebra from mutation. Each mutable edge $i \rightarrow k$ is contained in a unique quadrilateral with vertices $(i,j,k,\ell)$. The mutated seed contains the triangulation with $i\rightarrow k$ replaced with $j \rightarrow \ell$. The element assigned to $j\rightarrow \ell$ is given by the formula 
 \begin{equation}\label{eqn:NoncomMutation}
     X_{j\ell} = X_{jk}X_{ik}^{-1}X_{i\ell} + X_{ji}X_{ki}^{-1}X_{k\ell}
 \end{equation}

 The noncommutative surface cluster algebra, $C(S)$ associated to an initial seed $(T,X_1, \dots X_{|T|})$ is the subalgebra of the noncommutative fraction field of the $\mathbb{Q}(\{X_i\})$ generated by all of the variables $X_{\gamma}$ found in all seeds modulo the triangle conditions. 
 \end{definition}

   This algebra becomes a $(\mathbb{Q},\sigma)$-algebra once we define that $\sigma(X_{ij})=X_{ji}$ and extend $\sigma$ as an anti-involution. We will usually use the map $\sigma$ to eliminate the need to specify both $X_{ij}$ and $X_{ji}$. 

  \begin{remark}

    Using $\sigma$, we see that the triangle condition says that $$T_i^{jk} = T_i^{kj} = \sigma(T_i^{jk}),$$ so that each noncommutative angle is in $C(S)^\sigma$ 
    
    This mutation rule is equivalent to the requirement that the angles are additive.
     \begin{align*}
         T_{i}^{j\ell} =~& T_{i}^{jk}+T_{i}^{k\ell}\\
         {X_{ji}^{-1}}X_{j\ell}{X_{i \ell}^{-1}} =~& {X_{ji}^{-1}}X_{jk}{X_{ik}^{-1}} + {X_{ki}^{-1}}X_{k\ell}{X_{i \ell}^{-1}}
     \end{align*}
 \end{remark}

\subsection{The $\A$ points of a Noncommutative Cluster Algebra}\label{sec:ApointOfClusterAlgebra}

It is a natural question to ask for the $\A$ valued points of a noncommutative cluster algebra, i.e. to parameterize the $(\Q,\sigma)$-maps $C(S) \to \A$. The generative nature of cluster algebras makes it even more natural to begin by evaluating an initial seed at elements of $\A$ and then use the mutation rule to evaluate the image of all other cluster variables. In order to do this we must make some considerations. 

Firstly, our choice of elements of $\A$ for the initial seed must satisfy the triangle conditions for each triangle in our triangulation of $S$. Secondly, in order to perform a mutation, we must have that each cluster variable is mapped to an invertible element of $\A$. 

This last point may be difficult to achieve; one may start with a collection of invertible elements of $\A$ and after several mutations end up with a non-invertible element. In the case that $\A = \R$ the real numbers, one may remedy this problem by restricting to initial collections of \emph{positive} real numbers. Since the mutation rule only involves multiplication and addition we can be assured that the new mutated cluster variables will also be positive and hence invertible. 

\begin{definition}
    We call the map $\A \to \A $ sending $A$ to $\sigma(B)AB$  the \keyword{$\sigma$-congruence} of $A$ by $B$ for $A,B \in \A$.
\end{definition}

\begin{definition}
    A \keyword{positive structure} on an involutive ring $(\A,\sigma)$ is a choice of a proper open convex cone, $A_+^\sigma$ in $(\A^\sigma)^\times$ which 
    \begin{enumerate}
        \item contains 1.
        \item is closed under inversion.
        \item is closed under $\sigma$-congruence by elements of $\A^\times$.
    \end{enumerate}
\end{definition}

The prototypical examples of a positive structure are the positive real numbers sitting inside $\R$ and the positive definite symmetric matrices sitting inside $\Mat_n(\R)$. 
\begin{remark}
    One subtlety of this definition is that the cone is open in $\A^\sigma$ but not necessarily in $\A$. For example the positive real numbers sitting inside $\R$ inside of $\C$ is a positive structure for $(\C,\sigma)$ with $\sigma$ complex conjugation. This example also generalizes to the positive definite Hermitian matrices sitting inside $\Mat_n(\C)^\sigma$ with conjugate transpose as the involution $\sigma$.
\end{remark}

\begin{definition}
    The \keyword{positive $\A-$points}  of a noncommutative cluster algebra with respect to a given positive structure, $\A_+^\sigma$ are the $\A-$points for which all of the angles lie in $\A_+^\sigma$.  
\end{definition}

\begin{theorem}\label{thm:PostiveStructure}
    Any choice of elements of $\A$ for the cluster variables of an initial seed of $C(S)$ which are invertible and satisfy that the angles are in $\A_+^\sigma$ parameterize a positive $\A-$point.
\end{theorem}
\begin{proof}
    We simply need to check that each mutation produces a new seed satisfying the conditions of the theorem, then by induction each seed will. Let $\{X_{ij}\}$ be the cluster variables in our initial seed and let $\{A_{ij}\}$ be the elements at $\A$ we are evaluating them at. 

    Consider a single mutation from the initial seed. Locally this happens in a square in the triangulation of $S$. Lets label the endpoints of this square $1,2,3,4$ and suppose that the variables $X_{12},X_{23},X_{34},X_{14},X_{13}$ are in our initial cluster. We would like to show that $$A_{24} = A_{21}A_{31}^{-1}A_{34}+A_{23}A_{13}^{-1}A_{14}$$ is invertible and that the new noncommutative angles are in $\A_+^\sigma$. 

    The mutation rule implies that the new angle $ T_1^{24}=A_{21}^{-1}A_{24}A_{14}^{-1} = T_1^{23}+T_1^{34} $ is a sum of two elements of $\A_+^\sigma$ and is thus an element of $\A_+^\sigma$. We also see that the angle \[T_2^{14} = A_{12}^{-1}A_{14}A_{24}^{-1} = \sigma(A_{21}^{-1})(A_{21}^{-1}A_{24}A_{14}^{-1})^{-1}(A_{21})^{-1}\] is obtained from an angle in $\A_+^\sigma$ by inversion and  $\sigma$-congruence. Similarly all other new angles can be found to be back in $\A_+^\sigma$. Finally we see that $A_{24}$ is invertible since $T_1^{24}, A_{21}, A_{14}$ all are. 

    Thus we conclude that after a mutation we again have a seed consisting of invertible elements and with all angles positive. 
\end{proof}

The total angle around a marked point, $i$, of $S$ is an important invariant of an $\A$ point of the cluster algebra. We write $\Phi_m = \sum_{(j,k)}T_m^{jk}$ for this sum where $(j,k)$ runs over all triangles with vertex $m$. For a positive $\A$ point, this is an element of $A_+^\sigma$. One might wonder which values of this total sum are possible to achieve for all possible positive $\A$ points.
\begin{proposition}
    Let $B \in \A^\times$. The operation sending $A_{ij}$ to $A_{ij}'$ induced by
    $$ A_{ij}' = \begin{cases}
        \sigma(B)^{-1}A_{mj} & \textit{if } i = m \\
        A_{ij} & \textit{else}
    \end{cases}$$ 
    for all $j$ sends the positive $\A$ point $\{A_{ij}\}$ to a new positive $\A$ point $\{A_{ij}'\}$ and sends $\Phi_m$ to $B\Phi_m\sigma(B)$ and fixes all other angles. 
\end{proposition}

Thus given an involutive ring  $(\A,\sigma)$, we can compute the positive $(\A,\sigma)$-points of a noncommutative surface cluster algebra.  This data is equivalent to the choice of a positive decorated representation $\pi_1(S) \rightarrow \PSP_2(\A,\sigma)$ up to conjugation.
 \begin{definition}
     A \keyword{$\PSPA$ decorated representation} of a bordered surface is a homomorphism $\rho \colon \pi_1(S) \rightarrow \PSPA$ along with a $\rho$ equivariant choice of pairwise transverse vectors in $\vec{F}_i \in \A^2_{\omega}$ for each marked point. In particular if $\gamma$ is a loop homotopic to a puncture $i$, we have that $\rho(\gamma) \vec{F}_i = \vec{F}_i$
 \end{definition}
 The correspondence follows the same procedure as for $\PSL_2(\R)$ after ``upgrading'' the formulas with the appropriate uses of $\sigma$.\\
Given such a representation the coordinate $A_{ij}$ assigned to the edge $i\rightarrow j$ connecting vectors $\vec{F}$, $\vec{G}$ is given (up to a sign) by $\omega(\vec{F},\vec{G})$. Explicitly if $\vec{F} = (F_1,F_1)$ and $\vec{G} = (G_1,G_2)$ then $A_{ij} = \pm (\sigmaT{F_1}G_2-\sigmaT{F_2}G_1)$. \\
The inverse construction builds a decorated representation from the data of a seed. Here we describe the construction of $\rho$ from the seed data. The only difference from the commutative case is the definition of the matrices $\matSP{G}_i^{jk}$ and $\matSP{S}_{ij}$ assigned to the edges of the hexagon as in \Cref{fig:RepInTriangle}. Here we take the following matrices
\[ G_i^{kj} = \begin{bmatrix}1 & T_i^{kj}\\0 & 1 \end{bmatrix} = \begin{bmatrix} 1 & A_{ki}^{-1}A_{kj}A_{ij}^{-1}\\ 0 & 1\end{bmatrix}\hspace{3pc} S_{ij} = \begin{bmatrix}
    0 & -\sigmaT{A_{ij}^{-1}}\\
    A_{ij} & 0
\end{bmatrix} \]

\begin{remark}
    When $(\A,\sigma)$ is $\Mat_n(\R)$ with transpose this is essentially the construction of  \cite{AGRW-SymplecticClusters}. In this case the positive points correspond to maximal representations, which are in particular discrete and faithful, see also \cite{R-FramedMaximalReps,KR-FramedLocalSystems}.
\end{remark}

\subsection{Interpreting Symplectic Matrices as Seeds}
\begin{figure}[hb]
    \centering
    \begin{tikzpicture}
        \node[] at (4,4) (2) {2};
        \node[] at (0,0) (4) {4};
        \node[] at (0,4) (1) {1};
        \node[] at (4,0) (3) {3};
        \node[] at (1/1.14,4-1/1.14) (Aback) {};
        \node[] at (4-1/1.14,1/1/1.14) (Afront) {};
        \node[] at (0,4-1/1.14) (Bfront) {};
        \node[] at (0,1/1.14) (Bback) {};
        \node[] at (4,4-1/1.14) (Cback) {};
        \node[] at (4,1/1.14) (Cfront) {};
        \node[circle,draw,fill=red] at (4-1/1.14,0) (Botback) {};
        \node[] at (1/1.14,0) (Botfront) {};
        \node[circle,draw,fill=red] at (4-1/1.14,4) (Topback) {};
        \node[] at (1/1.14,4) (Topfront) {};
        %Square
        \path[] (4) edge [<-] node[left] {B}  (1);
        \path[] (1) edge [->] node[below] {A} (3);
        \path[] (1) edge [->] node[above] {1} (2);
        \path[] (2) edge [->] node[right] {C} (3);
        \path[] (3) edge [->] node[below] {1} (4); 

        %Hexagon 124
        \path[red,thick] (Aback) edge [->] (Afront);
        \path[red,thick] (Bback) edge [->] (Bfront);
        \path[red,thick] (Botback) edge [->] (Botfront);
        \path[red,thick] (Afront) edge [->] (Botback);
        \path[red,thick] (Botfront) edge [->] (Bback);
        \path[red,thick] (Bfront) edge [->] (Aback);

        %Hexagon 123
        \path[red,thick] (Aback) edge [->] (Afront);
        \path[red,thick] (Cback) edge [->] (Cfront);
        \path[red,thick] (Topback) edge [<-] (Topfront);
        \path[red,thick] (Afront) edge [<-] (Cfront);
        \path[red,thick] (Cback) edge [<-] (Topback);
        \path[red,thick] (Topfront) edge [<-] (Aback);
    \end{tikzpicture}
    \caption{Seed associated to a symplectic matrix.}
    \label{fig:SymplecticSeed}
\end{figure}
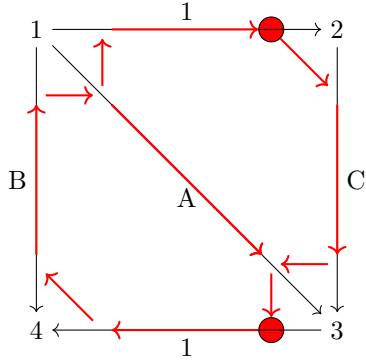

Given a matrix $\matSP{M} = \begin{bmatrix}A & B\\ C & D\end{bmatrix} \in \SPA$ we define a seed $S_\matSP{M}$ associated to a disk with 4 marked points. As the disk has trivial fundamental group the only data of a decorated representation is the choice of pairwise transverse vectors in $\A_\omega^2$. The seed $S_\matSP{M}$ corresponds to the decorated representation with the following choice of vectors:
\[\vec{F}_1 = \begin{bmatrix}0\\-1\end{bmatrix} \hspace{1pc} \vec{F}_2 = \begin{bmatrix}1 \\ 0 \end{bmatrix} \hspace{1pc} \vec{F}_3 = \begin{bmatrix}A \\ C \end{bmatrix}  \hspace{1pc} \vec{F}_4 = \begin{bmatrix}B \\ D \end{bmatrix}.\]
We take the triangulation with an edge from $1\rightarrow 3$. Then we compute the coordinate for each edge (See \Cref{fig:SymplecticSeed}):
\[A_{13} = A, A_{23} = C, A_{14} = B, A_{12} = A_{34} = 1.\]
The noncommutative angles are symmetric due to the fact that $\matSP{M}$ satisfies the symplectic conditions. For example $T_3^{12} = \sigmaT{A^{-1}}1C^{-1} = (C\sigmaT{A})^{-1}$.\\

\begin{remark}
    It is a simple computation to see that mutating at $1\rightarrow 3$ results in a seed with the edge $2\rightarrow 4$ and $A_{24} = D$.
\end{remark}

\begin{theorem}
    The algebra $\B = \CoordSPA[A^{-1},B^{-1},C^{-1},D^{-1}]$ is isomorphic to the noncommutative algebra associated to a square after specializing the variables $X_{12}=X_{21}=X_{34}=X_{43}=1$, under the $\sigma-$algebra map sending $$\{A,B,C,D\} \to \{X_{13},X_{14},X_{23},X_{24}\}.$$
\end{theorem}
\begin{proof}
    From the previous discussion we see that the seed produced from $\B$ is a valid seed. Given such a seed we must verify all the relations in $\B$ hold. Since $A$ is invertible, we only need to verify the preconditions of \Cref{thm:SymplecticEquationsAinv}. These follow exactly from the triangle conditions and the mutation rule.
\end{proof}

\begin{proposition}
    Given triangulation of a disk with 4 marked points with $A_{13} = A, A_{23} = C, A_{14} = B, A_{12} = A_{34} = 1$ the monodromy associated to the path $G_2^{13}S_{23}G_3^{21}G_3^{14}$ is $\begin{bmatrix}A & B\\ C & D\end{bmatrix}$.
\end{proposition}
\begin{proof}
    First we compute the prescribed monodromy:
    \begin{align*}
        \begin{bmatrix} 1 & AC^{-1}\\0 & 1 \end{bmatrix} &\begin{bmatrix} 0 & -\sigmaT{C^{-1}}\\C & 0 \end{bmatrix} \begin{bmatrix} 1 & C^{-1}\sigmaT{A^{-1}}\\0 & 1 \end{bmatrix}\begin{bmatrix} 1 & A^{-1}B\\0 & 1 \end{bmatrix}\\
        &\begin{bmatrix}A & -\sigmaT{C^{-1}}\\C & 0\end{bmatrix}\begin{bmatrix}1 & C^{-1}\sigmaT{A^{-1}}+A^{-1}B\\ 0 & 1\end{bmatrix}\\
        &\begin{bmatrix}
            A & B\\ C & \sigmaT{A^{-1}}+ CA^{-1}B
        \end{bmatrix}
    \end{align*}
    It remains to show that $D = \sigmaT{A^{-1}} + CA^{-1}B$. The key fact is that $CA^{-1} = \sigmaT{A^{-1}}\sigmaT{C}$ due to the symmetry of $T_3^{12}$. It is then clear the proposed formula for $D$ is equivalent to the requirement that the determinant of the matrix is 1.
    \[\sigmaT{A}D - \sigmaT{C}B = 1\]
\end{proof}

We can view the seed $S_\matSP{M}$ as a triangulation of cylinder by gluing the two identity sides together. This leaves two distinct marked points which we call $L$ and $R$. There are two natural functions on $S_M$; the sum of angles incident to each vertex. For our generic symplectic matrices these functions are:
\begin{align}\label{eqn:SymplecticSeedInvariants}
    \phi_R =& 1AC^{-1} + C^{-1}1\sigmaT{A^{-1}} + A^{-1}B1 = AC^{-1}+C^{-1}(\sigmaT{A^{-1}} + CA^{-1}B) = AC^{-1} + C^{-1}D\\
    \phi_L =& 1CA^{-1} + \sigmaT{A^{-1}}1B^{-1} + B^{-1}A1 = DB^{-1}+B^{-1}A
\end{align}

\begin{proposition}
    The functions $\phi_L$ and $\phi_R$ are invariant under mutation. Furthermore as commutative expressions both $\phi_L$ and $\phi_R$ are equal to the trace of the monodromy divided by a cluster variable.
\end{proposition}
\begin{proof}
    The total angle around a marked point is independent of triangulation and so these functions are invariant under mutation. Furthermore
    \[\phi_R = C^{-1}(CAC^{-1}+D) = (A+C^{-1}DC)C^{-1} \hspace{3pc} \phi_L = B^{-1}(BDB^{-1}+A) = (D+B^{-1}AB)B^{-1}\]
    The functions $CAC^{-1}+D$, $A+C^{-1}DC$, $BDB^{-1}+A$, and $D+B^{-1}AB$ are equal to the trace if $\A$ is commutative.
\end{proof}

\begin{proposition}
    The seed $S_{\matSP{M}\matSP{N}}$ is mutation equivalent to the seed given  by amalgamating $S_\matSP{M}$ and $S_\matSP{N}$ along an identity edge.
\end{proposition}
\begin{proof}
    This is ``clear'' via concatenating the representative path through the two seeds. However it is instructive to perform the mutation sequence  realizing the concatenated path as a path through a single square. This is done by first mutating at the newly interior identity edge then mutating at the edges labeled $A$ and $E$ (\Cref{fig:ConcatinationAsMultiplication}). \\
    The first mutation is easy to compute the new variable is $AE+BG$ as expected. However the other mutations require remembering how to express the lower right entry of each matrix in terms of the old coordinates. Recall that $D = CA^{-1}B + \sigmaT{A^{-1}}$. Then when we perform the mutation at $A$ we obtain:
    \begin{align*}
    CA^{-1}(AE+BG) + 1 \sigmaT{A^{-1}}G = CE + (CA^{-1}B + \sigmaT{A^{-1}})G = CE+DG.
    \end{align*}
\end{proof}
\begin{remark}
    The algebra generated by $S_{\matSP{M}\matSP{N}}$ is isomorphic to the algebra $T_2^*(\B)$. As such the construction of the previous proposition gives a geometric interpretation of the coproduct on $\CoordSPA$. 
\end{remark}
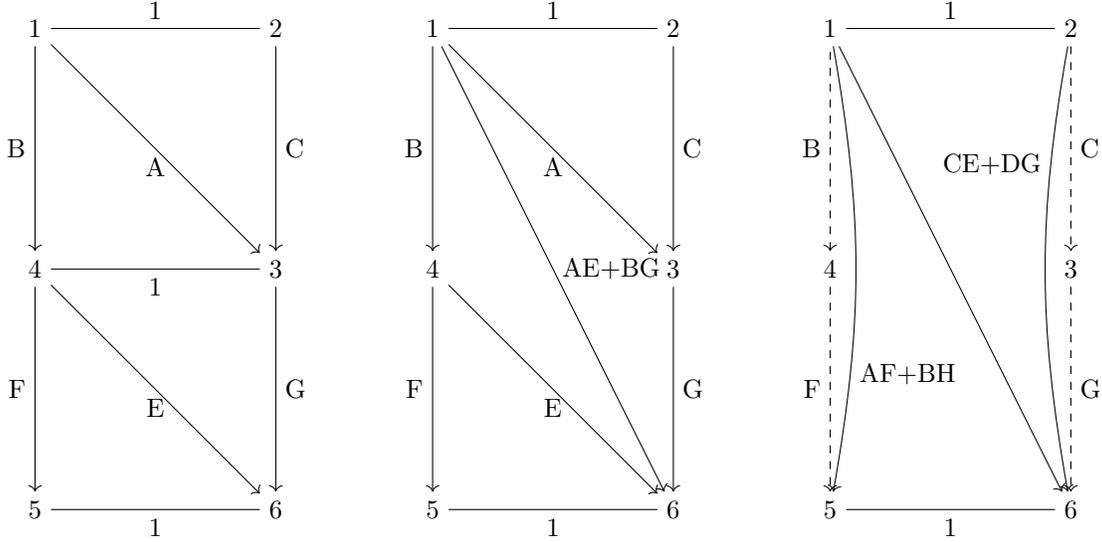
\begin{figure}[hb]
    \centering
    \begin{tikzpicture}[scale=0.8]
        \node[] at (4,4) (2) {2};
        \node[] at (0,0) (4) {4};
        \node[] at (0,4) (1) {1};
        \node[] at (4,0) (3) {3};
        \node[] at (0,-4) (5) {5};
        \node[] at (4,-4) (6) {6};
        
        %Square ABCD
        \path[] (4) edge [<-] node[left] {B}  (1);
        \path[] (1) edge [->] node[below] {A} (3);
        \path[] (1) edge [-] node[above] {1} (2);
        \path[] (2) edge [->] node[right] {C} (3);
        \path[] (3) edge [-] node[below] {1} (4); 

      %Square EFGH
        \path[] (4) edge [->] node[left] {F}  (5);
        \path[] (4) edge [->] node[below] {E} (6);
        \path[] (3) edge [->] node[right] {G} (6);
        \path[] (5) edge [-] node[below] {1} (6); 
    \end{tikzpicture}
    \hspace{2pc}
     \begin{tikzpicture}[scale=0.8]
        \node[] at (4,4) (2) {2};
        \node[] at (0,0) (4) {4};
        \node[] at (0,4) (1) {1};
        \node[] at (4,0) (3) {3};
        \node[] at (0,-4) (5) {5};
        \node[] at (4,-4) (6) {6};
        
        %Square ABCD
        \path[] (4) edge [<-] node[left] {B}  (1);
        \path[] (1) edge [->] node[below] {A} (3);
        \path[] (1) edge [-] node[above] {1} (2);
        \path[] (2) edge [->] node[right] {C} (3);
        \path[] (1) edge [->] node[right] {AE+BG} (6); 

      %Square EFGH
        \path[] (4) edge [->] node[left] {F}  (5);
        \path[] (4) edge [->] node[below] {E} (6);
        \path[] (3) edge [->] node[right] {G} (6);
        \path[] (5) edge [-] node[below] {1} (6); 
    \end{tikzpicture}
    \hspace{2pc}
     \begin{tikzpicture}[scale=0.8]
        \node[] at (4,4) (2) {2};
        \node[] at (0,0) (4) {4};
        \node[] at (0,4) (1) {1};
        \node[] at (4,0) (3) {3};
        \node[] at (0,-4) (5) {5};
        \node[] at (4,-4) (6) {6};
        
        %Square ABCD
        \path[] (4) edge [<-,dashed] node[left] {B}  (1);
        \path[bend right=10] (2) edge [->] node[left,near start] {CE+DG} (6);
        \path[] (1) edge [-] node[above] {1} (2);
        \path[] (2) edge [->,dashed] node[right] {C} (3);
        \path[] (1) edge [->] node[right] {} (6); 

      %Square EFGH
        \path[] (4) edge [->,dashed] node[left] {F}  (5);
        \path[bend left=10] (1) edge [->] node[right,near end] {AF+BH} (5);
        \path[] (3) edge [->,dashed] node[right] {G} (6);
        \path[] (5) edge [-] node[below] {1} (6); 
    \end{tikzpicture}
    \caption{Seeds associated to multiplying symplectic matrices.}
    \label{fig:ConcatinationAsMultiplication}
\end{figure}

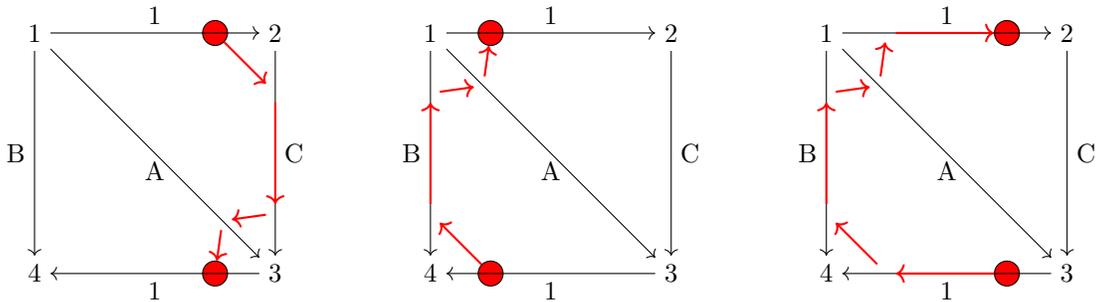
\begin{figure}[hb]
    \centering
    \begin{tikzpicture}[scale=0.8]
        \node[] at (4,4) (2) {2};
        \node[] at (0,0) (4) {4};
        \node[] at (0,4) (1) {1};
        \node[] at (4,0) (3) {3};
        \node[] at (1/1.14,4-1/1.14) (Aback) {};
        \node[] at (4-1/1.14,1/1/1.14) (Afront) {};
        \node[] at (0,4-1) (Bfront) {};
        \node[] at (0,1) (Bback) {};
        \node[] at (4,4-1) (Cback) {};
        \node[] at (4,1) (Cfront) {};
        \node[circle,draw,fill=red] at (4-1,0) (Botback) {};
        \node[] at (1,0) (Botfront) {};
        \node[circle,draw,fill=red] at (4-1,4) (Topback) {};
        \node[] at (1,4) (Topfront) {};
        %Square
        \path[] (4) edge [<-] node[left] {B}  (1);
        \path[] (1) edge [->] node[below] {A} (3);
        \path[] (1) edge [->] node[above] {1} (2);
        \path[] (2) edge [->] node[right] {C} (3);
        \path[] (3) edge [->] node[below] {1} (4);

        \path[red,thick] (Cback) edge [<-] (Topback);
        \path[red,thick] (Cback) edge [->] (Cfront);
        \path[red,thick] (Cfront) edge [->] (Afront);
        \path[red,thick] (Afront) edge [->] (Botback);
        
    \end{tikzpicture}
    \hspace{2pc}
    \begin{tikzpicture}[scale=0.8]
        \node[] at (0,0) (4) {4};
        \node[] at (0,4) (1) {1};
        \node[] at (4,0) (3) {3};
        \node[] at (4,4) (2) {2};
        \node[] at (1/1.14,4-1/1.14) (Aback) {};
        \node[] at (4-1/1.14,1/1/1.14) (Afront) {};
        \node[] at (0,4-1) (Bfront) {};
        \node[] at (0,1) (Bback) {};
        \node[] at (4,4-1) (Cback) {};
        \node[] at (4,1) (Cfront) {};
        \node[] at (4-1,0) (Botback) {};
        \node[circle,draw,fill=red] at (1,0) (Botfront) {};
        \node[] at (4-1,4) (Topback) {};
        \node[circle,draw,fill=red] at (1,4) (Topfront) {};
        %Square
        \path[] (4) edge [<-] node[left] {B}  (1);
        \path[] (1) edge [->] node[below] {A} (3);
        \path[] (1) edge [->] node[above] {1} (2);
        \path[] (2) edge [->] node[right] {C} (3);
        \path[] (3) edge [->] node[below] {1} (4); 

        \path[red,thick] (Botfront) edge [->] (Bback);
        \path[red,thick] (Bback) edge [->] (Bfront);
        \path[red,thick] (Bfront) edge [->] (Aback);
        \path[red,thick] (Aback) edge [->] (Topfront);
              
    \end{tikzpicture}
    \hspace{2pc}
        \begin{tikzpicture}[scale=0.8]
        \node[] at (0,0) (4) {4};
        \node[] at (0,4) (1) {1};
        \node[] at (4,0) (3) {3};
        \node[] at (4,4) (2) {2};
        \node[] at (1/1.14,4-1/1.14) (Aback) {};
        \node[] at (4-1/1.14,1/1/1.14) (Afront) {};
        \node[] at (0,4-1) (Bfront) {};
        \node[] at (0,1) (Bback) {};
        \node[] at (4,4-1) (Cback) {};
        \node[] at (4,1) (Cfront) {};
        \node[circle,draw,fill=red] at (4-1,0) (Botback) {};
        \node[] at (1,0) (Botfront) {};
        \node[circle,draw,fill=red] at (4-1,4) (Topback) {};
        \node[] at (1,4) (Topfront) {};
        %Square
        \path[] (4) edge [<-] node[left] {B}  (1);
        \path[] (1) edge [->] node[below] {A} (3);
        \path[] (1) edge [->] node[above] {1} (2);
        \path[] (2) edge [->] node[right] {C} (3);
        \path[] (3) edge [->] node[below] {1} (4); 

        \path[red,thick] (Botfront) edge [->] (Bback);
        \path[red,thick] (Bback) edge [->] (Bfront);
        \path[red,thick] (Bfront) edge [->] (Aback);
        \path[red,thick] (Aback) edge [->] (Topfront);

        \path[red,thick] (Botback) edge [->] (Botfront);
        \path[red,thick] (Topback) edge [<-] (Topfront);       
    \end{tikzpicture}
    \caption{Paths corresponding to the original matrix, its transpose, and its inverse.}
    \label{fig:SymplecticInversion}
\end{figure}

\begin{proposition}
    The seed  $S_{\sigmaT{\matSP{M}}^T}$ is given by rotating the seed $S_\matSP{M}$ by 180 degrees.
\end{proposition}
\begin{proof}
    This rotation reverses the orientation of the diagonal edge and swaps the two vertical edges while also reversing their orientation (\Cref{fig:SymplecticSeedTranspose}). To fix the orientation each variable is acted on by $\sigma$. The resulting matrix is 
    \[\begin{bmatrix}
        \sigmaT{A} & \sigmaT{C}\\
        \sigmaT{B} & \sigmaT{D}
    \end{bmatrix}\]
\end{proof}
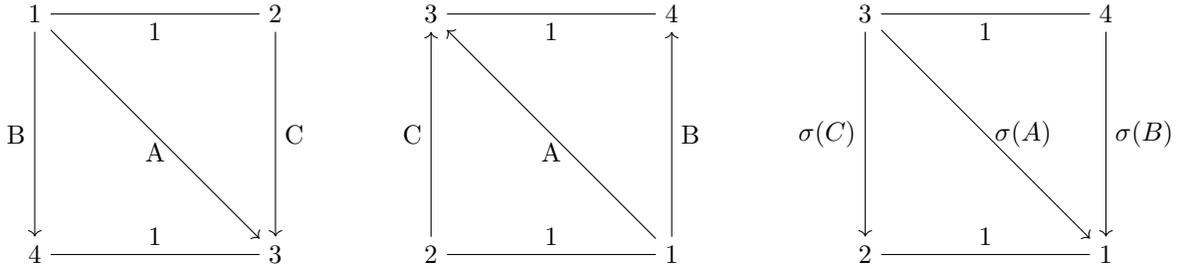
\begin{figure}[hb]
    \centering
    \begin{tikzpicture}[scale=0.8]
        \node[] at (4,4) (2) {2};
        \node[] at (0,0) (4) {4};
        \node[] at (0,4) (1) {1};
        \node[] at (4,0) (3) {3};
        
        %Square ABCD
        \path[] (4) edge [<-] node[left] {B}  (1);
        \path[] (1) edge [->] node[below] {A} (3);
        \path[] (1) edge [-] node[below] {1} (2);
        \path[] (2) edge [->] node[right] {C} (3);
        \path[] (3) edge [-] node[above] {1} (4); 
    \end{tikzpicture}
    \hspace{2pc}
    \begin{tikzpicture}[scale=0.8]
        \node[] at (4,4) (4) {4};
        \node[] at (0,0) (2) {2};
        \node[] at (0,4) (3) {3};
        \node[] at (4,0) (1) {1};
        
        %Square ABCD
        \path[] (4) edge [<-] node[right] {B}  (1);
        \path[] (1) edge [->] node[below] {A} (3);
        \path[] (1) edge [-] node[above] {1} (2);
        \path[] (2) edge [->] node[left] {C} (3);
        \path[] (3) edge [-] node[below] {1} (4); 
    \end{tikzpicture}
    \hspace{2pc}
    \begin{tikzpicture}[scale=0.8]
        \node[] at (4,4) (2) {4};
        \node[] at (0,0) (4) {2};
        \node[] at (0,4) (1) {3};
        \node[] at (4,0) (3) {1};
        
        %Square ABCD
        \path[] (4) edge [<-] node[left] {$\sigmaT{C}$}  (1);
        \path[] (1) edge [->] node[right] {$\sigmaT{A}$} (3);
        \path[] (1) edge [-] node[below] {1} (2);
        \path[] (2) edge [->] node[right] {$\sigmaT{B}$} (3);
        \path[] (3) edge [-] node[above] {1} (4); 
    \end{tikzpicture}
    \caption{Seed associated to transpose of a symplectic matrix.}
    \label{fig:SymplecticSeedTranspose}
\end{figure}
\begin{corollary}
    The transpose of a symplectic matrix is symplectic. 
\end{corollary}
\begin{corollary}
    The inverse of of a symplectic matrix $\matSP{M}$ is given by $\matSPgreek{\Omega}_1 \matSP{M}^{T}\matSPgreek{\Omega}_1$
\end{corollary}
\begin{proof}
    When rotating the seed 180 degrees we see that the path corresponding to the original matrix and the path corresponding to the transpose matrix differ by the segments along the identity edge $\matSP{S}_{12}$ and $\matSP{S}_{34}$ (See \Cref{fig:SymplecticInversion}). In this case both matrices are $-\matSPgreek{\Omega}_1$. Thus
    \[\matSP{M}^{-1} = (-\matSPgreek{\Omega}_1)\sigmaT{\matSP{M}}^T(-\matSPgreek{\Omega}_1) = \matSPgreek{\Omega}_1 \sigmaT{\matSP{M}}^T \matSPgreek{\Omega}.\]
    We note that the homotopy identifying the two paths passes through two triangles and so picks up a factor of $(-\id)^2$ resulting in equality at the level of $\SPA$ not just $\PSPA$.
\end{proof}

\subsection{Quantization}

We can easily define a quantization of the noncommutative cluster algebra associated to a square which agrees with the quantized $(k_q,\sigma)$-algebra defined earlier. We will now assume that $\A$ is a $(k,\sigma)$-algebra.

\begin{definition}
    Let $i,j,k$ be in cyclic order agreeing with the orientation of the surface. We define the \emph{quantized noncommutative angle} $T_{q,i}^{jk} = q^{1/2}T_i^{jk}$. If $i,j,k$  disagrees with orientation of the surface, then we define  $T_{q,i}^{jk} = q^{-1/2}T_i^{jk}$.
\end{definition}

Our triangle condition becomes $$ \sigma(T_{q,i}^{jk}) = T_{q,i}^{kj} $$ recalling that $\sigma(q)=q^{-1}$.

\begin{lemma}
    The symmetry of one angle in a triangle implies the symmetry of the other angles. 
\end{lemma}
\begin{proof}
    Consider a triangle with vertices cyclically labeled $1,2,3$ and edge functions $X_{ij}$. We assume that $\sigma(T_{q,1}^{23}) = T_{q,1}^{32}$, that is that $$q^{1/2}X_{12}^{-1}X_{32}X_{31}^{-1}= q^{-1/2}X_{13}^{-1}X_{23}X_{21}^{-1}.$$ Then we have that $$q^{1/2}X_{32}X_{31}^{-1}X_{21}$$ is also symmetric. This implies that $$q^{-1/2}X_{21}^{-1}X_{31}X_{32}^{-1} = T_{q,2}^{13}$$ is also symmetric.
\end{proof}

Additivity of the angles implies the mutation rule, which is the same as in the unquantized case. Thus, if we start with a configuration of variables satisfying the the triangle conditions, then the variables produced after any number of mutations will also satisfy the triangle conditions.

\begin{theorem}
    The quantized noncommutative cluster algebra associated with a square with variable $1\in k$ assigned to the top and bottom edges is isomorphic as $(k_q,\sigma)$-algebras to $\CoordSPA_q[A^{-1},B^{-1},C^{-1},D^{-1}]$
\end{theorem}
\begin{proof}
    Picking the triangulation of \Cref{fig:SymplecticSeed}, we can see that the symmetry of the angles implies the relations
    \begin{align*}
     \sigma(A)C-q^{-1}\sigma(C)A =& 0\\
     A\sigma(B)-q^{-1} B\sigma(A) =& 0
 \end{align*}
 Moreover, the exchange relation implies that $D = \sigma(A^{-1}) + CA^{-1}B$ so that $D\sigma(A) = 1+ qC\sigma(B)$ which implies the last of the relations.
\end{proof}

\begin{lemma}
        The quantized noncommutative cluster algebra associated with a triangle with $X_{21} = X, X_{23} =1, X_{31}= Y$ is isomorphic as a $(k_q,\sigma)$-algebra to the quantum plane $(\A^2_\omega)_q$.
\end{lemma}

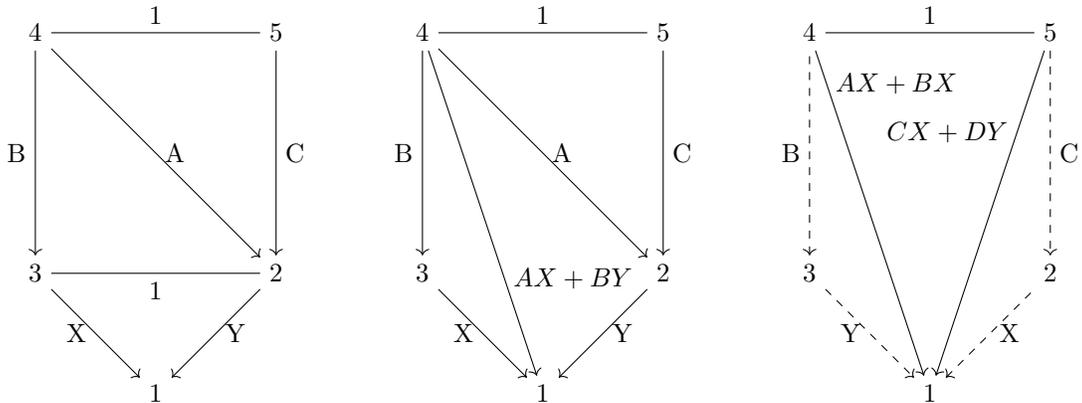
\begin{figure}[ht]
    \centering
    \begin{tikzpicture}[scale=0.8]
        \node[] at (2,-2) (1) {1};
        \node[] at (0,0) (3) {3};
        \node[] at (4,0) (2) {2};
        \node[] at (0,4) (4) {4};
        \node[] at (4,4) (5) {5};
        %Triangle
        \path[] (1) edge [<-] node[right] {Y}  (2);
        \path[] (2) edge [-] node[below] {1}  (3);
        \path[] (3) edge [->] node[left] {X}  (1);
        %Square
        \path[] (3) edge [<-] node[left] {B} (4);
        \path[] (4) edge [] node[above] {1} (5);
        \path[] (5) edge [->] node[right] {C} (2);
        \path[] (4) edge [->] node[right] {A} (2);
    \end{tikzpicture}\hspace{2pc}
    \begin{tikzpicture}[scale=0.8]
        \node[] at (2,-2) (1) {1};
        \node[] at (0,0) (3) {3};
        \node[] at (4,0) (2) {2};
        \node[] at (0,4) (4) {4};
        \node[] at (4,4) (5) {5};
        %Triangle
        \path[] (1) edge [<-] node[right] {Y}  (2);
        \path[] (4) edge [->] node[right,pos=0.7] {$AX+BY$}  (1);
        \path[] (3) edge [->] node[left] {X}  (1);
        %Square
        \path[] (3) edge [<-] node[left] {B} (4);
        \path[] (4) edge [] node[above] {1} (5);
        \path[] (5) edge [->] node[right] {C} (2);
        \path[] (4) edge [->] node[right] {A} (2);
    \end{tikzpicture}\hspace{2pc}
    \begin{tikzpicture}[scale=0.8]
        \node[] at (2,-2) (1) {1};
        \node[] at (0,0) (3) {3};
        \node[] at (4,0) (2) {2};
        \node[] at (0,4) (4) {4};
        \node[] at (4,4) (5) {5};
        %Triangle
        \path[] (1) edge [<-,dashed] node[right] {X}  (2);
        \path[] (4) edge [->] node[pos=.1,right] {$AX+BX$}  (1);
        \path[] (3) edge [->,dashed] node[left] {Y}  (1);
        %Square
        \path[] (3) edge [<-,dashed] node[left] {B} (4);
        \path[] (4) edge [] node[above] {1} (5);
        \path[] (5) edge [->,dashed] node[right] {C} (2);
        \path[] (5) edge [->] node[left,near start] {$CX+DY$} (1);
    \end{tikzpicture}
    \caption{Action of $\CoordSPA_q$ on $(\A^2_\omega)_q$}
    \label{fig:SeedQuandtumAction}
\end{figure}

\begin{theorem}
    The coaction of $\CoordSPA_q$ on $(\A^2_\omega)_q$ is given by gluing a square for $\CoordSPA_q$ above a triangle for $(\A^2_\omega)_q$ and performing two mutations to create a new triangle with top variable equal to 1 (See \Cref{fig:SeedQuandtumAction}).
\end{theorem}

\section{Noncommutative Markov Triples}\label{sec:NoncomMarkov}

There is a well known relationship between the cluster algebra associated to the once punctured torus and the Markov numbers see e.g. \cite{FG-DualTeichmuller},\cite{FWZ-IntroToClusters}. The Markov equation itself is an example of a mutation invariant function in the sense of \cite{K-MutationInvariants}, and its solutions are given by evaluating the initial cluster variables at $1$. 

Bringing this relationship to the noncommutative setting, we obtain generalizations of Markov numbers to any involutive algebra $(\A,\sigma)$. We accomplish this by computing a noncommutative analogue of a mutation invariant function and then by choosing some initial set of elements of $\A$ to evaluate the cluster variables at. 

First we recall the usual story of the Markov cluster algebra.
\subsection{Commutative Markov Triples}
The classical Markov equation is given by
\begin{equation}\label{eqn:MarkovClassic}
    a^2 + b^2 + c^2 = 3abc.
\end{equation}
The integer solutions to the Markov equation are called \keyword{Markov triples}. An initially surprising fact is that given a Markov triple $(a,b,c)$ the triple $(\frac{b^2+c^2}{a},b,c)$ is also a Markov triple. As proved by Markov, every Markov triple can be obtained from the ``trivial'' solution $(1,1,1)$ by repeatedly applying this rule at each entry.

This structure of the Markov numbers can be obtained from the cluster structure on a punctured torus. Each Markov triple corresponds to a seed of this cluster algebra (\Cref{fig:torusQuiver}). The cluster mutation rule for this seed states the mutated coordinate $a'$  satisfies the relation $a a' = b^2 + c^2$. This is exactly the mutation of Markov triples described above. There is a mutation invariant function on this cluster algebra $f(a,b,c) = 2\frac{a^2+b^2+c^2}{abc}$. When starting with the triple $(1,1,1)$ we see that $f(1,1,1) = 6$. Since this function is fixed by mutation, every other triple obtained by mutation also satisfies $f(a,b,c) = 6$, which is equivalent to the Markov equation. Furthermore, the fact that all triples obtained by mutation are positive integers is a consequence of the positivity of the cluster exchange relation.

We obtain a geometric interpretation of the Markov triples using the relationship between this cluster structure and  hyperbolic structures on the punctured torus with a choice of horocycle at the puncture. Under this interpretation each Markov number corresponds to an arc on the torus and each triple corresponds to a triangulation. Each Markov number is the ``lambda length'' of the corresponding arc \cite{FT-LambdaLengths}. There is also a geometric interpretation of the mutation invariant function as the sum of ``angles'' $t_i^{jk} = \frac{a_{jk}}{a_{ij}a_{ki}}$ in the triangulation. Each of these angles is the length of the horocycle  inside the triangle whose edges are labeled with the cluster variables $a_{ij},a_{jk},a_{ki}$ (\Cref{thm:MagicHorocycleLength}). On the punctured torus, the total length of the horocycle is then given by the \keyword{Markov function}:
\begin{equation}
    f(a,b,c) = 2\left(\frac{a}{bc}+\frac{b}{ac}+\frac{c}{ab}\right) = 2\frac{a^2+b^2+c^2}{abc}
\end{equation}
This length is clearly independent of the chosen triangulation. Furthermore the mutation of the Markov quiver leaves the adjacency of the triangles fixed. Thus exact form of the function $f$ is invariant under mutation.

\subsection{Noncommutative Generalization}
To define the noncommutative generalization of Markov triples to an involutive algebra $(\A,\sigma)$ we use the noncommutative cluster structure of \cite{BR-noncommutativeClusters} recalled in \Cref{sec:surfaceClusterStructure}. We consider the (noncommutative seed) associated to the punctured torus (\Cref{fig:torusQuiver}). The key difference from the classic case is that the lengths of horocycle arcs $t_i^{jk}$ are ``upgraded'' to the noncommutative angle $T_i^{jk}$. We require these angles to be fixed by the anti-involution $\sigma$.\\
Explicitly these angles are all cyclic permutations of:
\[\sigmaT{A^{-1}}B\sigmaT{C^{-1}} \in \A^\sigma \hspace{2pc} \sigmaT{A^{-1}}C\sigmaT{B^{-1}} \in \A^\sigma \]

\begin{figure}[hb]
    \centering
    \begin{tikzpicture}
    \node[] at (4,4) (2) {3};
\node[] at (0,0) (1) {1};
\node[] at (0,4) (3) {2};
\node[] at (4,0) (4) {4};
\path[] (2) edge [<-] node[below] {$A$}  (1);
\path[] (3) edge [<-] node[below] {$C$} (2);
\path[] (1) edge [<-] node[right] {$B$} (3);
\path[] (2) edge [->] node[left] {$B$} (4);
\path[] (4) edge [->] node[above] {$C$} (1);
\end{tikzpicture}
   \hspace{2pc}
%    \begin{tikzpicture}[]
% \node[circle,draw] at (0, 0) (2) [left
%  ] {$A$} ;
% \node[circle,draw] at (2, 3.464) (1) [] {$B$} ;
% \node[circle,draw] at (4,0) (3) [] {$C$} ;
% \path[{Square[normal,open,length=3mm]}-{Latex[transpose,open,length=3mm]}, shorten >=1.5mm,shorten <=1.5mm] (2.90) edge [] node {} (1.-150);
% \path[{Square[normal,length=3mm]}-{Latex[transpose,length=3mm]}, shorten >=1.5mm,shorten <=1.5mm] (2.30) edge [] node {} (1.-90);
% \path[{Square[normal,open,length=3mm]}-{Latex[transpose,open,length=3mm]}, shorten >=1.5mm,shorten <=1.5mm] (1.-30) edge [] node {} (3.90);
% \path[{Square[normal,length=3mm]}-{Latex[transpose,length=3mm]}, shorten >=1.5mm,shorten <=1.5mm] (1.-90) edge [] node {} (3.150);
% \path[{Square[normal,open,length=3mm]}-{Latex[transpose,open,length=3mm]}, shorten >=1.5mm,shorten <=1.5mm] (3.210) edge [] node {} (2.-30);
% \path[{Square[normal,length=3mm]}-{Latex[transpose,length=3mm]}, shorten >=1.5mm,shorten <=1.5mm] (3.150) edge [] node {} (2.30);
% \end{tikzpicture}
     \caption{Seed associated to a triangulation of a punctured torus}
     \label{fig:torusQuiver}
\end{figure}
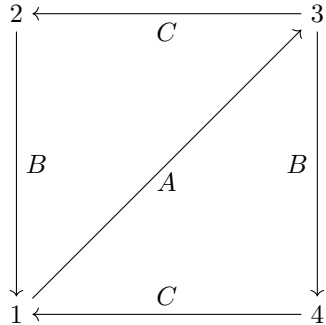
\begin{remark}\label{rem:AngleCheck}
    It suffices to verify the symmetry of one equivalence class of each kind. Each cyclic rotate can be written from the starting angle by $\sigma$-congruence as follows:
    \[\sigmaT{B^{-1}}C \sigmaT{A^{-1}} = \sigma\left[\sigmaT{A}\cdot \sigmaT{A^{-1}}B\sigmaT{C^{-1}} \cdot A\right]^{-1}\]
\end{remark}

\begin{definition}
    The \keyword{noncommutative mutation} at $A$ in \Cref{fig:torusQuiver} is given by the equation:
    \begin{equation}\label{eqn:MarkovMutation}
        A' = \sigma(C)A^{-1}\sigma(C) + B\sigma(A^{-1})B
    \end{equation}
\end{definition}

\begin{remark}
    The symmetry of the placements of $\sigma$ in the equation above is due to the fact that the arcs of the triangulation are oriented in cycles. To preserve this property after mutation we reverse the orientation of the arc labeled $B$. Thus the explicit mutation of triples becomes:
    \[(A,B,C) \mapsto (A',C,\sigma(B)).\]
    This move combined with cyclically rotating triples can be combined to obtain every seed. For example,  mutation along the path given by freezing $C$ corresponds to rotating right once after mutating. In other words iterate the following transformation:
    \[(A,B,C) \mapsto (\sigmaT{B},A',C).\]
\end{remark} 

\begin{lemma}\label{thm:GeneralizedMarkovFunction}
    There is a mutation invariant function $F(A,B,C)$ for the cluster structure on the punctured torus. 
\end{lemma}
\begin{proof}
    The sum of all six noncommutative angles is the analogue of the total length of the horocycle in the commutative case. In Theorem 3.40 of \cite{BR-noncommutativeClusters} they show the total angle at a marked point is independent of triangulation, and thus mutation invariant. Furthermore by fixing the cyclic orientation, $(A,B,C) \mapsto (A', C, \sigma(B))$ the expression for this total angle is identical before and after mutation. As such $F$ itself is mutation invariant.
\end{proof}
 
 \begin{definition}
     We call the function from \Cref{thm:GeneralizedMarkovFunction} the \keyword{generalized Markov function}. Explicitly on a seed $(A,B,C)$ whose triangles are cyclically oriented: 
     \begin{equation}\label{eqn:NonComMarkov}
     \begin{aligned}
         F(A,B,C) =& ~\hphantom{+}\sigma(A^{-1})B\sigma(C^{-1}) + \sigma(B^{-1})C\sigma(A^{-1})+\sigma(C^{-1})A\sigma(B^{-1})\\
         &+\sigma(C^{-1})B\sigma(A^{-1})+\sigma(A^{-1})C\sigma(B^{-1})+\sigma(B^{-1})A\sigma(C^{-1})
         \end{aligned}
     \end{equation}
 \end{definition}
\begin{corollary}
    Let $F_0 = F(A_0,B_0,C_0)$ be the value of the generalized Markov function on an initial seed $(A_0,B_0,C_0)$. Then every other triple $(A,B,C)$ obtained by mutation satisfies the generalized Markov equation $F(A,B,C) = F_0$.
\end{corollary}

 \begin{remark}
     If $\sigma$ is trivial then $a,b,c$ commute. This simplifies $F$ to $F(a,b,c) = 2\frac{a^2+b^2+c^2}{abc}$ which is the standard Markov function (\Cref{eqn:MarkovClassic}).
 \end{remark}
 \begin{remark}
     If the initial seed $(A,B,C)$ is chosen to be $(1,1,1)$ we recover the standard Markov numbers as the solutions to the generalized Markov equation. 
 \end{remark}

\subsection{The Fibonacci Branch}
 The classic Markov numbers that belong to a triple $(a,b,1)$ are odd indexed Fibonacci numbers $f_{2n+1}$. These triples are obtained by ``freezing'' one position labeled 1 in the initial seed and alternating mutation at the other two positions. We obtain a recurrence for these variables by combining the mutation rule and Markov function. The classic mutation at $a_{i-1}$ in the triple $(a_{i-1},a_{i},1)$ produces $a_{i+1} = \frac{a_{i}^2+1}{a_{i-1}}$. Then the Markov function can be simplified as follows 
 \[ F(a_{i-1},a_{i},1) = 2\frac{a_{i-1}^2+a_{i}^2+1}{a_{i-1}a_{i}} = \frac{2}{a_{i}}\left(a_{i-1}+ \frac{a_{i}^2+1}{a_{i-1}}\right) = \frac{2}{a_{i}}\left(a_{i-1}+ a_{i+1}\right) \]
 As the Markov function is always some constant $F_0$, we see that $a_i$ satisfies the recurrence $a_{i+1} = \frac{1}{2}F_0 a_{i} - a_{i-1}$. Starting with the triple $(2,1,1)$ means that $a_0 = 1 = f_1$,  $a_1 = 2 = f_3$ and $F_0 = 6$. The recurrence then simplifies to 
 \[a_{n+1} = 3a_{n}-a_{n-1}.\]
 This is exactly the recurrence satisfied by the odd indexed Fibonacci numbers.

 We can generalize this recurrence in the noncommutative case. The key step is the ability to rewrite the Markov function as a sum involving a variable and its mutation. So let
\[\halfInvariant(A,B,C) = \sigmaT{C^{-1}}A\sigmaT{B^{-1}} + \sigma(A^{-1})B\sigma(C^{-1}) + B^{-1}\sigma(C)A^{-1} = \sigma(C^{-1})A\sigma(B^{-1}) + B^{-1}A'\sigma(C^{-1})\] where $A'$ is the new coordinate obtained by mutating at $A$ in the triple $(A,B,C)$. Over a commutative ring $\halfInvariant(A,B,C)$ is equal to half of the Markov function and it will play the same role in obtaining the recurrence in the noncommutative case.

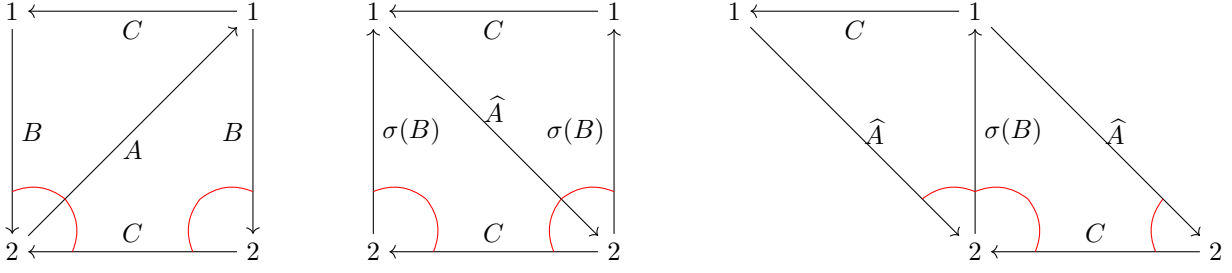
\begin{figure}
    \centering
    \begin{tikzpicture}[scale=0.8]
    %Square One
        \node[] at (4,4) (3) {1};
        \node[] at (0,0) (1) {2};
        \node[] at (0,4) (2) {1};
        \node[] at (4,0) (4) {2};
        \path[] (3) edge [<-] node[below] {$A$}  (1);
        \path[] (2) edge [<-] node[below] {$C$} (3);
        \path[] (1) edge [<-] node[right] {$B$} (2);
        \path[] (3) edge [->] node[left] {$B$} (4);
        \path[] (4) edge [->] node[above] {$C$} (1);
        \path[bend left,red] (0,1) edge [-] node[above] {} (1/1.14,1/1.14);
        \path[bend right,red] (1,0) edge [-] node[above] {} (1/1.14,1/1.14);
        \path[bend left,red] (3,0) edge [-] node[above] {} (4-1/1.14,1/1.14);
        \path[bend left,red] (4-1/1.14,1/1.14) edge [-] node[above] {} (4,1);
    %Mutated
        \node[] at (6+4,4) (3) {1};
        \node[] at (6+0,0) (1) {2};
        \node[] at (6+0,4) (2) {1};
        \node[] at (6+4,0) (4) {2};
        \path[] (4) edge [<-] node[above] {$\hat{A}$}  (2);
        \path[] (2) edge [<-] node[below] {$C$} (3);
        \path[] (1) edge [->] node[right] {$\sigmaT{B}$} (2);
        \path[] (3) edge [<-] node[left] {$\sigmaT{B}$} (4);
        \path[] (4) edge [->] node[above] {$C$} (1);
        \path[bend left,red] (6+0,1) edge [-] node[above] {} (6+1/1.14,1/1.14);
        \path[bend right,red] (6+1,0) edge [-] node[above] {} (6+1/1.14,1/1.14);
        \path[bend left,red] (6+3,0) edge [-] node[above] {} (6+4-1/1.14,1/1.14);
        \path[bend left,red] (6+4-1/1.14,1/1.14) edge [-] node[above] {} (6+4,1);
    %Reconfigured to set up next mutation
        \node[] at (16-4,4) (2) {1};
        \node[] at (16+0,0) (1) {2};
        \node[] at (16+0,4) (3) {1};
        \node[] at (16+4,0) (4) {2};
        \path[] (3) edge [<-] node[right] {$\sigmaT{B}$}  (1);
        \path[] (2) edge [<-] node[below] {$C$} (3);
        \path[] (1) edge [<-] node[right] {$\hat{A}$} (2);
        \path[] (3) edge [->] node[right] {$\hat{A}$} (4);
        \path[] (4) edge [->] node[above] {$C$} (1);
        \path[bend left,red] (16-1/1.14,1/1.14) edge [-] node[above] {} (16+0,1);
        \path[bend right,red] (16+1,0) edge [-] node[above] {} (16+1/1.14,1/1.14);
        \path[bend right,red] (16+1/1.14,1/1.14) edge [-] node[] {} (16+0,1);
        \path[bend left,red] (16+3,0) edge [-] node[above] {} (16+4-1/1.14,1/1.14);
    \end{tikzpicture}
    \caption{Mutation preserves the angles needed to compute $\halfInvariant$.}
    \label{fig:AnglePreservation}
\end{figure}

\begin{lemma}
    The function $\halfInvariant(A,B,C)$ is invariant on all triples obtained by mutation with $C$ frozen.
\end{lemma}
\begin{proof}
    Since we never mutate at $C$, at every triangulation obtained by mutation we can cut along $C$ to obtain the same cylinder with two marked points. The function $\halfInvariant(A,B,C)$ is the sum of noncommutative angles around the marked point labeled 2 in \Cref{fig:AnglePreservation}. Thus by \cite{BR-noncommutativeClusters}, this function is mutation invariant.
\end{proof}
From \Cref{fig:AnglePreservation} we see that the first form of $\halfInvariant(A,B,C)$ writes the sum of the angles in one triangulation, while the second form uses one angle from before and after mutation.
\begin{remark}
    When $C=1$, $\halfInvariant(A,B,1) = A\sigma(B^{-1}) + B^{-1}A'$ is related to the trace of $\begin{bmatrix}
        A & B\\B & A'
    \end{bmatrix}$. Over a commutative ring $\halfInvariant(A,B,1)$ is the trace divided by $B$.
\end{remark}

\begin{theorem}
If $C = \sigmaT{C}$ the functions occurring in the path $(F_{i-1},F_{i},C)$ satisfy the recurrence relation \[F_{i+1} = F_{i}\halfInvariant C - \sigmaT{F_{i-1}}.\]
\end{theorem}
\begin{proof}
    This follows by rearranging the equation for $\halfInvariant$.
    \begin{align*}
        \halfInvariant =~& \sigmaT{C^{-1}}F_{i-1}\sigmaT{F_{i}^{-1}}+F_{i}^{-1}F'_{i-1}\sigmaT{C^{-1}}\\
        \halfInvariant=~&\sigmaT{C^{-1}}F_{i-1}\sigmaT{F_{i}^{-1}}+F_{i}^{-1}F_{i+1}\sigmaT{C^{-1}}\\
        \halfInvariant-\sigmaT{C^{-1}}F_{i-1}\sigmaT{F_{i}^{-1}} =~&F_{i}^{-1}F_{i+1}\sigmaT{C^{-1}}\\
        \halfInvariant-\sigmaT{C^{-1}}F_{i-1}\sigmaT{F_{i}^{-1}} =~&C^{-1}\sigmaT{F_{i+1}}\sigmaT{F_{i}^{-1}}\\
        C\halfInvariant \sigmaT{F_{i}} - F_{i-1} =~& \sigmaT{F_{i+1}}\\
        F_{i}\halfInvariant C - \sigmaT{F_{i-1}} =~& F_{i+1}
    \end{align*}
\end{proof}
 \begin{remark}
     When $C= 1$ this recurrence simplifies to 
 \[ F_{i+1} = F_{i}\halfInvariant -\sigmaT{F_{i-1}}\]
     As such we see a noncommutative generalization of odd indexed Fibonacci numbers sitting inside the noncommutative Markov numbers.
 \end{remark}
%\begin{remark}
%    The seed for the markov quiver when $C= 1$ is exactly the seed associated to a symplectic matrix $M$ whose top left and bottom right entries are equal. This causes $M$ to be $\sigma$ stable.
%    \[M= \begin{bmatrix}
%        A & B\\
%        B & D
%    \end{bmatrix}\]
    
%    Furthermore after each mutation the seed can be rearranged to reform the standard symplectic seed (See \Cref{fig:AnglePreservation}). As such the mutation path along $C=1$ corresponds to sequence of  $\sigma$ stable symplectic matrices each sharing one row.
%    \[\dots \leftrightarrow \begin{bmatrix}\sigmaT{B^{-1}}+AB^{-1}A & A\\ A & B\end{bmatrix}
%            \leftrightarrow\begin{bmatrix} A & B\\ B & D \end{bmatrix} 
%            \leftrightarrow \begin{bmatrix}B & D\\ D & \sigmaT{B^{-1}}+DB^{-1}D \end{bmatrix}
%    \leftrightarrow \dots
%    \]
%    \end{remark}

\subsection{Generating Noncommutative Markov Triples}\label{sec:GeneratingNonComMarkovTriples}

We would now like to pick some involutive ring $\A$ and compute the noncommutative Markov numbers within this ring by choosing some initial seed. Even in the classic case of $\R$ it is nontrivial to specify the conditions that uniquely recover Markov numbers from the cluster theory. Any initial triple $(a_0,b_0,c_0)$ satisfies a Markov like equation $a^2+b^2+c^2= f(a_0,b_0,c_0) abc$. Furthermore from the cluster theory there is nothing special about the integers as a subring of $\R$. The Laurent phenomena for cluster algebras states that every number obtained by mutation can be expressed as a Laurent polynomial in the original triple $(a_0,b_0,c_0)$. Thus if $(a_0,b_0,c_0)$ are units in any subring all values obtained by mutation belong to that subring.\\
Nevertheless the triple $(1,1,1)$ is special as the numbers obtained from this triple have many applications in number theory. As such we look for initial triples in an analogue of the integers inside our larger rings and conjecture these triples are more interesting from a number theory perspective then an arbitrary choice. Recall the integral elements of a ring $\A$ are the roots of monic polynomials with integer coefficients in $\A$. For $\A$ noncommutative this set may not be a subring. Thus we will require that our starting triples are units in some subring of the integral elements. 
\begin{example}
    Over $\R$ any algebraic number is an integral element. However we focus on $\Z$ (the integral elements of $\Q$). In principle choosing invertible algebraic numbers as initial elements also yields interesting Markov-like numbers. 
\end{example}
\begin{example}
    Over $\Mat_n(\R)$ we focus on $\Mat_n(\Z)$. This subring is large enough to contain infinite families of units like $\begin{bmatrix}
        1 & t\\0 & 1
    \end{bmatrix}$ for any $t \in \Z$. Using these families we can consider nontrivial deformations of Markov numbers inside $\Mat_n(\R)$.
\end{example}

Another subtly when choosing an initial triple is that every element of the cluster algebra is invertible. Simply choosing a triple of invertible elements isn't enough to guarantee that all elements obtained by mutation are invertible. For example in $\C$ with the trivial involution, mutation at 1 in the triple $(1,1,i)$ results in the triple $(0,1,i)$. A sufficient solution to avoid these triples is to require all noncommutative angles in the initial seed are ``positive'' (\Cref{thm:PostiveStructure}). Over $\R$ it is clear what positive means, but over arbitrary involutive rings it is less clear. For \Cref{thm:PostiveStructure} to work, we define a positive structure on $(\A,\sigma)$ to be subset $P$ of $\A^\sigma$ containing 1 and be closed under  addition, inversion, and $\sigma$-congruence by invertible elements ($A \mapsto \sigma(B)AB$). 
\begin{remark}
    The requirement that the noncommutative angles are positive is different than the requirement the variables themselves are positive. Over $\R$ the triple $(-1,-1,1)$ has all positive angles even though some of the entries are negative. However when $\sigma$ is nontrivial requiring that the variables themselves are positive, would require the variables to be fixed by $\sigma$, dramatically reducing the possible starting triples.
\end{remark}

A final consideration is the value of the Markov function. As this is mutation invariant, if two seeds $(A,B,C)$ and $(A',B',C')$ have distinct values of the Markov function they produce distinct sets of Markov numbers. Over $\R$, fixing the value of the Markov function is enough to obtain a unique positive integer tree. We will see examples where distinct integer trees have the same value of the Markov function in larger rings.
\begin{definition}
    We fix a subring $R$ of integral elements of $\A$ and a positive structure $P$ on $(\A,\sigma)$. An \keyword{admissible starting triple} for $(R,P)$ is a triple $(a,b,c) \in \A^3$ such that $a,b,c$ are units in $R$ and all the noncommutative angles are in $P$.  
\end{definition}

\begin{lemma}\label{thm:MarkovInverseTriple}
     For any subring $R$ and $B$ a unit in $R$, the triple $(1,B,B^{-1})$ is admissible.
\end{lemma}
\begin{proof}
    It is clear that $1,B,B^{-1}$ are all units in $R$ so it suffices to check that the noncommutative angles are positive. We compute that the angles are 
    \[1,1,B\sigma(B),\sigma(B)B, B^{-1}\sigma(B^{-1}), \sigma(B^{-1})B^{-1}\]
    These are all clearly fixed by $\sigma$. Furthermore they are all obtained from $1$ via $\sigma$-congruence by an invertible element and thus must be in the positive cone.
\end{proof}

We now look at examples of admissible triples in key examples of involutive rings. In the complex numbers (\Cref{sec:ComplexMarkov}) and dual numbers (\Cref{sec:DualMarkov}) we observe the phenomena of ``shadow'' Markov numbers. Here the new Markov numbers split into a piece corresponding to the classic Markov numbers and a shadow with its own interesting mutation and dynamics.
\begin{remark}\label{rem:MarkovSymmetry}
    The usual Markov tree starting from $(1,1,1)$ has additional symmetry as mutation at any variable results in the triple $(1,1,2)$. Similarly mutation at either $1$ in this triple results in the triple $(1,2,5)$. When the initial triple has distinct entries the new mutation tree won't have this symmetry. As such for Markov numbers greater than 2, there is potential for at least six distinct shadows. If the conjecture that each Markov number appears once as the largest term in a triple then there will be at most six possible shadows.
\end{remark}

In the ring of matrices (\Cref{sec:MatrixMarkov}) we have our first examples of genuinely noncommutative Markov numbers. Many of these examples can be thought of as deformations as the standard Markov numbers, they are parameterized families of Markov triples extending from $(1,1,1)$. Finally we look at the group rings of finite cyclic groups in \Cref{sec:GroupRingMarkov}. Here we see a commutative example that doesn't split cleanly into a real and shadow part.

\subsection{Complex Markov Numbers}\label{sec:ComplexMarkov}
We begin by establishing the subring of integral elements we chose to obtain interesting new Markov-like trees. Here the full set of integral elements are in fact a subring of $\C$ as $\C$ is commutative. A key example of integral elements are the roots  of unity, which satisfy $\zeta^n -1 = 0$. From this we see this subring is much larger than the ring of Gaussian integers $\{a + b i \mid a ,b\in \Z\}$.\\
Over $\C$ we also have two reasonable choices of involution. The first is that $\C$ commutes so we can take $\sigma$ trivial. The other natural choice for $\sigma$ is complex conjugation. 
\begin{lemma}
    There is no possible choice of positive structure on $\C$ with the trivial involution.
\end{lemma}
\begin{proof}
    The positive structure consists of invertible elements and must contain 1 and be fixed under addition and $\sigma$-congruence. Thus $-1 = \sigma(i)1i$ is also positive. Then  $0 = -1+1$ must be in the positive structure and $0$ is not invertible.
\end{proof}

\begin{lemma}
   The positive real numbers are a positive structure for $\C$ with $\sigma$ is complex conjugation. 
\end{lemma}
\begin{proof}
   It is simple to check the positive reals are fixed by $\sigma$, contain 1 and are closed under addition and inversion. A quick calculation shows that $\sigma$-congruence by $a+b i$, corresponds to multiplication by $a^2+b^2$ which is positive as long as $a+b i$ is invertible.\\ 
\end{proof}

\begin{proposition}
    Let $z_1,z_2,z_3\in \C$ be three roots of unity whose product is 1. Then the triple $(z_1, z_2,z_3)$ is admissible for $\sigma$ complex conjugation.
\end{proposition}
\begin{proof}
    As the inverse of a root of unity is another root of unity, it is clear that $z_1,z_2,z_3$ are units in the subring of integral elements. We now compute the noncommutative angles:
    \[  \overline{z_1^{-1}}z_2\overline{z_3^{-1} } = \frac{1}{|z_1| \cdot |z_3|}z_1z_2z_3  \]
    By assumption $|z_1| = |z_2| = |z_3| = 1$ and $z_1z_2z_3 = 1 \in \R$, so every angle is 1 which is positive.
\end{proof}

We now see the mutation rule splits into separate mutation rules for the modulus and the argument. As such we write a triple in polar notation, $z_1 = Ae^{2\pi ia}, z_2 = Be^{2\pi ib},z_3 = Ce^{2\pi ic}$. The condition for the angles to be positive simplifies to the requirement that $a+b+c \in \Z$. We then compute mutation at $z_1$.
 \begin{align*}
     z'_1 =~& Ce^{2\pi i(-c)}\frac{1}{A}e^{2\pi i(-a)}Ce^{2\pi i(-c)} + Be^{2\pi ib}\frac{1}{A}e^{2\pi ia}Be^{2\pi ib}\\
     =~& \frac{C^2}{A}e^{2\pi i (-a-2c)} + \frac{B^2}{A}e^{2\pi i(a+2b)}\\
     =~& \frac{C^2}{A}e^{2\pi i (b-c)} + \frac{B^2}{A}e^{2\pi i(b-c)}\\
     =~& \frac{B^2+C^2}{A}e^{2\pi i(b-c)}
 \end{align*}
 We observe that the norms of the complex numbers follow the classical mutation rule, while the arguments evolve separately. 
 
\begin{figure}[hb]
    \centering
    \includegraphics[height=3in]{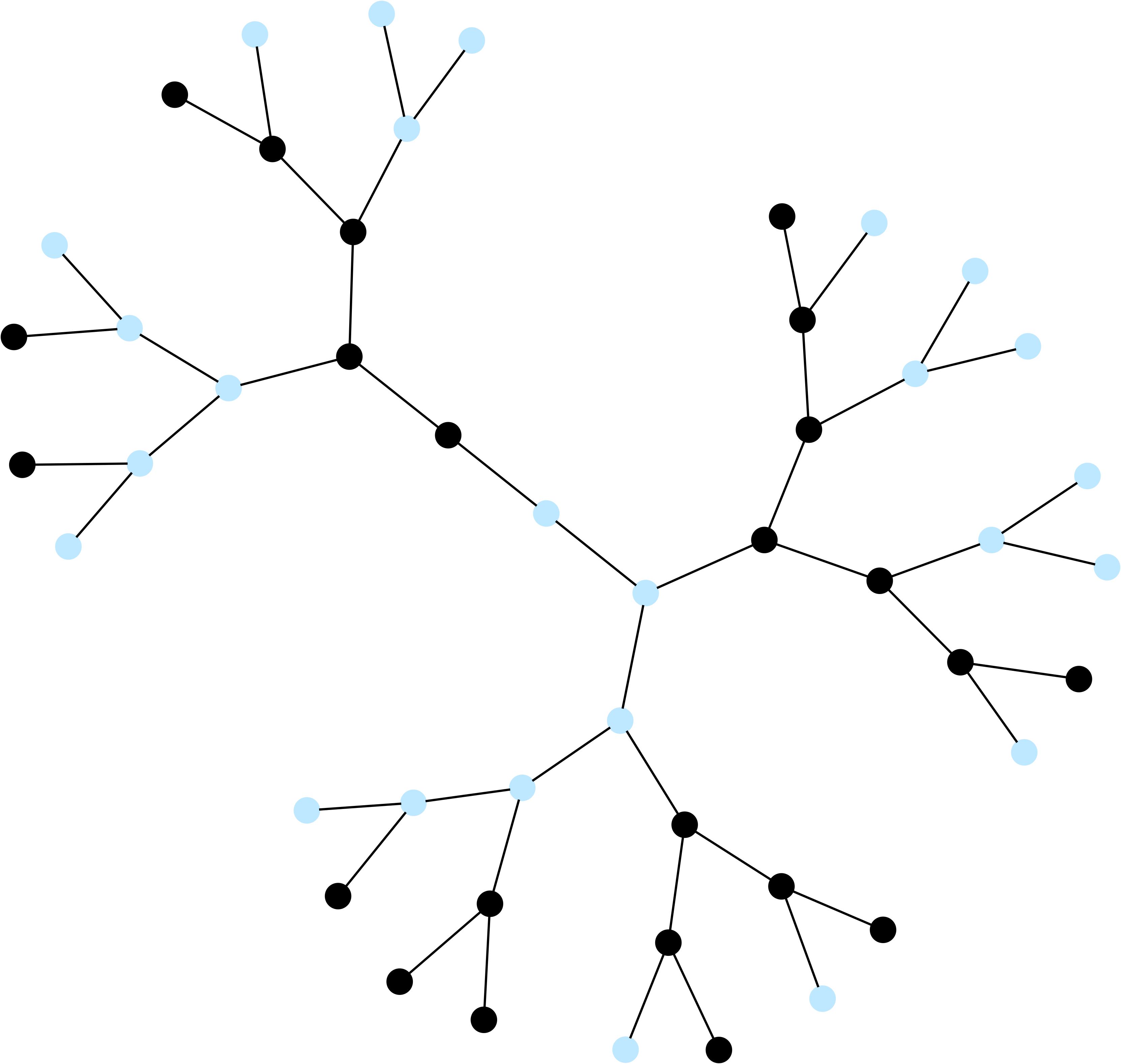}
    \caption{A complex Markov tree for fourth roots of unity. Black nodes and blue nodes distinguish the two possible patterns of arguments. The two central nodes of degree 2 correspond to $(1,i,-i)$ and $(-2,-i,-i)$ where mutation at either imaginary location produces the same triple.}
    \label{fig:ComplexMarkovTree}
\end{figure}

 \begin{remark}
     When the roots of unity are chosen to be nontrivial, this starting tuple gives an example of a distinct tree from the usual Markov numbers. However all the noncommutative angles are $1$ and thus the Markov function is $6$ in this case.
 \end{remark}
 Due to the splitting, we see that every triple obtained by mutation has norms corresponding to a standard Markov triple. However the arguments will be $2\pi\frac{a}{n}$ for some $0\leq a < n$. As in \Cref{rem:MarkovSymmetry} there are potentially (at least) six distinct angles associated to each Markov number. It is an interesting question to ask which ratios are associated with a given Markov number and what set of ratios appear due to mutation.
 As a concrete example consider the initial seed $(1,i,-i)$ consisting of fourth roots of unity. Its neighbors are the tuples  $(-2,-i,-i), (-2 i, 1,i), (-2i,1,i)$. In these tuples we see all possible arguments with denominator $4$. Furthermore we see two distinct patterns of arguments in a tuple: $(1,i,-i)$ and $(-1,-i,-i)$. In \Cref{fig:ComplexMarkovTree}, we color the nodes of the Markov tree black if they are of the first type and blue if they are of the second type to see the dynamics of the angle pattern. 
 \begin{conjecture}
      If the initial triple consists of 3 roots of unity $(z_1,z_2,z_3)$ such that $n$ is the minimal integer with $z_i^n = 1$ then every triple of $\frac{a}{n},\frac{b}{n},\frac{c}{n}$ that sums to $1$ is achieved as the argument of a complex Markov triple. Furthermore if $(z_1,z_2,z_3)$ are irrational angles then the image of the arguments of Markov triples is dense in the circle.
 \end{conjecture}
     Experimentally this conjecture holds up to $n=100$.

\subsection{Dual Markov Numbers}\label{sec:DualMarkov}
    Another interesting example is to consider the dual numbers $\{a + \alpha\epsilon | a,\alpha\in \R, \epsilon^2 = 0\}$. We focus on the subring $R$ of integral elements of the form $a + \alpha \epsilon$ with $a,\alpha \in \Z$. The units of this subring are elements of the form $1 + \alpha\epsilon$ for $\alpha \in \Z$. It will simplify computations to write invertible dual numbers in the form $a(1+\alpha\epsilon)$ with $a \neq 0$. This is analogous to the polar form of a complex number.
    As mentioned in \Cref{ex:involution} there are several natural choices of involution. \\  

   We start with case of trivial $\sigma$. Unlike the complex numbers there is a positive structure here.
   \begin{lemma}
       The set of dual numbers with positive real part is the positive structure with trivial involution.
   \end{lemma}
   \begin{proof}
       We compute the inverse and $\sigma$-congruence of an arbitrary invertible dual number $a(1+\alpha\epsilon)$ as these are the difficult conditions to check.
       \[(a(1 + \alpha\epsilon))^{-1} = \frac{1}{a}(1-\alpha \epsilon) \hspace{3pc} (b(1+\beta\epsilon))(a(1+\alpha\epsilon))(b(1+\beta\epsilon)) = b^2a(1+(\alpha+2\beta)\epsilon)\]
       We see in both cases if $a>0$ the results also has positive real part. 
   \end{proof}

   \begin{proposition}
       The triple $(1+\alpha \epsilon,1+\beta\epsilon,1+\gamma\epsilon)$ is admissible for any $\alpha,\beta,\gamma\in \R$.
   \end{proposition}
   \begin{proof}
       It is clear this triple consists of units in $R$, $(1+\alpha \epsilon)^{-1} = 1-\alpha \epsilon$. We compute a noncommutative angle:
       \[\sigma(1+\alpha \epsilon)^{-1}(1+\beta\epsilon)\sigma(1+\gamma\epsilon)^{-1} = (1 + (-\alpha+\beta-\gamma)\epsilon) \]
       The other angles are obtained via permuting $\alpha,\beta,\gamma$ and thus all have real part 1 and are positive.
   \end{proof}
   We compute the value of the Markov function for this triple:
   \[F(1+\alpha \epsilon,1+\beta\epsilon,1+\gamma\epsilon) = 6-2(\alpha+\beta+\gamma)\epsilon\]
   \begin{remark}
       If we take $\alpha+\beta+\gamma=0$ with at least one of $\alpha,\beta,\gamma$ nonzero we obtain another distinct tree of integer Markov-like numbers with the same value of Markov function as $(1,1,1)$.
   \end{remark}
   
\begin{remark}
    For arbitrary dual numbers $(a(1+\alpha\epsilon),b(1+\beta\epsilon),c(1+\gamma\epsilon))$, the noncommutative angles are cyclic rotates of $\frac{b}{ac}(1+(-\alpha+\beta-\gamma)\epsilon)$. As such a triple has all positive angles if an even number of $a,b,c$ are negative. Thus a triple of the form $(-1+\alpha\epsilon, -1+\beta\epsilon, 1+ \gamma\epsilon)$ is also admissible.
\end{remark}
    
   The mutation rule also splits into a real and shadow part. We compute
   \[ (a+\alpha \epsilon, b + \beta \epsilon, c+ \gamma \epsilon) \mapsto \left(\frac{b^2+c^2}{a} + \frac{2ab\beta + 2ac\gamma - (b^2 + c^2)\alpha}{a^2} ,c+\gamma \epsilon, b+ \beta \epsilon\right)\]
   Here we see the real part mutates like the classic Markov rule. The shadow rule is more complex and depends linearly on $\alpha,\beta,\gamma$. While the initial vector of shadows $(\alpha,\beta,\gamma)$ has no effect on admissibility, we can try to identify vectors with additional properties. \\
   This case is studied in \cite{Bonin-DualMarkovNumbers} without the language of noncommutative cluster algebras. They observe that choosing $\alpha=\beta=\gamma$ is not particularly interesting as this results in the usual Markov numbers multiplied by $(1+\alpha\epsilon)$. Thus the most interesting triples come from choosing $(\alpha,\beta,\gamma)$ as a linearly independent vector from $(1,1,1)$. In \cite{Bonin-DualMarkovNumbers} they examine the basis $(0,1,1)$ and $(0,1,0)$. They observe that in one third of the tree for $(1, 1+\epsilon, 1+\epsilon)$ the shadow part is also positive. However the shadows in the rest of the tree may be negative. For example the terms with real part 169 in the tree are $169 - 411 \epsilon, 169 - 172 \epsilon$ and $169 + 921\epsilon$. This makes sense as we only consider the positive structure with real part positive. In fact the set of dual numbers with positive real and shadow part is not preserved by $\sigma$-congruence and so is not a possible choice of positive structure.\\
   %The case with $\sigma$ trivial is studied in \cite{Bonin-DualMarkovNumbers} without the language of noncommutative cluster algebras. There they claim the natural nontrivial starting triple is $(1,1+\epsilon,1+\epsilon)$. We recover this choice in the language of noncommutative cluster algebras by considering the requirement that the noncommutative angles are positive. Here the noncommutative angles are the cyclic rotates of:
  % \[\sigma(1+\alpha \epsilon)^{-1}(1+\beta\epsilon)\sigma(1+\gamma\epsilon)^{-1} = 1 + (-\alpha+\beta-\gamma)\epsilon \]
  % The condition that all these angles are positive corresponds the cone $\alpha\geq\beta+\gamma,\beta \geq \alpha + \gamma,\gamma\geq\alpha+\beta $. The extremal rays of the cone are $(0,1,1)$ and cyclic rotates, which corresponding to the starting choice of \cite{Bonin-DualMarkovNumbers}.\\
  % We can also compute the Markov function in this case is $2(3-(\alpha+\beta+\gamma)\epsilon)$. 
  % \[1 + (-\alpha + \beta\]
  % \begin{remark}
  %     This starting choice gives a distinct tree from the standard Markov numbers as the coefficient of $\epsilon$ is nontrivial here.
  % \end{remark}

   If $\sigma$ is taken to be the conjugation map $a+\alpha\epsilon \mapsto a -\alpha \epsilon$ we obtain other interesting dual Markov triples. 
   In this case, the $\sigma$ fixed points are isomorphic to the real numbers and so the positive reals are a natural choice of positive structure. This cone is fixed by $\sigma$-congruence as $\sigma(a+\alpha\epsilon)(a+\alpha \epsilon) = a^2$.
\begin{proposition}
    The triple $(1+\alpha\epsilon,1+\beta\epsilon,1+\gamma\epsilon)$ with $\alpha + \beta + \epsilon = 0$ is admissible for $\sigma$ ``dual conjugation''.
\end{proposition}
   \begin{proof}
       The interesting condition to check is that the noncommutative angles are positive. We compute the angle for an arbitrary triple.
       \[\sigma(a(1+\alpha \epsilon))^{-1}b(1+\beta\epsilon)\sigma(c(1+\gamma\epsilon))^{-1} = \frac{b}{ac}(1 + (\alpha+\beta+\gamma)\epsilon) \]
       Here we observe the shadow part is more symmetric than the trivial involution case. Moreover by assumption this simplifies to $1$ which is clearly in the positive structure.
   \end{proof}
\begin{remark}
    This starting triple has Markov function $6$. So when one of $\alpha,\beta,\gamma$ is nonzero, we obtain another example of a distinct tree with the same value of the Markov function as $(1,1,1)$.
\end{remark}
   The mutation rule also splits more cleanly into real and shadow parts in this case. Using the fact that $(\alpha+\beta+\gamma)= 0$ we have 
   \[ (a(1+\alpha \epsilon), b(1 + \beta \epsilon), c(1+ \gamma \epsilon)) \mapsto \left(\frac{b^2+c^2}{a}\left(1+(\beta-\gamma)\epsilon\right) ,c(1+\gamma\epsilon), b(1-\beta\epsilon)\right) \]
   This mirrors the complex case where the shadows $\alpha$ mutate like the argument of the complex numbers and the norm/real part mutates like classic Markov numbers.

\subsection{Matrix Markov Numbers}\label{sec:MatrixMarkov}
Here we consider $\A = \Mat(n,\R)$ with $\sigma$ the standard matrix transpose. We chose $\Mat_n(\Z)$ as the subring of integral elements to focus on. The positive structure consists of the set of positive definite symmetric matrices.\\

Here there are many complicated examples of admissible triples, even in $n=2$. As such we provide examples some we find particularly interesting. 
 \begin{example}
Let $B_t = \begin{bmatrix}
     1 & t\\ 0 & 1
 \end{bmatrix}$ and consider the tuple $(1,B_t,B_t^{-1})$. 
The value of the Markov function is 
\[\begin{bmatrix}
    6+2t^2 & 0\\
    0 & 6 + 2t^2
\end{bmatrix}\]
 \end{example}
 Already in this example the new Markov numbers are complicated. For instance the matrices obtained by freezing the identity and mutating alternately at the other two edges are:
 \begin{equation*}
 \begin{bmatrix} 1 & t \\ 0 & 1\end{bmatrix}\leftrightarrow
     \begin{bmatrix}t^2+2&-t\\ -t^3-3t&t^2+2\end{bmatrix}\leftrightarrow
\begin{bmatrix}t^4+5t^2+5&-t^5-6t^3-8t\\ -t^3-3t&t^4+5t^2+5\end{bmatrix} \leftrightarrow
\begin{bmatrix}t^6+8t^4+19t^2+13&-t^5-6t^3-8t\\ -t^7-9t^5-25t^3-21t&t^6+8t^4+19t^2+13\end{bmatrix}
%\begin{bmatrix}t^8+11*t^6+42*t^4+65*t^2+34&-t^9-12*t^7-51*t^5-90*t^3-55*t\\ -t^7-9*t^5-25*t^3-21*t&t^8+11*t^6+42*t^4+65*t^2+34\end{bmatrix}
 \end{equation*}
 This is a particularly nice path in the tree as it results in matrices with positive monic polynomials as diagonal entries. However more complicated paths result in more complex matrices. For example the new matrices along the path given by mutating at $1$ then $B_t$ then $B_t^{-1}$ are
 \begin{equation*}
     \begin{bmatrix}2 & 2t\\ -2t & 2\end{bmatrix} \hspace{2pc}
     \begin{bmatrix}5 & -4t\\ 4t^3+11t & 5 \end{bmatrix} \hspace{2pc}
     \begin{bmatrix}-16t^4-32t^2+29 & -16t^3-52t\\ 36t^3+93t & -16t^4-32t^2+29 \end{bmatrix}
 \end{equation*}
 We emphasize that although these matrices have negative entries the angles at every point are positive definite symmetric matrices.
 
 To get a better handle on these matrices we can compute the determinant or trace. 
 \begin{lemma}
     The trace of each ``matrix Markov number'' in this case is a polynomial $p(t)$ with $p(0) = 2m$ for $m$ a standard Markov number. Similarly the determinant is a polynomial $q(t)$ with $q(0) = m^2$.
 \end{lemma}
 \begin{proof}
     At $t=0$ these matrices must match the case $(1,1,1)$ for which each matrix is $m\cdot 1$ for $m$  Markov number. Then clearly $\trace (m\cdot 1) = 2m$ and $\det(m\cdot 1) = m^2$ as needed. 
 \end{proof}
 \begin{remark}
     As each matrix Markov number $m$ is invertible for all $t$ the determinant has fixed sign for all $t$. Furthermore from these previous lemma we see that at 0 the determinant is a classic Markov number squared. Thus the determinant is always positive. Experimentally we have observed the stronger fact that each determinant is an even function with positive integer coefficients.  
 \end{remark}
 \begin{remark}
     The same classic Markov number appears in matrices with distinct trace and determinant as seen by 5 appearing as the constant term of the diagonal in both of the following matrices obtained from $(1,B_t,B_t^{-1})$ by mutation.
     \[\begin{bmatrix}t^4+5t^2+5&-t^5-6t^3-8t\\ -t^3-3t&t^4+5t^2+5\end{bmatrix} \hspace{2pc} \begin{bmatrix}5 & -4t\\ 4t^3+11t & 5 \end{bmatrix}\]
 \end{remark}
 \begin{remark}
     For $t \neq 0$ these matrices produce a tree distinct from the standard tree starting from $(1,1,1)$. This is clear as every matrix in the standard tree is diagonal.
 \end{remark}

 \begin{remark}
     This example can be generalized by taking $1$ and $t$ to be members of larger rings. One interesting such choice is the $n\times n$ identity matrix and $\Diag{1}{t_1,\dots,t_n}$ to obtain example Markov numbers in $\GL(2n,\R)$. In this case, the matrix Markov numbers can be separated into $n$ copies of the $2\times 2$ example above.
 \end{remark}

If we consider a slightly larger subring of integral elements, $\Mat_n(\Z[t,t^{-1}])$ we have another interesting example. 
 \begin{example}
     Let $B_t = \begin{bmatrix}
         t & 0\\
         0 & t^{-1}
     \end{bmatrix}$ and consider the triple $(1,B_t,B_t^{-1})$. Note that in this example we recover the Markov numbers at $t=1$. The value of the invariant function is 
     \[\begin{bmatrix}
         \frac{2}{t^2}+2+2t^2 & 0\\
         0 & \frac{2}{t^2} + 2 + 2t^2
     \end{bmatrix}\]
    
 \end{example}
  Experimentally the traces along the path where we freeze $A = 1$  are \[\sum_{k=1}^{2n-1} \frac{c_k}{t^{2k-1}}+c_kt^{2k-1}\] 
     where $c_k$ is the number of $01$ sequences of length $2n-1$ with exactly $k$ ones and all odd runlengths (\href{https://oeis.org/A078821}{oeis:A078821}). When $t=1$ this sum is twice an odd Fibonnaci number as expected.\\

For larger $n$ we can define families with more parameters, leading to even greater flexibility and complexity.
\begin{example}
    Let $\Diag{1}{t_1,\cdots,t_{n-1}}$ be the $n\times n$ matrix with super diagonal $t_1,\cdots,t_{n-1}$. Then let $B_{\vec{t}} = \exp(\Diag{1}{t_1,\cdots,t_n})$ and consider $(1,B_{\vec{t}},B_{\vec{t}}^{-1})$. As before when all $t_i =0$ we recover the standard Markov numbers. However now the value of the Markov function is more complicated. We give the values for $n=3,4$ below
    \[\begin{bmatrix}
        6+2t_1^2+\frac{1}{2}t_1^2 t_2^2 & 0 & 2t_1 t_2\\
        0 & 6+ 2t_1^2+2t_2^2 & 0\\
        2t_1 t_2 & 0 & 6 + 2t_2^2 + \frac{1}{2}t_1^2t_2^2
    \end{bmatrix}\]
    \[\begin{bmatrix}
        6+2t_1^2(1+\frac{1}{4}t_2^2+\frac{1}{36}t_2^2t_3^2) & 0 & 2t_1 t_2(1 + \frac{1}{6}t_3^2)& 0\\
        0 & 6+ 2t_1^2+2t_2^2(1+ \frac{1}{4}t_3^2) & 0  & 2t_2t_3(1+\frac{1}{6}t_1^2)\\
        2t_1t_2(1+ \frac{1}{6}t_3^2) & 0 & 6 + 2t_2^2(1+\frac{1}{4}t_1^2) + 2t_3^2 & 0\\
        0 & 2t_2 t_3(1 + \frac{1}{6}t_1^2) & 0 & 6 + 2t_3^2( 1 + \frac{1}{4}t_2^2 + \frac{1}{36}t_2^2t_3^2)
    \end{bmatrix}\]
\end{example}

\begin{remark}
    Via the standard embedding of the complex numbers in $\Mat_2(\R)$, we can via the examples from \Cref{sec:ComplexMarkov} as examples of matrix Markov numbers as well.
\end{remark}

\subsection{Group Ring Markov Numbers}\label{sec:GroupRingMarkov}
Our first choice is of the choice of $\Z[G]$ for the subring of integer elements of $\R[G]$. Identifying the units in this ring is difficult for general $G$. However for $G$ finite, Hurley \cite{hurley-GroupRings} gives an embedding of  $\R[G]$ in $\Mat_{|G|}(\R)$ via the formula
\[\sum_{g\in G} a_g g \mapsto \begin{bmatrix}
    a_{g_1^{-1}g_1} & a_{g_1^{-1}g_2} & a_{g_1^{-1}g_3}& \dots & a_{g_1^{-1}g_n}\\
    a_{g_2^{-1}g_1} & a_{g_2^{-1}g_2} & a_{g_2^{-1}g_3}& \dots & a_{g_2^{-1}g_n}\\
    \vdots & \vdots & \vdots & \vdots & \vdots\\
    a_{g_n^{-1}g_1} & a_{g_1^{-1}g_2} & a_{g_1^{-1}g_3}& \dots & a_{g_n^{-1}g_n}\\
\end{bmatrix}. \]
\begin{lemma}
    The map above is an involutive ring homomorphism where $\sigma(\sum a_g g) = \sum a_g g^{-1}$.
\end{lemma}
\begin{proof}
    In \cite{hurley-GroupRings}, they show the map is a ring homomorphism. From the structure of the image we see that transpose exchanges $a_{g_i^{-1}g_j}$ with $a_{g_{j^{-1}}g_i}$. As $g_i^{-1}g_j$ is the inverse of $g_j^{-1}g_i$ we see this maps intertwines the involutions as well. 
\end{proof}
Since the matrix embedding respects the involution, we can pullback the positive structure on matrices to define a positive structure on $\R[G]$. Concretely an element $\sum a_{g} g$ is in the positive cone if the corresponding matrix is symmetric positive definite. Furthermore the choice of $\Z[G]$ as the integer subring corresponds to the choice of matrices with integer entries. 
\begin{proposition}
    Let $g \in G$. The tuple $(g^a,g^b,g^c)$ with $g^{a+b+c} =1$ is an admissible triple.
\end{proposition}
\begin{proof}
    The element $\pm h$ for any $h\in G$ is a unit of $\Z[G]$. Thus it remains to check the noncommutative angles are positive. We see that $\sigma(g^a) = g^{-a} = (g^{a})^{-1}$ and every noncommutative angle is $g^{a+b+c}$. By assumption this is 1 and so is positive as needed.
\end{proof}
The admissible triples above are equivalent to the primitive roots of unity discussed in \Cref{sec:ComplexMarkov}.

For small cyclic groups these are the only examples of admissible triples, as the only units of $\Z[G]$ are of the form $g^a$. However in $G=\Z/5\Z$ the are elements of the integer subring with more than one nonzero term. Let $Z$ be the generator of $G$. The first such nontrivial example has at least 3 nonzero terms. The following table shows all units of $\Z/5\Z$ with 3 nonzero terms up to multiplication by $\pm 1$. The bottom row of the table is the inverse of the top and elements related by $\sigma$ are grouped together.
\begin{center}
\begin{tabular}{c | c c |c c}
    $-1 + Z^2 + Z^3$ & $1 - Z + Z^2$ & $1 + Z^3 - Z^4$ & $1 + Z - Z^3$ &  $1 - Z^2 + Z^4$  \\ 
    $-1 + Z + Z^4$ & $Z + Z^2 - Z^4$ & $-Z + Z^3 + Z^4$ & $Z - Z^2 + Z^3$  & $Z^2 - Z^3 + Z^4$  
\end{tabular} 
\end{center}
There are of course other units with more nonzero terms, for example $3 + Z - 2 Z^2 - 2 Z^3 + Z^4 = (-1+Z^2+Z^3)^2$. However the list above is enough to produce many interesting admissible triples. In particular there are triples with nontrivially distinct elements, i.e. no element is a power or $\sigma$ of a power of any other element. 
\begin{example}
    The triple $(Z, 1 - Z + Z^2, 1 + Z - Z^3)$ is an admissible triple whose elements are nontrivially distinct. 
\end{example}
To see this it suffices to compute the angles for the triple above are:
    \[1,3-2(Z+Z^4)+(Z^2+Z^3),3+(Z+Z^4)-2(Z^2+Z^3)\]
    Note that although these have negative coefficients the corresponding symmetric matrices are positive definite. The corresponding value of Markov function is $14 - 2 Z - 2Z^2 - 2Z^3 - 2Z^4 = 2(3+4-Z-Z^2-Z^3-Z^4)$. 
    
    \begin{remark}
        We can view this admissible triple as nontrivial deformations of the Markov numbers. 
    \end{remark}
    First consider these expressions as polynomials in $Z$. If we evaluate the initial triple at $Z=1$ we obtain $(1,1,1)$. Thus every group ring element evaluated at $Z=1$ will be a classic Markov number. As in \Cref{rem:MarkovSymmetry} in the full Markov tree we will see (at least) six deformations of each Markov number greater than 2. We give the deformations of $2$ and $5$ we observed below.
\begin{center}
\small
    \begin{tabular}{c|c c c }
    2 & $7 + 7 Z - 2 Z^2 - 8 Z^3 - 2 Z^4$ & $-1 - Z + 6 Z^2 - Z^3 - Z^4$ &$-8 +7 Z - 2 Z^2 - 2 Z^3 + 7 Z^4$\\
    \hline
    5 & $-18 + 53 Z + 53 Z^2 - 18 Z^3 - 65 Z^4$ & $46 - 11 Z - 11 Z^2 + 46 Z^3 - 65 Z^4$ & $53 - 18 Z - 18 Z^2 + 53 Z^3 - 65 Z^4$\\
      & $-65 - 11 Z + 46 Z^2 + 46 Z^3 - 11 Z^4$ & $173 + 54 Z - 138 Z^2 - 138 Z^3 + 54 Z^4$ &  $-138 + 173 Z - 138 Z^2 + 54 Z^3 + 54 Z^4$ \\
\end{tabular}
\end{center}

%%%%%%%%%%%%

\bibliography{bibliography.bib}
\bibliographystyle{alpha}
\end{document}